\documentclass[sn-mathphys-num]{sn-jnl}

\usepackage{color}
\usepackage{latexsym}
\usepackage{graphicx}
\usepackage{caption}
\usepackage{subcaption}
\usepackage{enumitem}
\usepackage{bbm}
\usepackage{mathrsfs}
\usepackage{amsmath, amssymb}  
\usepackage{mathtools}         
\usepackage{booktabs}          
\usepackage{amsfonts}

\usepackage{amsmath}
\numberwithin{equation}{section}
\usepackage{amsthm}
\usepackage{amssymb}

\definecolor{darkblue}{rgb}{0.0,0,0.5}
\usepackage{algorithm,algorithmic}
\usepackage{caption}

\newcommand{\R}{\mathbb{R}}
\newcommand{\E}{\mathbb{E}}

\renewcommand{\H}{\mathcal{H}}

\newcommand{\kernel}{K}

\newcommand{\norm}[1]{\left\|{#1}\right\|} 

\providecommand{\minimize}{\mathop{\rm minimize}}

\newcommand{\matrixnorm}[1]{\left|\!\left|\!\left|{#1}
  \right|\!\right|\!\right|} 

\newcommand{\normal}{\mathsf{N}}  
\newcommand{\N}{\mathbb{N}}
\renewcommand{\P}{\mathbb{P}}

\newtheorem{lemma}{Lemma}

\newtheorem{assumption}{Assumption}[section]
\newtheorem{innercustomgeneric}{\customgenericname}
\providecommand{\customgenericname}{}

\renewcommand{\H}{\mathbb{H}}

\newcommand{\Cov}{{\rm Cov}}

\renewcommand{\H}{\mathcal{H}}

\newcommand{\eps}{\epsilon}

\newcommand{\Range}{{\rm Range}}

\newcommand{\grad}{\nabla}

\newcommand{\rank}{{\rm rank}}

\newcommand{\T}{\mathcal{T}}

\newcommand{\profile}{\psi}

\newcommand{\C}{\mathbb{C}}

\newcommand{\Sym}{{\rm \mathrm{Sym}}}

\providecommand{\minimize}{\mathop{\rm minimize}}

\theoremstyle{thmstyleone}%
\newtheorem{theorem}{Theorem}
\theoremstyle{thmstyletwo}%
\newtheorem{example}{Example}%
\newtheorem{remark}{Remark}%

\theoremstyle{thmstylethree}%
\newtheorem{definition}{Definition}%

\begin{document}

\title[A Theory of Feature Learning in Kernel Models]{A Theory of Feature Learning in Kernel Models}


\author*{\fnm{Yunlu} \sur{Chen}}\email{yunluchen2026@u.northwestern.edu}
\author*{\fnm{Yang} \sur{Li}}\email{yl454@cam.ac.uk}
\author*{\fnm{Keli} \sur{Liu}}\email{keli.liu25@gmail.com}
\author*{\fnm{Feng} \sur{Ruan}}\email{fengruan@northwestern.edu}

\affil*{{\normalsize \orgdiv{Department of Statistics and Data Science}, \orgname{Northwestern University}}}
\affil*{{\normalsize \orgdiv{Department of Pure Mathematics and Mathematical Statistics}, \orgname{Cambridge University}}}
\affil*{{\normalsize \orgname{Company E}}}
\affil*{{\normalsize \orgdiv{Department of Statistics and Data Science}, \orgname{Northwestern University}}}



\abstract{
\phantom{.....}
We study feature learning in a compositional variant of kernel ridge regression in which the predictor is applied to a learnable linear transformation of the input. When the response depends on the input only through a low-dimensional predictive subspace, we show that all global minimizers of the population objective for the linear transformation annihilate directions orthogonal to this subspace, and in certain regimes, exactly identify the subspace. Moreover, we show that global minimizers of the finite-sample objective inherit the exact same low-dimensional structure with high probability, even without any explicit penalization on the linear transformation. 
}

\maketitle


\section{Introduction} 
\indent\indent
Modern learning systems, including neural networks, are empirically
observed to form internal representations that are 
\emph{low-dimensional} and \emph{predictive}. Rather than relying on manually
designed features, these models appear to automatically extract low-dimensional
representations that capture task-relevant structure in the data. This
phenomenon is often referred to as \emph{feature learning}~\cite{BengioCoVi13}.  

Feature learning is commonly associated with two key ingredients.  First, the model is
\emph{compositional}, meaning that predictions are formed as compositions of mappings of
the form $x \mapsto g(h(x))$, where the intermediate representation $h(x)$ is not fixed
\emph{a priori}. Second, these representations are learned by \emph{optimizing}
a data-driven loss function tied to predictive performance. While this combination is ubiquitous 
in practice, the mechanisms by which compositional structure and optimization give rise 
to predictive low-dimensional representations remain incompletely understood.

%


In this work, we study feature learning in a compositional kernel model obtained by a \emph{minimal} extension of classical kernel ridge regression (KRR)~\cite{CuckerSm02}, in a setting that allows structural mechanisms underlying feature learning to be analyzed in a mathematically precise way. Classical kernel ridge regression provides a canonical framework for supervised learning, in which prediction is formulated as the minimization of a single variational objective over a function space. Its simplicity makes it a natural \emph{baseline} for \emph{isolating} structural effects introduced by representation learning.

As a point of reference, we first recall the classical KRR formulation~\cite{CuckerSm02}. 
Let $(X,Y) \sim \mathbb P$ with $X \in \mathbb R^d$ and
$Y \in \mathbb R$, and let $H$ be a reproducing kernel Hilbert space (e.g., Sobolev space).
Classical KRR learns a predictor $f \in H$ by minimizing
\begin{equation*}
	\min_f  \E[(Y - f(X))^2] + \lambda \norm{f}_H^2.
\end{equation*}
The first term measures approximation error, while the second term
$\lambda \|f\|_H^2$ controls the regularity of $f$
through the RKHS norm, which helps generalization to unseen data. In classical KRR, the 
representation is fixed. 

To incorporate representation learning, we introduce a linear transformation of the input and 
consider predictors of the form $x \mapsto f(Ux)$, where $U \in \R^{d \times d}$ is 
learned jointly with $f \in H$. This leads to the variational problem
\begin{equation}
\label{eqn:variational-objective-f}
		\min_U \min_f \E[(Y - f(U X))^2] + \lambda \norm{f}_H^2.
\end{equation}
Relative to classical KRR, introducing the linear transformation $U$ allows the 
input representation to be adapted during training, so the predictor $f$ acts on 
$UX$ rather than directly on $X$. In this sense, the matrix $U$ serves as a 
representation-learning component of the compositional model $x \mapsto f(Ux)$.  

In practice, however, the joint minimization in~\eqref{eqn:variational-objective-f}
is carried out using finitely many samples. Given i.i.d.\ observations
$\{(X_i,Y_i)\}_{i=1}^n \sim \mathbb P$, the population expectation is replaced by
its empirical counterpart, leading to the finite-sample objective
\begin{equation}
\label{eqn:variational-objective-f-empirical}
	\min_U \min_f \E_n[(Y - f(U X))^2]  + \lambda \|f\|_H^2,
\end{equation}
where $\E_n$ denotes the empirical expectation: $\E_n [F(X, Y)] := \frac{1}{n} \sum_{i=1}^n F(X_i, Y_i)$.

This paper studies feature learning in the compositional kernel model~\eqref{eqn:variational-objective-f} 
and its finite sample counterparts~\eqref{eqn:variational-objective-f-empirical},
understood as the emergence of predictors that depend on fewer effective 
degrees of freedom than the ambient input. We formalize this notion by assuming that $Y$
 depends on $X$ only through a low-dimensional projection:
there exists a subspace $S \subset \R^d$ such that
\begin{equation}
\label{eqn:subspace-low-dimension}
	\mathbb{E}[Y | X] = \E[Y| \Pi_S X]
\end{equation} 
where $\Pi_S$ denotes the orthogonal projection onto $S$. 

Identifiability of the predictive subspace $S$ is ensured by the notion of the 
\emph{central mean subspace} in classical statistics, defined as the smallest subspace satisfying this 
conditional independence (Definition~\ref{definition:central-mean-subspace})~\cite{CookLi02}. Under this identification, feature learning 
corresponds to the property that the learned transformation $U$ preserves 
directions along $S$ while suppressing directions orthogonal to it.

\vspace{.5em} 
\begin{remark} 
\emph{
The compositional formulation in
\eqref{eqn:variational-objective-f} admits a formal analogy to two-layer neural networks, 
in the sense that the map $x \mapsto Ux$ may  be interpreted as a linear first layer, and the function $f$ acts 
as a second-layer predictor. This observation is not used in the analysis and is mentioned 
only for \emph{intuition}. 
Further discussion of this analogy is deferred to Section~\ref{sec:analogy-to-two-layer-NN}.
}
\end{remark}

\subsection{A Motivating Numerical Experiment} 
\label{sec:puzzling-numerical-experiment}
\indent\indent
Before turning to the mathematical results, we describe a simple numerical 
observation that motivates the theoretical study. The phenomenon we report 
arises in a setting that is deliberately simple and does not incorporate any 
explicit mechanism designed to promote low-dimensional structure.

\begin{figure}[!htb]
\centering
  \begin{subfigure}{\linewidth}
    \centering
    \includegraphics[width=0.45\linewidth]{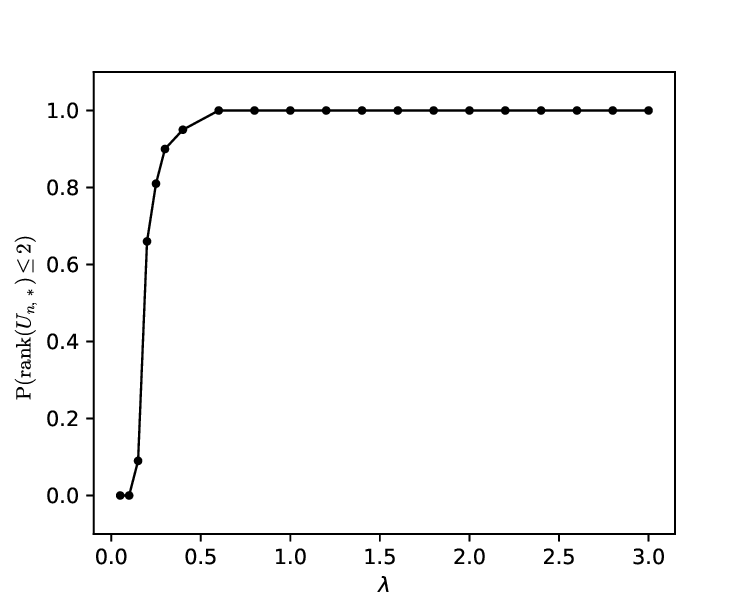}
    \includegraphics[width=0.45\linewidth]{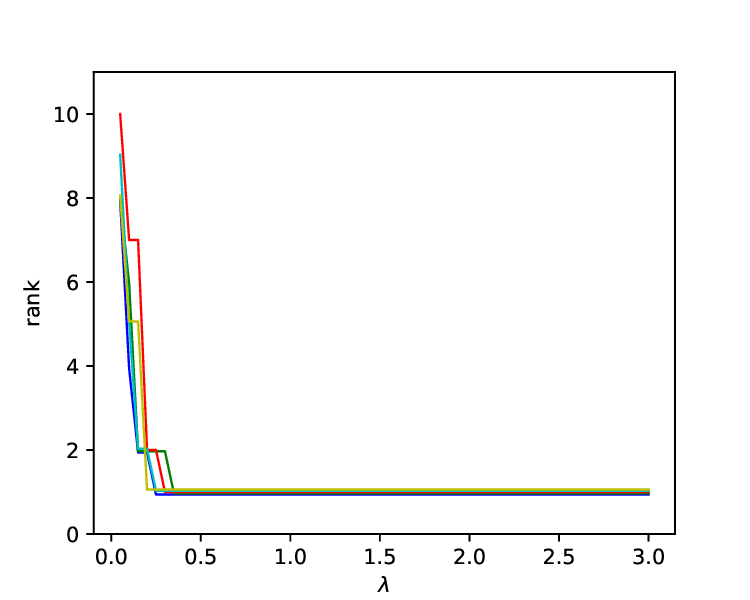}
  \end{subfigure}
      \caption{
The covariate $X\sim \normal(0,I)$, and the response $Y=F(X)+\epsilon$, where 
$F(x)=0.1(x_1 + x_2 + x_3)^3 + \tanh(x_1 + x_3 + x_5)$ and 
$\epsilon\sim N(0,\sigma^2)$ with $\sigma=0.1$. Here $n = 300$ and $d = 50$. 
On the left panel, we plot the empirical probability 
of $\rank(U_{n, *}) \le 2$ against $\lambda$. On the right panel, we display 
how the rank of the solution $U_{n, *}$ changes with different $\lambda$ values, 
using $5$ examples of $(X, y)$.}
\label{plot:intro} 
\end{figure}

In the experiment shown in Figure~\ref{plot:intro}, we generate $n = 300$ samples
with covariates $X \sim \normal(0, I)$ in $\mathbb{R}^{50}$. The response $Y$
depends on $X$ through a combination of nonlinear functions applied to two
linear projections, with additive Gaussian noise. The reproducing kernel Hilbert space $H$
is chosen to be the Gaussian RKHS. 
We numerically minimize the finite sample objective~\eqref{eqn:variational-objective-f-empirical}
and record the resulting solution matrix
$U_{n, *}$ across repeated trials and a range of values of the ridge parameter
$\lambda$. 

Across repeated experiments, the learned matrix $U_{n, *}$
consistently exhibits a pronounced tendency toward low rank. In particular, for
many $\lambda$, we observe that $\mathrm{rank}(U_{n, *})
\leq 2$ with high empirical frequency (Figure~\ref{plot:intro}).

At first glance, this observation is difficult to reconcile with standard
intuition. The standard intuition is 
that recovering a distribution's intrinsic low-dimensional structure from 
finite samples often requires regularization mechanisms such 
as low-rank penalties or early stopping. In the present setting,
no such mechanism is present, yet the learned transformation $U$ is \emph{exactly} low-dimensional.

This raises a natural question: is the observed low-rank behavior a finite-sample or algorithmic artifact, or does it reflect a more fundamental property of the underlying population objective~\eqref{eqn:variational-objective-f}?

As we show in the sequel of the paper, the latter is the case. Our main results identify geometric properties 
of the population objective~\eqref{eqn:variational-objective-f} that explain how low-dimensional 
structure is encoded in the variational formulation and inherited by minimizers of the finite-sample 
problem~\eqref{eqn:variational-objective-f-empirical}.

\subsection{Main Results} 
\indent\indent
We are interested in whether, and how the latent subspace structure
\eqref{eqn:subspace-low-dimension} is encoded in, and recoverable from, its variational
formulations \eqref{eqn:variational-objective-f} and
\eqref{eqn:variational-objective-f-empirical}, with an eye towards explaining the feature learning behavior observed in the numerical experiments of Section~\ref{sec:puzzling-numerical-experiment}. Specifically, we ask 

\vspace{.5em} 
\begin{itemize}
\item[($\mathsf{P1}$)] Does the population objective admit minimizers, whose associated
  transformation $U$ is low-dimensional, in the sense that it
 annihilates directions orthogonal to the true subspace $S$? More stringently, when 
  does the variational formulation \emph{exactly} identify the subspace $S$? 
\vspace{.2em} 
\item[($\mathsf{P2}$)] Can the predictive subspace $S$---including its dimension---be recovered from
finite samples by solving the associated variational optimization?
\end{itemize} 

\vspace{.5em} 
Answering these questions requires a careful analysis of the population objective, particularly its 
behavior at and near its minimizers. 
Our main results provide affirmative answers to both questions under natural assumptions on the data 
distribution and the RKHS $H$, which we state below. 
These conclusions rely on structural properties of the RKHS $H$ and do not extend to linear RKHSs; see Appendix~\ref{sec:linear-vs-nonlinear} for a discussion of the linear case.


\subsubsection{Assumptions} 
\indent\indent 
Throughout this paper, we assume the basic moment conditions
$\E[\norm{X}^8] < \infty$ and $\E[Y^2] < \infty.$ 

\vspace{.5em} 
\noindent\noindent 
\emph{$\mathsf{(A1)}$ Predictive Subspace and Identifiability}. \vspace{.3em} 

We formalize the 
dependence of $Y$ on a low-dimensional projection of $X$, as
suggested by~\eqref{eqn:subspace-low-dimension}, using the notion of 
\emph{central mean subspace} from statistics literature---a canonical and minimal object 
capturing all directions of $X$ relevant for predicting $Y$~\cite[Definition 2]{CookLi02}.

\vspace{.5em}

\begin{definition}[Central Mean Subspace]
\label{definition:central-mean-subspace}
Define \vspace{-.5em}
\[
S := \bigcap \left\{\, T \;\middle|\; 
T \text{ is a linear subspace of $\R^d$ and } \E[Y | X] = \E[Y | \Pi_T X] \,\right\}. \vspace{-.5em}
\]
If $S$ itself satisfies $\E[Y|X] = \E[Y|\Pi_S X]$, it is called the central mean subspace.
\end{definition}

\vspace{.5em}
The central mean subspace need not exist in complete generality: although multiple
subspaces may satisfy the defining conditional independence, their intersection
need not inherit this property. To rule out this degeneracy, we impose a mild
regularity condition on the input distribution.

\vspace{.5em}
\begin{assumption}
\label{assumption:X-continuous}
The distribution of $X$ is continuous.
\end{assumption}

\vspace{.5em}
Assumption~\ref{assumption:X-continuous} guarantees the existence of the central
mean subspace $S$~\cite{CookLi02}. It also ensures that the variational objective
\eqref{eqn:variational-objective-f} admits finite minimizers in $U$ (see 
Theorem~\ref{theorem:existence-of-global-minimizers}); in the absence
of such regularity, one can construct examples where  minimizers ``escape to infinity"
(see Section~\ref{sec:discrete-counterexample}).

\vspace{.8em}
\noindent\noindent
\emph{$\mathsf{(A2)}$ Non-Degeneracy of the Regression Problem}. 

\vspace{.3em}
\begin{assumption}
\label{assumption:non-degeneracy}
The conditional expectation $\E[Y|X] \neq 0$ with positive probability.
\end{assumption}
 
\vspace{.5em}
Without Assumption~\ref{assumption:non-degeneracy} (i.e., $\E[Y|X] \equiv 0$), the inner minimization over $f$ in
\eqref{eqn:variational-objective-f} is attained at $f \equiv 0$ for every $U$,
so the resulting objective is constant  in $U$ and feature learning is
vacuous.

\vspace{.8em}
\noindent\noindent
\emph{$\mathsf{(A3)}$ Noise Directions}.  

\vspace{.3em}
Let $S^\perp$ denote the orthogonal complement of $S$ in $\R^d$. 
Accordingly, the input $X$ admits 
the decomposition $X = \Pi_S X + \Pi_{S^\perp} X$.

\vspace{.5em} 
\begin{assumption}
\label{assumption:independence}
The random vectors $\Pi_S X$ and $\Pi_{S^\perp} X$ are probabilistically independent.  
Moreover, the random vector $\Pi_{S^\perp} X$ is not degenerate (i.e., $\Cov(\Pi_{S^\perp} X)$ is not zero). 
\end{assumption} 
\vspace{.5em} 
Assumption~\ref{assumption:independence} formalizes the idea that components of
$X$ in $S^\perp$ carry no information relevant for prediction, and may be regarded as pure noise. 
This has a natural interpretation: it reflects the idea that truly
predictive features in a system---such as object shape in an image---should vary independently
of background clutter. In Section~\ref{sec:invariance-and-relaxation}, we slightly relax this assumption,
for example allowing $X \sim \normal(\E[X], \Cov(X))$ with general covariance by exploiting additional
invariance in~\eqref{eqn:variational-objective-f}.

\vspace{.8em}
\noindent\noindent
\emph{$\mathsf{(A4)}$ Translation-and Rotational Invariance of $H$}.  

\vspace{.3em} 
In the absence of prior information that privileges a specific origin or coordinate system in $\R^d$, the learning problem itself should be modeled so as to be invariant under translations and rotations of the input space. We therefore focus on RKHSs that are translation- and rotation-invariant.

\vspace{.5em} 
\begin{assumption} 
\label{assumption:symmetry-of-RKHS}
The RKHS $H$ is associated with a kernel of the form $k(x,x') = \psi(\|x-x'\|^2)$, where 
$\psi(z) = \int_0^\infty e^{-tz}\mu(dt)$ for some probability measure $\mu$ that is not concentrated at zero, 
and that the right derivative $\psi'_+(0)$ exists.
\end{assumption} 

\vspace{.5em} 
This class includes standard examples such as Gaussian and Sobolev RKHSs.
The representation of $\psi$ is motivated by a classical result of Schoenberg~\cite{Schoenberg38},
which characterizes all functions $\psi$ for which $k(x,x') = \psi(\|x - x'\|^2)$ 
is a positive definite kernel on $\R^d$ for all $d$.
In our analysis, the invariance of $H$ is inherited by the objective
(see Section~\ref{sec:rotation-invariance-and-reparameterization}), while the existence of $\psi_+'(0)$ 
implies differentiability of $\psi$ on $(0, \infty)$ and enables first-order variational analysis.

\subsubsection{Formal Statements} 
\indent\indent
Our first main result is concerned with the population objective; existence of minimizers 
is given in Theorem~\ref{theorem:existence-of-global-minimizers}. For a matrix $U \in \R^{d \times d}$, 
we write $\ker(U):=\{x: Ux = 0\}$ for its null space. 

\vspace{.5em} 
\begin{theorem}[see Theorem~\ref{thm:population-recovery}]
\label{theorem:population-result-intro}
Let Assumptions~\ref{assumption:X-continuous}--~\ref{assumption:symmetry-of-RKHS} 
hold.  Then \vspace{.5em} 
\begin{itemize}
\item[$\mathsf{(a)}$] 
For every $\lambda > 0$, 
\begin{equation*}
	\ker(U_*) \supset S^\perp,~~~\forall (U_*, f_*)~\text{minimizing~\eqref{eqn:variational-objective-f}}.
\end{equation*} 
\vspace{0em} 
\item[$\mathsf{(b)}$] 
There exists $\lambda_0 > 0$ such that for all $0 < \lambda \le \lambda_0$,
\begin{equation*}
	\ker(U_*) = S^\perp,~~~\forall (U_*, f_*)~\text{minimizing~\eqref{eqn:variational-objective-f}}.
\end{equation*} 
\end{itemize}
\end{theorem} 

Theorem~\ref{theorem:population-result-intro} shows that, at the population level, all global minimizers annihilate directions orthogonal to the predictive subspace $S$. Moreover, when the regularization parameter $\lambda$ is sufficiently small, this annihilation is \emph{exact}, in the sense that no additional directions are annihilated beyond $S^\perp$. In this sense, the predictive subspace $S$ is \emph{identifiable} from the population minimizers.

Our second main result concerns the finite-sample objective; existence of minimizers is discussed in
Theorem~\ref{theorem:existence-of-global-minimizers-empirical}. We write $o_P(1)$ for random quantities
converging to zero in probability as $n \to \infty$.

\vspace{.5em} 
\begin{theorem}[see Theorem~\ref{thm:finite-sample-recovery}]
\label{theorem:finite-sample-result-intro}
Let Assumptions~\ref{assumption:X-continuous}--~\ref{assumption:symmetry-of-RKHS} 
hold.  Then \vspace{.5em} 
\begin{itemize}
\item[$\mathsf{(a)}$] 
For every $\lambda > 0$, 
\begin{equation*}
	\lim_{n \to \infty} \P( \rank(U_{n, *}) \le \dim(S),~~~\forall (U_{n, *}, f_{n, *})~\text{minimizing~\eqref{eqn:variational-objective-f-empirical}}) = 1.
\end{equation*}
\item[$\mathsf{(b)}$] 
For a given $0 < \lambda \le \lambda_0$, where $\lambda_0$ is specified in Theorem~\ref{theorem:population-result-intro}$\mathsf{(b)}$
\begin{equation*}
	\lim_{n \to \infty} \P( \rank(U_{n, *}) = \dim(S),~~~\forall (U_{n, *}, f_{n, *})~\text{minimizing~\eqref{eqn:variational-objective-f-empirical}}) = 1.\end{equation*} 
\item[$\mathsf{(c)}$]
For a given $0 < \lambda \le \lambda_0$, where $\lambda_0$ is specified in Theorem~\ref{theorem:population-result-intro}$\mathsf{(b)}$, as $n \to \infty$,
\begin{equation*}
	\sup\left\{\norm{\Pi_{\ker(U_{n, *})}- \Pi_{S^\perp}} \mid (U_{n, *}, f_{n, *})~\text{minimizing~\eqref{eqn:variational-objective-f-empirical}} \right\} = o_P(1).
\end{equation*}  
\end{itemize}
\end{theorem} 

Theorem~\ref{theorem:finite-sample-result-intro} provides a finite-sample analogue of 
Theorem~\ref{theorem:population-result-intro}. Theorem~\ref{theorem:finite-sample-result-intro}$\mathsf{(a)}$ 
shows that empirical minimizers have rank at most $\dim(S)$, while
Theorem~\ref{theorem:finite-sample-result-intro}$\mathsf{(b)}$ 
establishes exact recovery of the intrinsic subspace dimension $\dim(S)$ in the small-$\lambda$
regime, both holding with high probability. Theorem~\ref{theorem:finite-sample-result-intro}$\mathsf{(c)}$ 
strengthens this conclusion by showing that the null space of $U_{n, *}$ converges, 
in probability, to $S^\perp$. Thus, $\ker(U_{n, *})^\perp$ yields consistent recovery of the predictive subspace $S$ and its dimension $\dim(S)$.

We conclude with two remarks.  \vspace{.5em}
\begin{itemize}
\item A notable feature of Theorem~\ref{theorem:finite-sample-result-intro} is the
consistency it establishes for the intrinsic subspace dimension $\dim(S)$.
In statistical learning literature, explicit regularization mechanisms—such
as nuclear-norm penalization or early stopping—are commonly employed to learn the 
intrinsic low-dimensional structure of the data~\cite{HastieTiFr09}. Theorem~\ref{theorem:finite-sample-result-intro}
shows that, for the compositional kernel model, recovery of the data's intrinsic
dimension does not rely on such mechanisms.

 \vspace{.5em} 
\item Theorem~\ref{theorem:finite-sample-result-intro}$\mathsf{(a)}$ does not follow directly from
Theorem~\ref{theorem:population-result-intro}$\mathsf{(a)}$ together with convergence of
empirical minimizers $U_{n, *}$ to population minimizers $U_*$. Indeed, convergence in norm does not, in
general, imply convergence of rank. Establishing the low-dimensional structure of
empirical minimizers $U_{n, *}$ therefore requires additional arguments beyond consistency of the
minimizers themselves.
\end{itemize} 
\vspace{.5em}

Section~\ref{sec:technical-insight} elucidates the key geometric properties of the population
objective near its minimizers---the \emph{sharpness} property---that are 
responsible for this subspace dimension consistency.

\newcommand{\col}{{\rm col}}
\subsection{Technical Insights} 
\label{sec:technical-insight}
\indent\indent 
This section explains the main technical idea underlying the subspace dimension consistency in 
Theorem~\ref{theorem:finite-sample-result-intro}. We restrict attention to the regime $0 < \lambda \le \lambda_0$, where $\lambda_0$ is specified in Theorem~\ref{theorem:population-result-intro}. We explain proof mechanism underlying Theorem~\ref{theorem:finite-sample-result-intro}$\mathsf{(b)}$, namely,
\begin{equation*}
	\lim_{n \to \infty} \P( \rank(U_{n, *}) \le \dim(S),~~~\forall (U_{n, *}, f_{n, *})~\text{minimizing~\eqref{eqn:variational-objective-f-empirical}}) = 1.
\end{equation*}

We define the population objective 
\begin{equation}
\label{eqn:definition-of-J}
	\mathcal{J}(U, \lambda) := \min_{f \in H}
		\E \big[(Y - f(UX))^2\big] + \lambda \norm{f}_{H}^2 .
\end{equation}
Since $H$ is rotationally invariant, one can show that $\mathcal{J}$ depends on $U$ only through $U^\top U$.
Accordingly, we can reparametrize the objective as
$J(\Sigma,\lambda) := \mathcal{J}(U,\lambda)$
where $U$ is any matrix satisfying $\Sigma = U^\top U$. 
If $(U_{*}, f_{*})$ minimizes \eqref{eqn:variational-objective-f}, then
$\Sigma_{*} := U_{*}^\top U_{*}$ minimizes
\begin{equation}
	\minimize_\Sigma J(\Sigma, \lambda)~~~\text{subject to}~~~\Sigma \succeq 0. 
\end{equation}
Here, $\Sigma \succeq 0$ denotes the positive semidefinite constraint. By Theorem
\ref{theorem:population-result-intro}$\mathsf{(b)}$, $\Range(\Sigma_*) = S$.

Similarly, we can define the empirical objective $J_n(\Sigma,\lambda)$ by replacing $\E$ with $\E_n$.
If $(U_{n, *}, f_{n, *})$ minimizes \eqref{eqn:variational-objective-f-empirical}, then
$\Sigma_{n, *} := U_{n, *}^\top U_{n, *}$ minimizes
\begin{equation}
	\minimize_\Sigma J_n(\Sigma, \lambda)~~~\text{subject to}~~~\Sigma \succeq 0.
\end{equation} 
Let's denote 
\begin{equation}
\label{eqn:definition-of-sym-+-d-M}
	\Sym_+^d := \{\Sigma: \Sigma \succeq 0\}
	~~\text{and}~~
	\mathcal{M} := \{\Sigma: \Sigma \succeq 0,~\rank(\Sigma) = \dim(S)\}.
\end{equation} 
Note $\mathcal{M}$ is a submanifold of the Euclidean space of all symmetric matrices $\Sym^d:=\{\Sigma: \Sigma = \Sigma^\top\}$.
\vspace{.5em}

\noindent\noindent
\textbf{Main Properties}~~
The argument rests on three key properties of $J$ and $J_n$. \vspace{.5em}
\begin{itemize}
\item (Sharpness Property) 	
	For every minimizer $\Sigma_*$ of $J$, there is $\rho > 0$ such that 
	for every matrix $W$ in the tangent cone of  $\Sym_+^d$ at $\Sigma_*$
	with $\Range(W) \subseteq \Range(\Sigma_*)^\perp$, 
	\begin{equation}
	\label{eqn:sharpness-property-intro}
		\lim_{s \to 0^+} \frac{1}{s} (J(\Sigma_* + sW, \lambda) - J(\Sigma_*, \lambda)) \ge \rho \norm{W}.
	\end{equation} 
	
	We give two interpretations of~\eqref{eqn:sharpness-property-intro}. Geometrically, the normal space to the submanifold 
	$\mathcal{M}$ with respect to the ambient space $\Sym^d$ at any $\Sigma \in \mathcal{M}$
	consists of all symmetric matrices $W$ satisfying $\Range(W) \subseteq \Range(\Sigma)^\perp$~\cite{HiriartMa12}. Hence, 
	the sharpness property~\eqref{eqn:sharpness-property-intro} means that
	$J$ grows locally linearly when deviating perpendicularly away from the manifold $\mathcal{M}$. 
	Statistically, since $\Range(\Sigma_*) = S$, this sharpness property~\eqref{eqn:sharpness-property-intro} 
	says that perturbation of $\Sigma_*$ in any direction $W$ orthogonal to the predictive 
	subspace $S$ would induce a local linear growth. 
	

\vspace{.5em}
\item (Gradient Regularity) The gradient $\grad J$ exists and is continuous on $\Sym_+^d$
(see Section~\ref{sec:first-variation} for formal definition of $\grad J$).


\vspace{.5em} 
\item (Uniform Convergence) 
	The function value and gradient of the empirical objective $J_n$ converge uniformly 
	to those of the population objective $J$ on the positive semidefinite cone 
	$\Sym_+^d$, i.e., 
	\begin{equation}
	\begin{split} 
		\sup_{\Sigma \in \Sym_+^d} |J_n(\Sigma, \lambda) - J(\Sigma, \lambda)| \to 0,~~~~~
		\sup_{\Sigma \in \Sym_+^d} \norm{\grad J_n(\Sigma, \lambda) - \grad J(\Sigma, \lambda)} \to 0
	\end{split} 
	\end{equation} 
	almost surely as the sample size $n \to \infty$. 
\end{itemize}

\vspace{.3cm}
\noindent\noindent
\textbf{Main Argument}~~
With these properties in mind, our main argument, which is largely inspired by
sensitivity analysis of minimizers' set identifiability in the optimization 
literature~\cite{DrusvyatskiyLe14}, is as follows. \vspace{.3cm}
\begin{itemize}
\item The uniform convergence property of objective values ensures that the empirical 
	minimizers converge to their population counterparts as $n \to \infty$.
	Since all population minimizers $\Sigma_*$ fall within $\mathcal{M}$, this implies the 
	empirical minimizers $\Sigma_{n, *}$ converge to $\mathcal{M}$.
	
	\vspace{.3cm}
\item To further show the empirical minimizers $\Sigma_{n, *}$ fall within $\mathcal{M}$, 
	the sharpness property is crucial. The sharpness property ensures the empirical 
	objective, whose gradient is uniformly close to that of the population objective, 
	is still with high probability sharp, i.e., increasing at least linearly when moving in perpendicular directions away 
	from $\mathcal{M}$ near the population minimizers. Consequently, $\Sigma_{n,*} \in \mathcal{M}$ 
	with probability tending to one.
\end{itemize}
\vspace{.3cm}

In summary, while uniform convergence of objective values implies convergence 
of minimizers, the sharpness property offers the stability that preserves the minimizers' 
low rankness property---a case of set identifiability---under perturbations. 
Formal proofs, including the establishment of sharpness properties, gradient regularity, and 
uniform convergence arguments---are given in Sections~\ref{sec:preliminaries}--\ref{sec:empirical-minimizers}.

\subsection{Related Works}
\indent\indent
This paper contributes to a growing line of work on \emph{feature learning}, understood as the emergence of low-dimensional, predictive representations in compositional models through the minimization of data-driven objectives~\cite{BengioCoVi13}. Classical kernel methods provide a well-established variational framework for regression, but because the kernel is fixed \emph{a priori}, they do not themselves exhibit representation learning. Recent extensions therefore modify classical kernel methods by introducing learnable representations, leading to compositional kernel models of the form $x \mapsto f(Ux)$. This formulation allows the representation $U$ and predictor $f$ to be optimized jointly within a single variational objective and and has attracted recent interest, in part due to an analogy to two-layer neural networks~\cite{FukumizuBaJo09, Adit24A, FollainBa24, LiRu25, LiRub25}.

Within this framework, we study feature learning by focusing on the \emph{variational geometry} of the population objective and its \emph{implications under finite sampling}. From this perspective, our contribution differs significantly from existing approaches to $x \mapsto f(Ux)$ models, a distinction we clarify below by contrasting our analysis with (i) variational formulations with structural constraints and (ii) algorithmic methods that investigate feature learning  through iterative dynamics.

\vspace{.5em} 
\begin{itemize}
\item (Variational formulations with constraints). An earlier line of work by Fukumizu, Bach and Jordan studied an analogous compositional kernel formulation called \emph{kernel dimension reduction}~\cite{FukumizuBaJo04, FukumizuBaJo09}. In their original formulation, the representation $U$ is subject to an explicit structural constraint: its rank is fixed to equal the dimension of the true predictive subspace $S$. As a result, these approaches presuppose prior knowledge of the intrinsic subspace dimension $\dim(S)$ and encode directly a low-dimensionality constraint into the variational objective.

By contrast, our analysis imposes \emph{no} constraints on $U$, and in particular \emph{no} rank constraints. Our contributions are two-fold. First, under the setting considered in this paper, we show that the predictive subspace can be recovered at the global minimizers of the population objective. In particular, this low-dimensional structure is not enforced \emph{a priori} through the model specification. Second, we identify a sharpness property of the population objective near its minimizers, which ensures that the low-dimensional structure is stable and inherited by finite-sample minimizers without the need for explicit regularizations (e.g., nuclear norm penalty). Together, these results provide a variational explanation of feature learning in the compositional kernel model, whereby low-dimensional predictive structure arises from the geometry of the objective rather than from explicit structural constraints. 

\item (Algorithmic methods for feature learning)
A complementary line of work studies \emph{algorithmic approaches} to feature learning in kernel methods, most notably the Average Gradient Outer Product (AGOP) method~\cite{Adit24A}. These approaches can be viewed as operating on the same underlying variational objective considered in this paper, but through iterative, data-driven updates of the representation $U$. AGOP-type methods rely on iterative, data-driven updates of the representation ($U^{(t+1)} = {\rm AGOP}(U^{(t)})$) and have been studied theoretically—yielding generalization guarantees that improve upon standard kernel methods—and empirically, with strong performance in multiple application domains such as tabular data prediction and image classification, matching state-of-the-art performance in representative benchmarks~\cite{Adit24A, ZhuDaDrFa25, Beagleholeetal25}. 

Our work takes a different, but complementary, perspective. Rather than analyzing iterative algorithmic schemes, we study the variational geometry of the objective itself and characterize the structure of its global minimizers, together with the stability of this structure under finite sampling. This analysis provides fundamental population-level and finite-sample guarantees for the emergence of low-dimensional predictive structure, \emph{independently of any particular algorithmic dynamics}. Interestingly, these results are consistent with empirical observations in AGOP-based methods, where representations frequently display effective low-dimensionality in practice~\cite{RadhakrishnanBeDr25}.  At the same time, because AGOP is not derived from an explicit variational principle, its precise relationship to the underlying variational landscape—and, in particular, to the structure of its minimizers—remains to be fully understood~\cite{RadhakrishnanBeDr25}. Clarifying this connection represents a natural and promising direction for future work.
\end{itemize}

\subsection{Technical Contributions} 
The technical contributions of this work are two-fold.

First, we identify sharpness of the population objective, formalized in~\eqref{eqn:sharpness-property-intro}, as a key mechanism preserving the low-dimensional structure of population minimizers. This explains why low-rank solutions persist under finite sampling without explicit penalizations, providing an alternative, but compatible, explanation within statistical learning. A broader discussion of this variational-analytic perspective and its connection to statistics perspective on regularization is deferred to Section~\ref{sec:variational-perspective-of-regularization}. 

Second, we characterize when such sharpness holds for the compositional kernel objective at the population level. A key ingredient of the analysis is a Fourier-analytic expansion of the bivariate map $(x,x') \mapsto \psi'(\|x-x'\|^2)$ appearing in directional derivatives of the objective. This expansion establishes a negative definiteness property of the map, which controls the sign of directional derivatives and underlies the proof of the sharpness property. See Section~\ref{sec:sharpness-property} for details.

\textbf{Note}. After the initial posting of this manuscript on arXiv, the authors prepared a separate, self-contained work~\cite{LiRu25} documenting basic mathematical facts concerning the population model~\eqref{eqn:variational-objective-f}. References to that work are confined to the preliminary Section~\ref{sec:preliminaries} and are not essential to the main arguments, which are otherwise self-contained.

\subsection{Notation} 
Let $x \in \mathbb{R}^d$ be a vector, written as $x = (x_1,\dots,x_d)^\top$. We write $\langle x, y \rangle := \sum_{i} x_i y_i$ for the Euclidean inner product between vectors $x,y \in \mathbb{R}^d$, and denote the associated Euclidean norm by $\|x\| := \sqrt{\langle x,x \rangle}$. Similarly, let $A,B \in \mathbb{R}^{d \times d}$ be matrices, where $A_{ij}$ and $B_{ij}$ denote their $(i,j)$-th entries. We write $\langle A, B \rangle := \sum_{i,j} A_{ij} B_{ij}$, and denote the associated Frobenius norm by $\|A\| := \sqrt{\langle A,A \rangle}$. 
For a vector $x \in \mathbb{R}^d$, and a positive semidefinite $\Sigma \in \Sym_+^d$, we denote 
$\|x\|_{\Sigma} := \sqrt{\langle x, \Sigma x \rangle}$. 
%

The Fourier transform of a function \( f : \R^d \to \C \) is denoted by \( \hat{f} \). We use the convention:
\begin{equation*}
	\hat{f}(\omega) = \int_{\R^d} f(x) e^{-2\pi i\langle x, \omega\rangle} dx,~~~f(x) = \int_{\R^d} \hat{f}(\omega) e^{2\pi i \langle x, \omega \rangle} d\omega. 
\end{equation*} 
 The reason for the $2\pi i$ choice is that it makes the least appearance of dimensional constants.


\section{Preliminaries} 
\label{sec:preliminaries} 
\subsection{Translation and Rotationally Invariant RKHS} 
We give a quick review about reproducing kernel Hilbert spaces, and construct the RKHS 
associated with the translation-and rotationally invariant kernel $(x,x') \mapsto \psi(\|x-x'\|^2)$.
More general references on RKHS can be found in~\cite{Aronszajn50, CuckerSm02}. 
  
Let $H$ be a Hilbert space of real-valued function on $\R^d$, defined via the Fourier 
transform by 
\begin{equation}
\label{eqn:Hilbert-space-construction}
	H := \left\{f: f \in  L_2(\R^d), 
		\norm{f}_H^2 = \int_{\R^d} \frac{|\hat{f}|^2(\omega)}{k(\omega)} d\omega < \infty\right\},
\end{equation}
where we require $k \in L_1(\R^d)$ to be strictly positive, even and normalized so that $\int_{\R^d} k(\omega) d\omega = 1$.
The following result is standard; we include a proof in Appendix~\ref{sec:proof-of-lemma-continuous-embedding} for completeness.
\vspace{.5em} 
\begin{lemma}
\label{lemma:continuous-embedding}
For every $f \in H$, we have $\norm{f}_{L_\infty} \le \norm{f}_H$, and $f$ has a continuous representative, and $\lim_{x \to \infty} f(x) = 0$.
\end{lemma} 

\vspace{.5em} 
Henceforth, we identify $H$ with its image in $C(\R^d)$ via this continuous embedding, and work with the continuous representatives of its elements. In particular, $H$ is a reproducing kernel Hilbert space. We then define the kernel function 
\begin{equation*}
	\kernel(x) = \int_{\R^d} k(\omega) e^{2\pi i \langle x, \omega\rangle} d\omega.
\end{equation*} 
Note that $\kernel$ is real-valued as $k$ is even, and continuous as $k \in L_1(\R^d)$. 
Let $\kernel(x, x') = \kernel(x-x')$.
The continuous function $\kernel(x, \cdot)$ lies in $H$, and has Fourier transform 
$k(\omega) e^{-2\pi i\langle \omega, x\rangle}$. These functions $\kernel(x, \cdot)$ 
are known as \emph{reproducing kernel}. Each $f \in H$ can be represented as 
\begin{equation*}
	f(x) = \langle f, \kernel(x, \cdot)\rangle_H,~~~\forall x\in \R^d.
\end{equation*}
The RKHS $H$ satisfies the following standard density property; a proof is included in
Appendix~\ref{sec:proof-of-lemma-denseness}.

\vspace{.5em}
\begin{lemma}
\label{lemma:denseness-of-H-in-L-2}
For any probability measure $\nu$ on $\R^d$, the Hilbert
space $H$ is dense in $L_2(\nu)$.
\end{lemma} 

\vspace{.5em}

The Fourier construction above makes $H$ invariant under translations.
Throughout this paper, we further restrict attention to rotationally invariant kernels. 
Specifically, we assume there exists a function $\psi : [0,\infty) \to \R$ such that
$\kernel(x) = \psi(\|x\|^2)$, and 
	$\kernel(x,x') = \psi(\|x-x'\|^2)$, where 
\begin{equation}
\label{eqn:psi-definition}
	\psi(z) = \int_0^\infty e^{-tz} \mu(dt),~~~\forall z \ge 0
\end{equation}
for some probability measure $\mu$ on $[0, \infty)$ that is not fully concentrated at zero.
This expression is motivated by a classical theorem of Schoenberg, which says
$\psi(\|\cdot\|^2)$ defines a positive-definite kernel on $\R^d$ for all dimensions $d$ if and only if 
$\psi$ is non-constant and admits a representation of~\eqref{eqn:psi-definition}~\cite{Schoenberg38}.

Under this condition, the kernel $\kernel$ and its Fourier pair $k$ have the expression
\begin{equation}
\begin{split} 
	\kernel(x) = \profile(\norm{x}_2^2) = \int_0^\infty e^{-t\norm{x}_2^2}\mu(dt),\quad\quad
	k(\omega) = \int_0^\infty 
		(\pi t^{-1})^{d/2} e^{- \pi^2 t^{-1} \norm{\omega}_2^2}\mu(dt).
\end{split} 
\end{equation}
This completes the construction of the RKHS $H$ associated with the kernel
$(x,x') \mapsto \psi(\|x-x'\|_2^2)$. The space $H$ is constructed explicitly via~\eqref{eqn:Hilbert-space-construction} and satisfies the 
properties stated in Lemmas~\ref{lemma:continuous-embedding} and~\ref{lemma:denseness-of-H-in-L-2}.

\subsubsection{Examples} 
Here are two canonical examples of translation and rotationally invariant  RKHS.

\vspace{.5em}
\begin{example} [Gaussian RKHS]
The Gaussian RKHS corresponds to the case where
\begin{equation*}
	\kernel(x) = e^{-\norm{x}_2^2},\quad~\text{and}~\quad k(\omega) = \pi^{d/2} e^{- \pi^2 \norm{\omega}_2^2}.
\end{equation*} 
The measure $\mu$ is the Dirac measure at $1$, with $\psi(z) = \int e^{-tz}\mu(dt) = e^{-z}$. $\clubsuit$
\end{example} 

\vspace{.5em}

\begin{example} [Sobolev RKHS]
Let $\gamma > 0$. The Sobolev RKHS corresponds to the case where
\begin{equation*}
	\kernel(x) = \frac{1}{4^\gamma \Gamma(\gamma)} \int_0^\infty t^{-\gamma-1} e^{-\frac{1}{4t}} e^{- t\norm{x}^2} dt,\quad~\text{and}~\quad k(\omega) = \frac{(4\pi)^{d/2}\Gamma(\gamma+d/2)}{\Gamma(\gamma)(1+ 4\pi^2 \norm{\omega}^2)^{\gamma+d/2}}.
\end{equation*} 
In the above, $\Gamma(\cdot)$ denotes the Gamma function. 
The measure $\mu$ is $\mu(dt) = \frac{1}{4^\gamma \Gamma(\gamma)} t^{-\gamma-1} e^{-\frac{1}{4t}} dt$, which coincides with the Inverse-Gamma distribution with shape parameter $\gamma$ and rate parameter $1/4$. $\clubsuit$
\end{example}

\subsection{Kernel Ridge Regression} 
Let $X \in \R^d$ be a random vector, and $Y \in \R$ be a random variable with $\E[Y^2] < \infty$. 
We define 
\begin{equation}
	\mathcal{I}(f, U, \lambda) :=  \E[(Y- f(U X))^2] + \lambda \norm{f}_H^2
\end{equation}
and consider the following minimization problem, which we will refer to as \emph{kernel ridge regression}: 
\begin{equation}
	\mathcal{J}(U, \lambda) = \min_{f \in \H} \mathcal{I}(f, U, \lambda).
\end{equation}
The existence of a unique minimizer is well-known~\cite{CuckerSm02}. We denote this unique minimizer by 
$f_{U, \lambda}$. 

Below we give the Euler-Lagrange characterization satisfied by the minimizer.
\vspace{.5em}
\begin{lemma}
\label{lemma:euler-lagrange-identity}
The minimizer $f_{U, \lambda}$ to this problem satisfies 
\begin{equation*}
	\E[(Y- f_{U, \lambda}(U X)) g(U X)] = \lambda \langle f_{U, \lambda}, g\rangle_H,~~~\forall g\in H.
\end{equation*}
\end{lemma} 
\begin{proof}
This follows by setting the first variation of $f \mapsto \mathcal{I}(f,U,\lambda)$ equal to zero.
\end{proof} 

We record a trivial but useful upper bound on $\mathcal{J}(U,\lambda)$.
\vspace{.5em}
\begin{lemma}
\label{lemma:upper-bound-U}
We always have $\E[Y^2] \ge \mathcal{J}(U, \lambda)$. The inequality is strict when $\E[Y|UX] \neq 0$.
\end{lemma}
\begin{proof}
Plug in the test function $f \equiv 0$, we get $\E[Y^2] = \mathcal{I}(0, U, \lambda) \ge \mathcal{J}(U, \lambda)$.
We now characterize the equality case $\mathcal{J}(U,\lambda)=\E[Y^2]$. Then the minimizer $f_{U,\lambda}$  must be identically zero, and Lemma~\ref{lemma:euler-lagrange-identity} yields $\E[Yg(UX)]=0$ for all $g\in H$. Since $H$ is dense in $L_2(\P_{UX})$ by Lemma~\ref{lemma:denseness-of-H-in-L-2}, where $\P_{UX}$ denotes the distribution of $UX$, this identity extends to all $g\in L_2(\P_{UX})$, implying $\E[Y| UX]=0$. 
\end{proof} 

We then record a lower bound on $\mathcal{J}(U,\lambda)$.
\vspace{.5em}
\begin{lemma}
\label{lemma:lower-bound-U}
Let $(X', Y')$ denote an independent copy of $(X, Y)$. We always have 
\begin{equation}
	\mathcal{J}(U, \lambda) \ge \E[Y^2 ] - \lambda^{-1} \E[YY' K(UX, UX')].
\end{equation} 
\end{lemma} 
\begin{proof} 
We first expand the objective: 
\begin{equation*}
	\mathcal{I}(f, U, \lambda) = \E[Y^2] - 2\E[Y f(U X)] + \E[f(U X)^2] + \lambda \norm{f}_H^2.
\end{equation*}
Applying the reproducing property and Cauchy–Schwarz inequality yields
\begin{equation*}
\begin{split}
	\E[Y f(U X)] &= \E[Y \langle \kernel(U X, \cdot), f \rangle_H] \\
		&= \langle \E[Y  \kernel(U X, \cdot)], f \rangle_H
		\le \norm{\E[Y  \kernel(U X, \cdot)]}_H \norm{f}_H.
\end{split}
\end{equation*}
The squared RKHS norm admits an explicit expression: 
\begin{equation*}
	\norm{\E[Y  \kernel(U X, \cdot)]}_H^2 
		= \langle \E[Y  \kernel(U X, \cdot)], \E[Y'  \kernel(U X', \cdot)]\rangle_H
		= \E[YY'\kernel(U X, U X')]
\end{equation*} 
where $(X', Y') \sim \P$ is an independent copy of $(X, Y)$. 
Combining the above gives a lower bound 
\begin{equation*}
\begin{split} 
	\mathcal{I}(f, U, \lambda) 
		&\ge \E[Y^2] - 2 \sqrt{\E[YY' \kernel(U X, U X')]} \norm{f}_H + \lambda \norm{f}_H^2.
\end{split} 
\end{equation*} 
Minimizing the right-hand side over $\norm{f}_H \in [0, \infty)$ gives the bound as desired. 
\end{proof}

\subsection{Rotation Invariance and Reparameterization}
\label{sec:rotation-invariance-and-reparameterization}
\indent\indent
The rotational invariance of the RKHS $H$ greatly simplifies the geometry of $\mathcal{J}$. In particular, 
$\mathcal{J}(U, \lambda)$ is invariant under the rotation group action.

\vspace{.4em} 
\begin{lemma}
\label{lemma:rotation-invariance}
For every orthogonal matrix $O$ and every matrix $U$, we have
$
\mathcal{J}(U,\lambda)=\mathcal{J}(OU,\lambda)
$.
Moreover, the minimizers satisfy $f_{U, \lambda}(Ux) = f_{OU, \lambda}(OUx)$ for all $x \in \R^d$. 
\end{lemma} 
\begin{proof}
Define the operator $T_O$ on functions by $(T_O f) (x) = f(Ox)$. Since $\widehat{T_Of}(\omega)=\hat f(O^{\top}\omega)$, we get 
\[
\|T_Of\|_{H}^{2}
=\int \frac{|\hat f(O^{\top}\omega)|^{2}}{k(\omega)}\,d\omega
=\int \frac{|\hat f(\xi)|^{2}}{k(O\xi)}\,d\xi
=\int \frac{|\hat f(\xi)|^{2}}{k(\xi)}\,d\xi
=\|f\|_{H}^{2},
\]
where we used $|\det O|=1$ and $k(O\xi)=k(\xi)$ since $K$ is radial. Hence,
$
	\mathcal{I}(T_O f, U, \lambda) 
	 =  \mathcal{I}(f, OU, \lambda)$ for all $f$. 
By taking infimum over $f$ on both sides, we get $\mathcal{J}(U, \lambda) = \mathcal{J}(OU, \lambda)$, 
and the minimizers are related by $T_Of_{OU,\lambda} = f_{U,\lambda}$. Hence, 
$f_{U, \lambda}(Ux) =f_{OU,\lambda}(OUx)$.
\end{proof} 

\vspace{.5em}

As a result of rotational invariance, we can reparameterize $\mathcal{J}$ in terms of $\Sigma = U^\top U$. 
Formally, for every positive semidefinite matrix $\Sigma$, we define
\begin{equation}
\label{eqn:definition-of-J}
  J(\Sigma,\lambda) := \mathcal{J}(U,\lambda),
\end{equation}
where $U$ on the right-hand side is any matrix such that $U^{\top}U = \Sigma$. Lemma~\ref{lemma:rotation-invariance}
shows that this is well-defined: if $U_1^{\top}U_1 = U_2^{\top}U_2$, then $U_1 = OU_2$ for some orthogonal $O$, and thereby
$\mathcal{J}(U_1,\lambda) = \mathcal{J}(U_2,\lambda)$. Similarly, for every positive semidefinite matrix $\Sigma$, we define
the function $F_{\Sigma, \lambda}$ by
\begin{equation}
\label{eqn:definition-of-F}
	F_{\Sigma, \lambda}(x) := f_{U, \lambda}(Ux), 
\end{equation}
where $U$ is any matrix such that $U^\top U=\Sigma$; by the same invariance argument, this is well-defined.

It is therefore convenient to translate Lemmas~\ref{lemma:upper-bound-U} and~\ref{lemma:lower-bound-U}
into the $\Sigma$-formulation.

\vspace{.5em}
\begin{lemma}
\label{lemma:bounds-on-Sigma}
For every $\Sigma \in \Sym_+^d$, we have 
\begin{equation}
	\E[Y^2] \ge J(\Sigma, \lambda) \ge \E[Y^2] - \lambda^{-1} \E[YY' \psi(\|X-X'\|_\Sigma^2)].
\end{equation}
\end{lemma}

In a nutshell, these reparameterizations from the $U$–space to the $\Sigma$–space remove the redundancy created by rotational invariance. Accordingly, our subsequent analysis is carried out primarily in the $\Sigma$ formulation.

\subsection{First Variation} 
\label{sec:first-variation}
\indent\indent 
It is straightforward to show that $J$ is continuous in its parameter $\Sigma$
whenever $\E[Y^2] < \infty$; see our companion work~\cite[Proposition 2.8]{LiRu25}.
We now turn to the first variation of $J$ with respect to $\Sigma$.
The detailed study is given in our companion paper; here we only recall 
the setup and result.
Since $\Sigma$ is restricted to the positive semidefinite cone $\Sym_+^d$,
admissible perturbations are restricted. Accordingly, we define 
\begin{equation*}
	T_{\Sym_+^d}(\Sigma) := \left\{W \in \R^{d \times d}: \exists s_0 > 0~\text{such that}~\Sigma + s W \in \Sym_+^d~\text{for all}~s \in [0, s_0)\right\}.
\end{equation*}
For every $\Sigma \in \Sym_+^d$ and $W \in T_{\Sym_+^d}(\Sigma)$, the directional 
derivative of $J$ at $\Sigma$ along $W$ is defined by 
\begin{equation*}
\mathrm{D} J(\Sigma, \lambda)[W] := \lim_{s \to 0^+} \frac{1}{s} \left( J(\Sigma + s W, \lambda) - J(\Sigma, \lambda) \right)
\end{equation*}
whenever the limit exists. For readers familiar with convex analysis, the set $T_{\Sym_+^d}(\Sigma)$ coincides with
the tangent cone of $\Sym_+^d$ at $\Sigma$~\cite{RockafellarTyWe98}. 

We first provide sufficient conditions under which $J$ is G\^ateaux differentiable (i.e., directionally differentiable), and 
provide an explicit formula for its directional derivative. Define for $\Sigma \in \Sym_+^d$,
\begin{equation*}
	r_{\Sigma,\lambda}(x,y) := y - F_{\Sigma,\lambda}(x).
\end{equation*}
When $\Sigma = U^\top U$, it satisfies $r_{\Sigma,\lambda}(x,y) = y - f_{U,\lambda}(Ux)$.
Lemma~\ref{theorem:first-variation-formula} follows from~\cite[Theorem 4.2]{LiRu25}.

\vspace{.5em} 
\begin{lemma}
\label{theorem:first-variation-formula}
The directional derivative $\mathrm{D} J(\Sigma, \lambda)[W]$ exists 
for every $\Sigma \in \Sym_+^d$, and every matrix $W \in T_{\Sym_+^d}(\Sigma)$, and satisfies 
\begin{equation}
\label{eqn:first-variation-formula}
	\mathrm{D} J(\Sigma, \lambda)[W] =
		- \frac{1}{\lambda} \E \left[r_{\Sigma, \lambda}(X, Y) r_{\Sigma, \lambda}(X', Y') \cdot \psi'\!\bigl(\|X - X'\|_\Sigma^2\bigr)\,(X - X')^\top W (X - X')\right].
\end{equation}
In the above, $(X', Y')$ denote an independent copy of $(X, Y)$. 
\end{lemma} 

\vspace{.5em} 
An important observation is that the directional 
derivative $W \mapsto \mathrm{D} J(\Sigma, \lambda)[W]$ is linear in $W$. Following Lemma~\ref{theorem:first-variation-formula},
it admits the representation
\begin{equation*}
	\mathrm{D}J(\Sigma,\lambda)[W]
	= \bigl\langle \nabla J(\Sigma,\lambda),\, W \bigr\rangle,
\end{equation*}
where $\nabla J(\Sigma,\lambda)$ is the symmetric matrix (which we will refer to as \emph{gradient} of $J$) given by
\begin{equation}
\label{eqn:gradient-of-J}
	\nabla J(\Sigma,\lambda)
	:=
	- \frac{1}{\lambda}
	\E\!\left[
		r_{\Sigma,\lambda}(X,Y)\, r_{\Sigma,\lambda}(X',Y')\,
		\psi'\!\bigl(\|X - X'\|_\Sigma^2\bigr)\,
		(X - X')(X - X')^\top
	\right].
\end{equation}
We next record a continuity property of $\nabla J$, established in~\cite[Corollary 4.10]{LiRu25}. 

\vspace{.5em} 
\begin{lemma}
For every $\lambda > 0$, the map $\Sigma \mapsto \nabla J(\Sigma,\lambda)$ is continuous. 
\end{lemma}

%
\vspace{.5em} 

To aid intuitions, we illustrate the directional derivative formula
$\mathrm{D}J(\Sigma,\lambda)[W]$ \eqref{eqn:first-variation-formula} and the gradient 
formula $\grad J(\Sigma, \lambda)$ \eqref{eqn:gradient-of-J} with concrete examples
in Section~\ref{sec:illustration-under-finite-atomic-measures}. 


\subsection{Law of Large Numbers}
\indent\indent
Let $\{(X_i, Y_i)\}_{i=1}^\infty$ be i.i.d.\ samples drawn from the distribution of $\P$.
For each $n \in \N$, the empirical expectation over $n$ samples is defined by 
\begin{equation*}
	\E_n[F(X, Y)] := \frac{1}{n} \sum_{i=1}^n F(X_i, Y_i).
\end{equation*}
Suppose that $\E[|F(X, Y)|] < \infty$, then the law of large numbers imply that almost surely 
\begin{equation*}
	\lim_{n \to \infty} \E_n[F(X, Y)] = \E[F(X, Y)].
\end{equation*}
By replacing $\E$ with the empirical expectation $\E_n$, we can define the empirical version 
$I_n(f, U, \lambda)$ of $I(f, U, \lambda)$, whose minimizer is $f_{U, \lambda, n}$ and the minimal 
value is $\mathcal{J}_n(U, \lambda)$. Proceeding as in the population setting, we further define the empirical counterparts: 
the objective $J_n(\Sigma,\lambda)$,  its directional derivative
$\mathrm{D}J_n(\Sigma,\lambda)[W]$, and gradient $\nabla J_n(\Sigma, \lambda)$ for all $\Sigma \in \Sym_+^d$ and $W \in \T_{\Sym_+^d}(\Sigma)$.

We state a \emph{uniform} law of large numbers, taken from our companion work~\cite[Proposition 4.4]{LiRu25}.

\vspace{.5em} 
\begin{lemma}
\label{lemma:law-of-large-number-constrained-set}
Let $\lambda > 0$ and $M < \infty$. Almost surely as $n \to \infty$: 
\begin{equation*}
\begin{split}
	&\sup_{\Sigma \in \Sym_+^d, \norm{\Sigma} \le M} 
	\bigl|J_n(\Sigma, \lambda) - J(\Sigma, \lambda)\bigr|
	\to 0,~~~\qquad \sup_{\Sigma \in \Sym_+^d, \norm{\Sigma} \le M} 
	\bigl\|
	\nabla J_n(\Sigma, \lambda)
	-
	\nabla J(\Sigma, \lambda)
	\bigr\|
	\to  0.
\end{split} 
\end{equation*}
\end{lemma} 
The restriction $\{\Sigma: \|\Sigma\| \le M\}$ ensures that the parameter space is compact,
which is typically required to upgrade pointwise convergence to uniform convergence over $\Sigma$.
In Section~\ref{sec:uniform-convergence}, we show how this restriction can be removed and uniform convergence extended
to the entire positive semidefinite cone $\Sym_+^d$ when $X$ has a continuous distribution.

\newcommand{\y}{\mathbf{y}}
\newcommand{\K}{\mathbf{K}}
\subsection{Illustration} 
\label{sec:illustration-under-finite-atomic-measures}
\indent\indent
This section provides intuition for the preceding results and illustrates the associated formulas.
In particular, we will demonstrate the reparameterization from $U$ to $\Sigma$, the first-variation
formulas $\mathrm{D}J(\Sigma,\lambda)[W]$, and the gradients $\nabla J(\Sigma,\lambda)$ through
an example in which $\P$ is a uniform measure on finite atoms.
Readers already familiar with these results may skip Section~\ref{sec:illustration-under-finite-atomic-measures}.

\subsubsection{Finite Atomic Measure Setup}
Consider the case where $\P$ denotes a uniform distribution over $m$ atoms $\{(x_i, y_i)\}_{i=1}^m$, namely, 
\begin{equation*}
	\P((X, Y) = (x_i, y_i)) = 1/m,~~~\forall i = 1, 2, \ldots, m.
\end{equation*}
By the representer theorem~\cite[Proposition 8]{CuckerSm02}, the objective under this measure $\P$ admits the formula
\begin{equation*}
	\mathcal{J}(U, \lambda) = \lambda \y^\top (\K_U+ m\lambda I_m)^{-1} \y
\end{equation*}
where $\y = (y_1,\ldots,y_m)^\top \in \R^m$,
$\K_U \in \R^{m\times m}$ has entries $(\K_U)_{ij} = \kernel(Ux_i, Ux_j)$,
and $I_m \in \R^{m\times m}$ denotes the identity matrix. The minimizer $f_{U, \lambda}(Ux_i)$ 
can be computed
and satisfies 
\begin{equation*}
	\bigl(
		f_{U,\lambda}(Ux_1),
		\ldots,
		f_{U,\lambda}(Ux_m)
	\bigr)^\top
	=
	\K_U\,
	\bigl(\K_U + m\lambda I_m\bigr)^{-1}
	\y .
\end{equation*}

\subsubsection{Rotation Invariance and Reparameterization} 
Given that $\kernel(x,x')=\psi(\norm{x-x'}^2)$, the kernel matrix $\K_U$ depends 
on $U$ only through $\Sigma = U^\top U$. Consequently, the objective depends on $U$ only through 
$\Sigma$ and can be written as
\begin{equation}	
\label{eqn:J-formula-finite-atoms}
	J(\Sigma, \lambda) =
	\lambda\,
	\y^\top
	\bigl(\K_\Sigma + m\lambda I_m\bigr)^{-1}
	\y,
\end{equation}
where $\K_\Sigma  \in \R^{m\times m}$ has entries $(\K_\Sigma)_{ij} = \psi\!\bigl(\norm{x_i-x_j}_\Sigma^2\bigr)$. 
Similarly, we find $f_{U, \lambda}(Ux_i)$ depends on $U$ only through $\Sigma$. The residual 
$r_\Sigma(x, y) = Y - f_{U, \lambda}(Ux)$ satisfies 
\begin{equation}
\label{eqn:r-sigma-finite-m}
	\bigl(
		r_{\Sigma,\lambda}(x_1, y_1),
		\ldots,
		r_{\Sigma,\lambda}(x_m, y_m)
	\bigr)^\top
	= m \lambda
	\bigl(\K_\Sigma + m\lambda I_m\bigr)^{-1}
	\y .
\end{equation}

\subsubsection{First Variation} 
Let $W \in T_{\Sym_+^d}(\Sigma)$. If $\psi$ is differentiable, then
$\Sigma \mapsto \K_\Sigma$ is directionally differentiable along $W$, with directional derivative
denoted by 
$
\mathrm{D}\K_\Sigma[W]
:=
\lim_{s\to 0^+}\frac{1}{s}\bigl(\K_{\Sigma+sW}-\K_\Sigma\bigr)$. 
Moreover, the derivative is given entrywise by
\begin{equation}
\label{eqn:D-K-Sigma-W}
	\bigl(\mathrm{D}\K_\Sigma[W]\bigr)_{ij}
	=
	\psi'\!\bigl(\|x_i -x_j\|_\Sigma^2\bigr)\,
	(x_i -x_j)^\top W (x_i -x_j),
	\qquad i,j=1,\ldots,m.
\end{equation}
As a result, $J$ has directional derivative, and by the chain rule, it satisfies 
\begin{equation} 
\label{eqn:direction-derivative-finite-m}
\mathrm{D}J(\Sigma,\lambda)[W] = -\lambda\, \y^\top (\K_\Sigma+m\lambda I_m)^{-1}\, \mathrm{D}\K_\Sigma[W]\, (\K_\Sigma+m\lambda I_m)^{-1}\y . 
\end{equation}
Substituting the formula from~\eqref{eqn:r-sigma-finite-m} and~\eqref{eqn:D-K-Sigma-W} 
into~\eqref{eqn:direction-derivative-finite-m}, we get 
\begin{equation}
\label{eqn:directional-derivative-formula-finite-m}
	\mathrm{D}J(\Sigma,\lambda)[W] =
	-\frac{1}{m^2\lambda}
	\sum_{i,j=1}^m
	r_{\Sigma,\lambda}(x_i,y_i)\,
	r_{\Sigma,\lambda}(x_j,y_j)\,
	\psi'\!\bigl(\|x_i -x_j\|_\Sigma^2\bigr)\,
	(x_i -x_j)^\top W (x_i -x_j).
\end{equation}
Note it admits the representation
$
	\mathrm{D}J(\Sigma,\lambda)[W]
	= \bigl\langle \nabla J(\Sigma,\lambda),\, W \bigr\rangle,
$
where $\nabla J(\Sigma,\lambda)$ is given by
\begin{equation}
\label{eqn:gradient-formula-finite-m}
	\nabla J(\Sigma,\lambda)
	=
	-\frac{1}{m^2\lambda}
	\sum_{i,j=1}^m
	r_{\Sigma,\lambda}(x_i,y_i)\,
	r_{\Sigma,\lambda}(x_j,y_j)\,
	\psi'\!\bigl(\|x_i -x_j\|_\Sigma^2\bigr)\,
	(x_i -x_j) (x_i -x_j)^\top.
\end{equation} 

Now, if we take independent copies $(X,Y)$ and $(X',Y')$ drawn from $\P$, then
\[
\P\bigl((X,Y),(X',Y') = (x_i,y_i),(x_j,y_j)\bigr)= 1/m^2,~~~\forall i, j = 1, 2 \ldots, m.
\]
Thus, the directional derivative and gradient formula~\eqref{eqn:directional-derivative-formula-finite-m} 
and~\eqref{eqn:gradient-formula-finite-m} can be written equivalently as
\begin{equation*}
\begin{split}
	\mathrm{D}J(\Sigma,\lambda)[W]
	&=
	-\frac{1}{\lambda}\,
	\E\!\left[
		r_{\Sigma,\lambda}(X,Y)\,
		r_{\Sigma,\lambda}(X',Y')\,
		\psi'\!\bigl(\|X-X'\|_\Sigma^2\bigr)\,
		(X-X')^\top W (X-X')
	\right] \\
	\nabla J(\Sigma,\lambda)
	&=
	-\frac{1}{\lambda}\,
	\E\!\left[
		r_{\Sigma,\lambda}(X,Y)\,
		r_{\Sigma,\lambda}(X',Y')\,
		\psi'\!\bigl(\|X-X'\|_\Sigma^2\bigr)\,
		(X-X') (X-X')^\top
	\right]
\end{split} .
\end{equation*}
These formula agrees with the first-variational formula in~\eqref{eqn:first-variation-formula}
and gradient formula in~\eqref{eqn:gradient-of-J}.

\section{Population Minimizers}
\label{sec:global-minimizers}

\subsection{Existence} 
\indent\indent
We first show that when $X$ admits a continuous distribution, the function $J$ attains its minimum. 
This rules out degenerate situations where minimizers ``escape to infinity". 

\vspace{.5em} 
\begin{theorem}
\label{theorem:existence-of-global-minimizers}
Let Assumption~\ref{assumption:X-continuous} hold. Then, for every $\lambda > 0$,
\begin{equation*}
	\lim_{\Sigma \to \infty} J(\Sigma, \lambda) = \E[Y^2]. 
\end{equation*}
Consequently, $J$ admits a minimizer $\Sigma_* \in \Sym_+^d$ satisfying $J(\Sigma_*, \lambda) = \inf_{\Sigma \in \Sym_+^d} J(\Sigma, \lambda)$.
\end{theorem}

\begin{proof} 
Choose $U$ such that $\Sigma = U^\top U$. By definition, $J(\Sigma,\lambda) = \mathcal{J}(U,\lambda)$. 

Let $(X', Y')$ denote an independent copy of $(X, Y)$. By definition, $K(UX, UX') = \psi(\norm{X-X'}_\Sigma^2)$. 
By Lemma~\ref{lemma:bounds-on-Sigma}, we have the upper and lower bounds
\begin{equation*}
	\E[Y^2] \ge J(\Sigma, \lambda) \ge \E[Y^2] - \lambda^{-1} \E[YY' \psi(\norm{X- X'}_\Sigma^2)],~~~\forall \Sigma \in \Sym_+^d.
\end{equation*}
Since $X$ is continuous, 
$\mathbb{P}(X \neq X') = 1$. Because $\psi(r) \to 0$ as $r \to \infty$, we have
$\psi(\|X - X'\|_\Sigma^2) \to 0$ almost surely as $\Sigma \to\infty$. 
Moreover, since $\psi$ is uniformly bounded, Lebesgue’s dominated convergence theorem yields
$\E[YY' \psi(\norm{X- X'}_\Sigma^2)] \to 0$ as $\Sigma \to \infty$. 
This implies $\lim_{\Sigma \to \infty} J(\Sigma, \lambda) = \E[Y^2]$.

Since $J$ is continuous and bounded above by $\E[Y^2]$, 
it must attain its infimum over $\Sym_+^d$.
\end{proof} 


In Section~\ref{sec:discrete-counterexample}, we illustrate that the continuity assumption on $X$ in
Theorem~\ref{theorem:existence-of-global-minimizers} is necessary by exhibiting
a discrete distribution $\P$ for which the objective $J(\cdot,\lambda)$ does not attain
its infimum.

\subsection{Subspace Recovery} 
\indent \indent
We then study recovery of the true subspace $S$ at the population level.
It turns out for all $\lambda > 0$, the minimizers at the population level span subspaces contained in $S$, 
and for all small enough $\lambda$, they \emph{exactly} recover the true predictive subspace.
 
\vspace{.5em} 
\begin{theorem}
\label{thm:population-recovery}
Let Assumptions~\ref{assumption:X-continuous}--~\ref{assumption:symmetry-of-RKHS} 
hold. \vspace{.5em}
\begin{itemize}
\item[$\mathsf{(a)}$] 
For every $\lambda > 0$, any minimizer $\Sigma_*$ of $J(\cdot, \lambda)$ satisfies
\begin{equation*}
	\Range(\Sigma_*) \subset S.
\end{equation*}
\item[$\mathsf{(b)}$] 
There exists $\lambda_0 > 0$ such that for all $0 < \lambda \le \lambda_0$,
any minimizer $\Sigma_*$ of $J(\cdot, \lambda)$ satisfies 
\begin{equation*}
	\Range(\Sigma_*) = S.
\end{equation*} 
\end{itemize}
\end{theorem} 
The $U$-formulation, Theorem~\ref{theorem:population-result-intro}, follows directly from the reparameterization $\Sigma = U^\top U$.

\subsubsection{Proof of Subspace Containment} 

\indent \indent
We first prove part $\mathsf{(a)}$ in Theorem~\ref{thm:population-recovery}. 
The main idea of our analysis rests on a simple observation 
(Lemma~\ref{lemma:projection-reduce-function-value}): discarding directions 
outside $S$ can only improve (and, in many cases, strictly improve) the objective value. 
Consequently, any population minimizer cannot have a nontrivial component in $S^\perp$, 
so its range must be contained in $S$.

\vspace{0.5em}
\begin{lemma}
\label{lemma:projection-reduce-function-value}
For every $\lambda > 0$
\begin{equation*}
	J(\Sigma, \lambda) \ge J(\Pi_S \Sigma \Pi_S, \lambda),~~~\forall \Sigma \in \Sym_+^d.
\end{equation*} 
The inequality is strict when $\Range(\Sigma) \not \subset S$ and $J(\Sigma, \lambda) < \E[Y^2]$.
\end{lemma} 

\begin{proof} 
Fix $U$ such that $\Sigma = U^\top U$. 
For every $\xi$ in $\R^d$, we define the translated function $\tau_\xi f$ by
\begin{equation*}
	(\tau_\xi f)(x):= f(x+\xi),~~~\forall x \in \R^d.
\end{equation*}
By the translation-invariance of the RKHS $H$, we have $\norm{\tau_\xi f}_H = \norm{f}_H$. 
By the minimizing property of $\mathcal{J}$, we deduce for every vector $\xi$ in $\R^d$, and $f \in H$ that 
\begin{equation*}
\begin{aligned}
  \E[(Y - f(U\Pi_S X + \xi))^2] + \lambda \|f\|_H^2
   &= \E[(Y - (\tau_\xi f)(U\Pi_S X))^2] + \lambda \|\tau_\xi f\|_H^2 \\
   &= \mathcal{I}(\tau_\xi f, U\Pi_S, \lambda)
   \;\ge\; \mathcal{J}(U\Pi_S, \lambda).
\end{aligned}
\end{equation*}
Define $c := \E[(Y - \E[Y|X])^2]$. Since $\E[Y|X] = \E[Y|\Pi_S X]$, we can decompose the squared error as
\[
\E[(Y - f(U \Pi_{S} X + \xi))^2] = \E\left[(\E[Y|\Pi_{S} X] - f(U \Pi_{S} X + \xi))^2\right] + c.
\]
Substituting into the previous inequality gives
\[
\E\left[(\E[Y|\Pi_{S} X] - f(U \Pi_{S} X + \xi))^2\right] + \lambda \|f\|_H^2 + c \ge \mathcal{J}(U\Pi_S, \lambda).
\]

By Assumption~\ref{assumption:independence}, $\Pi_{S} X$ and $\Pi_{S^\perp} X$ are independent. Setting $\xi := U \Pi_{S^\perp} X$, we get
\[
  \E\left[(\E[Y|\Pi_{S} X] - f(U \Pi_{S} X + U\Pi_{S^\perp} X))^2 \mid  \Pi_{S^\perp} X\right] + \lambda \|f\|_H^2 + c \ge \mathcal{J}(U\Pi_S, \lambda).
\]
Taking expectation on both sides, and using $I = \Pi_{S} + \Pi_{S^\perp}$, we get 
\begin{equation*}
\begin{split} 
	\mathcal{I}(f, U, \lambda) 
		&=  \E[(Y - f( U X))^2] + \lambda \norm{f}_H^2  \\
		&= \E[(Y - f( U \Pi_{S} X + U \Pi_{S^\perp} X))^2] + \lambda \norm{f}_H^2 \\
		&=   \E\left[(\E[Y|\Pi_{S} X] - f(U \Pi_{S} X + U \Pi_{S^\perp} X))^2\right] + \lambda \|f\|_H^2 + c 
			\ge  \mathcal{J}(U\Pi_S, \lambda).
\end{split} 
\end{equation*}
By minimizing over $f \in H$, we get 
$\mathcal{J}(U, \lambda)= \min_f \mathcal{I}(f, U, \lambda)  \ge \mathcal{J}(U\Pi_S, \lambda)$. 
Since $\Sigma = U^\top U$, 
\begin{equation*} J(\Sigma,\lambda) = \mathcal{J}(U,\lambda) \ge \mathcal{J}(\Pi_S U,\lambda)
  = J(\Pi_S \Sigma \Pi_S,\lambda).
 \end{equation*} 

Our derivation shows equality 
$J(\Sigma, \lambda) = J(\Pi_{S} \Sigma \Pi_S, \lambda)$ holds
only if the minimizer $f_{U, \lambda}$ obeys 
\begin{equation*}
(\tau_{\xi}) f_{U, \lambda} = (\tau_{\xi'}) f_{U, \lambda}~~\text{for almost all}~\xi= U\Pi_{S^\perp} X, \, \xi' = U\Pi_{S^\perp}X'
\end{equation*}
where $X, X'$ are independent copies. If $\Pi_{S^\perp}\Sigma\Pi_{S^\perp}\neq 0$, then $U\Pi_{S^\perp}\neq 0$, and thus $\operatorname{Cov}(U\Pi_{S^\perp}X)\neq 0$. Therefore, with positive probability $\xi\neq \xi'$. Fix such a realization and set $z=\xi-\xi'\neq 0$. Then
\begin{equation*}
	f_{U, \lambda}(x + z)  = f_{U, \lambda}(x)~~~\forall x\in \R^d.
\end{equation*}
However, since $f_{U, \lambda} \in H$, it satisfies $\lim_{x \to \infty} f_{U, \lambda}(x) = 0$ 
by Lemma~\ref{lemma:continuous-embedding}. In turn, this implies that $f_{U, \lambda}$ must be identically zero. 
Consequently, $\mathcal{J}(U, \lambda) = \E[Y^2]$, and therefore $J(\Sigma, \lambda) = \E[Y^2]$.
\end{proof} 

%

\vspace{.5em} 

Let $\Sigma_*$ denote a 
minimizer of $J(\cdot, \lambda)$. Lemma~\ref{lemma:projection-reduce-function-value}
guarantees 
\[
J(\Sigma_*,\lambda)\ge J(\Pi_S\Sigma_*\Pi_S,\lambda),
\]
with strict inequality whenever
$\Range(\Sigma_*)\not\subset S$ and
$J(\Sigma_*,\lambda)<\mathbb{E}[Y^2]$. On the other hand, 
$J(\Sigma_*,\lambda) \le J(\Pi_S \Sigma_* \Pi_S,\lambda)$ by optimality of $\Sigma_*$. 
Hence equality must hold, and thereby, it must be either
\[
\Range(\Sigma_*)\subset S
\qquad\text{or}\qquad
J(\Sigma_*,\lambda)=\mathbb{E}[Y^2].
\]
By Assumption~\ref{assumption:non-degeneracy}, $\E[Y|X] \neq 0$ with positive probability. 
Lemma~\ref{lemma:upper-bound-U} then implies $\mathcal{J}(I,\lambda) < \E[Y^2]$ for the identity matrix $I$. Thereby, 
$J(I,\lambda)=\mathcal{J}(I,\lambda)$, and by optimality, $J(\Sigma_*, \lambda) < \E[Y^2]$. 

As a result, $\Range(\Sigma_*)\subset S$ must hold. This completes the proof.

\subsubsection{Proof of Exact Recovery} 

\indent\indent
We then prove part $\mathsf{(b)}$ of Theorem~\ref{thm:population-recovery}. By part $\mathsf{(a)}$, already established in the previous subsection, 
any global minimizer must satisfy $\Range(\Sigma_*) \subset S$. Thus it remains to determine whether $\Range(\Sigma_*)$ 
can be a proper subspace of $S$. 

To do so, we exploit the fact that $S$ is the minimal subspace such that 
$\E[Y | X] = \E[Y | \Pi_S X]$. 
In particular, projecting $X$ onto any strict subspace $T \subsetneq S$ necessarily increases the prediction error, and the next lemma shows that this increase is uniformly bounded away from zero. The proof is deferred to Appendix~\ref{section:proof-uniform-gap}.
\vspace{.5em} 
\begin{lemma}
\label{lem:uniform-gap}
\begin{equation*}
	\inf_{T:  T\subsetneq S} \E\bigl[(Y - \E[Y | \Pi_T X])^2\bigr] >  \E\bigl[(Y - \E[Y | X])^2\bigr]. 
\end{equation*} 
\end{lemma}

\vspace{.5em}

Next, we relate this to the objective $J$. By the optimality of conditional expectation,  
\begin{equation*}
\mathcal{I}(f, U, \lambda) \ge \E[(Y- f(UX))^2] \ge \E[(Y - \E[Y|UX])^2]
\end{equation*}
and therefore, after minimizing over $f$, we get 
\begin{equation*}
	\mathcal{J}(U, \lambda) \ge \E[(Y - \E[Y|UX])^2] = \E[(Y - \E[Y|\Pi_{\Range(U^\top)}X])^2].
\end{equation*}
Passing this inequality to the $\Sigma$–parameterization via $\Sigma = U^\top U$ gives for all $\lambda > 0$,
\begin{equation*}
	J(\Sigma, \lambda) \ge \E[(Y - \E[Y|\Pi_{\Range(\Sigma)}X])^2],~~~\forall \Sigma \in \Sym_+^d.
\end{equation*} 
In particular, Lemma~\ref{lem:uniform-gap} implies the existence of $\varepsilon>0$ such that
\begin{equation*}
	\inf \{J(\Sigma, \lambda) \mid \Range(\Sigma) \subsetneq S \} \ge  \E[(Y - \E[Y|X])^2] + \eps. 
\end{equation*}  

We now contrast this with what is achievable when $\lambda$ is small. Since $H$ is dense in $L_2(\P_X)$ (here
$\P_X$ is the marginal distribution of $X$), there exists $f_\varepsilon\in H$ with
\begin{equation*}
	\E[(Y-f_\varepsilon (X))]^2 < \E[(Y - \E[Y|X])^2] + \eps/2.
\end{equation*}
Letting $I$ denote the identity matrix, we obtain for sufficiently small $\lambda>0$:
\begin{equation*}
	J(I, \lambda) = \mathcal{J}(I, \lambda) \le 
		\mathcal{I} (f_\varepsilon, I, \lambda) = \E[(Y-f_\varepsilon (X))]^2 + \lambda \norm{f_\varepsilon}_H^2 < 
		\E[(Y - \E[Y|X])^2] + \eps.
\end{equation*}
Thus
\begin{equation*}
	\inf \{J(\Sigma, \lambda) \mid \Range(\Sigma) \subsetneq S \} >  J(I, \lambda) \ge 
		 \inf \{J(\Sigma, \lambda) \mid \Sigma \in \Sym_+^{d} \}.
\end{equation*} 
Consequently, for all sufficiently small $\lambda > 0$, no minimizer $\Sigma_*$ of $J(\cdot,\lambda)$ can satisfy $\Range(\Sigma_*) \subsetneq S$. Since Part $\mathsf{(a)}$ (proved in the previous section) already implies that any global minimizer must obey $\Range(\Sigma_*) \subset S$, the only possibility is $\Range(\Sigma_*) = S$.

\newcommand{\proj}{\Pi}
\newcommand{\tangent}{\mathcal{T}}
\newcommand{\tr}{{tr}}
\section{Sharpness Property}
\label{sec:sharpness-property}
\indent\indent
In variational analysis, sharpness is a key mechanism behind stability and strong identification of solution sets 
under perturbations~\cite{RockafellarTyWe98, Lewis02}. In our setting, sharpness originates from a first-order variational property of the objective $J$.
Specifically, Theorem~\ref{theorem:sharpness-general} shows that at any $\Sigma$ whose range is
 in $S$, any perturbation $W$ that introduces components in $S^\perp$ cannot decrease 
 the objective at first order, and strictly increases it whenever $J(\Sigma, \lambda) < \E[Y^2]$.

\vspace{.5em} 
\begin{theorem}
\label{theorem:sharpness-general} 
Let Assumptions~\ref{assumption:X-continuous}--~\ref{assumption:symmetry-of-RKHS} 
hold. Let $\lambda > 0$. 
For every $\Sigma \in \Sym_+^d$ with $\Range(\Sigma) \subset S$, and $W \in T_{\Sym_+^d}(\Sigma)$ with $\Range(W) \subset S^\perp$, 
\begin{equation*}
	\mathrm{D} J(\Sigma, \lambda)[W] \ge 0.
\end{equation*}
The inequality is strict whenever $W \neq 0$ and $J(\Sigma, \lambda) < \E[Y^2]$. 
\end{theorem} 

\vspace{.5em} 
An important consequence of Theorem~\ref{theorem:sharpness-general} is 
Theorem~\ref{theorem:sharpness}, which shows that the objective 
$J$ is \emph{sharp} at each of its minimizers. This sharpness property means that 
at any minimizer, every perturbation orthogonal to $S$ strictly increases the objective, and does so at a \emph{linear} rate.

\vspace{.5em} 
\begin{theorem}
\label{theorem:sharpness}
Let Assumptions~\ref{assumption:X-continuous}--~\ref{assumption:symmetry-of-RKHS} 
hold. Let $\lambda > 0$. For every minimizer $\Sigma_*$ of $J(\cdot,\lambda)$,
there exists $\rho > 0$ such that
\begin{equation*}
	\mathrm{D} J(\Sigma_*, \lambda)[W] \ge \rho \norm{W}, \qquad \forall W \in T_{\Sym_+^d}(\Sigma_*)~\text{with}~ \Range(W) \subset S^\perp.
\end{equation*}
\end{theorem} 

Theorem~\ref{theorem:sharpness} underlies the remainder of the paper, as the sharpness of $J$ provides the key 
mechanism that preserves the low-rankness of minimizers when passing to the empirical objective $J_n$.



\subsection{Proof of Theorem~\ref{theorem:sharpness-general}}
\label{sec:proof-of-theorem-sharpness-general}

\indent\indent
Let $\Sigma$ satisfy $\Range(\Sigma) \subset S$.
We first characterize matrices $W \in T_{\Sym_+^d}(\Sigma)$ whose range lies in $S^\perp$.
Lemma~\ref{lemma:basic-convex-geometry} follows from basic convex analysis; its proof is in Appendix~\ref{sec:proof-of-lemma-basic-convex-geometry}.

\vspace{.5em}
\begin{lemma}
\label{lemma:basic-convex-geometry}
Suppose $\Sigma$ satisfies $\Range(\Sigma) \subset S$. Then 
\begin{equation*}
	\left\{W: W \in T_{\Sym_+^d}(\Sigma) \, , \, \Range(W)\subset S^\perp\right\}
		= \left\{W: W \in \Sym_+^d \, , \, \Range(W)\subset S^\perp\right\}.
\end{equation*}
\end{lemma} 
By Lemma~\ref{lemma:basic-convex-geometry}, any $W\in T_{\Sym_+^d}(\Sigma)$ with $\Range(W)\subset S^\perp$
must be positive semidefinite. 
By the spectral theorem, there exist eigenvalues $\lambda_i\ge 0$ and orthonormal $\{v_i\} \subset S^\perp$ such that
\begin{equation}
  W = \sum_i \lambda_i v_i v_i^\top.
\end{equation} 
By linearity of $\mathrm{D} J(\Sigma,\lambda)[\cdot]$,
\begin{equation}
\label{eqn:W-into-rank-one-perturbations}
  \mathrm{D} J(\Sigma,\lambda)[W]
  = \sum_i \,\lambda_i \mathrm{D} J(\Sigma,\lambda)[v_i v_i^\top].
\end{equation}
This motivates us to compute the directional derivative $\mathrm{D} J(\Sigma, \lambda)[vv^\top]$ with $v \in S^\perp$.
\vspace{.5em} 
\begin{lemma}
\label{lemma:factorization-lemma}
Suppose $\Sigma$ satisfies $\Range(\Sigma)\subset S$ and $v\in S^\perp$. Then 
\begin{equation}
\label{eqn:factorization-formula}
	\mathrm{D} J(\Sigma, \lambda)[vv^\top] = -
		\frac{1}{\lambda} \E\bigl[
      r_{\Sigma,\lambda}(X, Y)\,
      r_{\Sigma,\lambda}(X', Y')\,
      \psi'\bigl(\|X-X'\|_{\Sigma}^2\bigr)
    \bigr]\, \cdot \,
    \E\bigl[(v^\top (X-X'))^2\bigr].
\end{equation} 
\end{lemma} 

\begin{proof} 
By Lemma~\ref{theorem:first-variation-formula}, 
\begin{equation*}
	\mathrm{D} J(\Sigma, \lambda)[vv^\top]  =  -\frac{1}{\lambda} 
		\E \left[r_{\Sigma, \lambda}(X, Y) r_{\Sigma, \lambda}(X', Y')  \psi^\prime(\norm{X-X'}_{\Sigma}^2)\, \cdot \, (v^\top (X-X'))^2\right].
\end{equation*} 
Pick $U$ with $\Sigma = U^\top U$. Then $r_{\Sigma, \lambda}(X, Y) = Y - f_{U, \lambda}(U X)$. 
By a direct computation,
\begin{equation*}
\begin{split} 
	&\E[r_{\Sigma, \lambda}(X, Y) r_{\Sigma, \lambda}(X', Y')  \psi^\prime(\norm{X-X'}_{\Sigma}^2)\mid X, X'] \\
	&= (\E[Y|X] - f_{U, \lambda}(U X)) (\E[Y'|X'] - f_{U, \lambda}(U X'))  \psi^\prime(\norm{U (X-X')}^2).
\end{split} 
\end{equation*} 
Since $\ker(U) = \Range(\Sigma)^\perp \supset S^\perp$, we have $UX = U\Pi_S X$. 
Moreover, $\E[Y | X] = \E[Y | \Pi_S X]$.  On the other hand, since $v \in S^\perp$, we have 
$v^\top X = v^\top\Pi_{S^\perp}X$. By Assumption~\ref{assumption:independence},  $\Pi_S X$ and $\Pi_{S^\perp} X$ are independent. 
It follows that $UX$ and $\E[Y | X]$ are independent of $v^\top X$ (and
likewise for $X'$). Thus
\begin{equation*}
	\E[r_{\Sigma, \lambda}(X, Y) r_{\Sigma, \lambda}(X', Y')  \psi^\prime(\norm{X-X'}_{\Sigma}^2)| X, X']
	~~\text{and}~~(v^\top(X-X'))^2
\end{equation*} 
are independent. Thereby, the expectation of their products equals the products of their expectations: 
\begin{equation*}
\begin{split}  
		& \E \left[r_{\Sigma, \lambda}(X, Y) r_{\Sigma, \lambda}(X', Y')  \psi^\prime(\norm{X-X'}_{\Sigma}^2)\, \cdot \, (v^\top (X-X'))^2\right] \\
		&= \E \left[r_{\Sigma, \lambda}(X, Y) r_{\Sigma, \lambda}(X', Y')  \psi^\prime(\norm{X-X'}_{\Sigma}^2)\right] \, \cdot \, \E\left[ (v^\top (X-X'))^2\right].
\end{split} 
\end{equation*}
This completes the proof. 
\end{proof} 

The second factor in~\eqref{eqn:factorization-formula} is clearly nonnegative: 
\begin{equation}
	\E\bigl[(v^\top (X-X'))^2\bigr] \ge 0.
\end{equation}
Under assumption~\ref{assumption:independence}, it is strictly positive when $0 \neq v \in S^\perp$. We now determine the sign of the first factor in~\eqref{eqn:factorization-formula}.
The key observation is that $(x, x') \mapsto \psi^\prime(\norm{x-x'}_{\Sigma}^2)$ is a negative definite kernel. This property can be established via a Fourier-transform representation, and leads to the following result.

\vspace{.5em} 
\begin{lemma}
\label{lemma: first-factor-negative}
For every $\Sigma$, we have 
 \begin{equation}
 \label{eqn:first-factor-negative}
 	 \E \left[r_{\Sigma, \lambda}(X, Y) r_{\Sigma, \lambda}(X', Y')  \psi^\prime(\norm{X-X'}_{\Sigma}^2)\right ] \le 0.
 \end{equation} 
 The inequality is strict when $J(\Sigma, \lambda) < \E[Y^2]$. 
\end{lemma} 

\begin{proof} 
Pick $U$ such that $\Sigma = U^\top U$ and define
\[
g(z) := \E\big[r_{\Sigma,\lambda}(X,Y)\mid UX = z\big].
\]
Since $\|X-X'\|_{\Sigma}^2 = \|U(X-X')\|^2$, conditioning on $(UX,UX')$ yields
\[
\begin{aligned}
&\E\bigl[
  r_{\Sigma,\lambda}(X,Y)\,r_{\Sigma,\lambda}(X',Y')\,
  \psi'(\|X-X'\|_{\Sigma}^2)
\bigr] 
=
\E\bigl[g(UX)\,g(UX')\,\psi'(\|U(X-X')\|^2)\bigr].
\end{aligned}
\]
By definition, $\psi'(z)=-\int_0^\infty e^{-tz}\,t\,\mu(dt)$. Denote $\zeta_t(\omega):=\frac{1}{(4\pi t)^{d/2}}e^{-\|\omega\|^2/(4t)}$. Then 
\[
\begin{aligned}
\E[g(UX)g(UX')\psi'(\|U(X-X')\|^2)]
&=
-\int_0^\infty
\E\big[g(UX)g(UX')e^{-t\|UX-UX'\|^2}\big]\,
t\,\mu(dt) \\
&=
-\int_0^\infty\,
\E\big[g(UX)g(UX')\int_{\R^d} e^{i\omega^\top(UX-UX')}\zeta_t(\omega) d\omega \big]\,
\,t\,\mu(dt) \\
&=
-\int_0^\infty\!\!\int_{\R^d}
\E\big[g(UX)g(UX')e^{i\omega^\top(UX-UX')}\big]\,
\zeta_t(\omega)\,t\,d\omega\,\mu(dt) \\
&=
-\int_0^\infty\!\!\int_{\R^d}
\bigl|\E\big[g(UX)e^{i\omega^\top UX}\big]\bigr|^2
\zeta_t(\omega)\,t\,d\omega\,\mu(dt),
\end{aligned}
\]
where the second line uses the Fourier representation of Gaussian: 
$e^{-t\|z-z'\|^2}
=
\int
e^{i\omega^\top(z-z')}\zeta_t(\omega)\,d\omega$, and the last line uses that
$UX$ and $UX'$ are independent and identically distributed. This shows
\[
\E\!\left[
r_{\Sigma,\lambda}(X,Y)\,r_{\Sigma,\lambda}(X',Y')\,
\psi'(\|X-X'\|_{\Sigma}^2)
\right]
\le 0.
\]

Finally, if equality holds, then the integrand must vanish. Given that
$\zeta_t(\omega)>0$ for all $t,\omega$, this forces
$\E[g(UX)e^{i\omega^\top UX}]=0$ for almost all $\omega$, and thus
$g(UX)=0$ almost surely. As a result, we get $0 = g(UX) = \E[r_{\Sigma, \lambda}(X, Y)|UX] = 
\E[Y-f_{U,\lambda}(UX)| UX]$. 
Combining these with Lemma~\ref{lemma:euler-lagrange-identity} yields
\begin{equation*}
\lambda\|f_{U,\lambda}\|_H^2
=\E[(Y-f_{U,\lambda}(UX))\,f_{U,\lambda}(UX)] = 0.
\end{equation*}
Thus $f_{U,\lambda}\equiv 0$. Consequently $\mathcal{J}(U,\lambda)=\E[Y^2]$, and hence
$J(\Sigma,\lambda)=\E[Y^2]$.

\end{proof} 


Combining Lemma~\ref{lemma:factorization-lemma} and~\ref{lemma: first-factor-negative}, we get
for every $\Range(\Sigma)\subset S$ and $v\in S^\perp$, 
\begin{equation}
\label{eqn:strict-complementary-slackness-condition}
	\mathrm{D} J(\Sigma, \lambda)[vv^\top] \ge 0
\end{equation} 
Moreover, the inequality is strict whenever $v \neq 0$ and $J(\Sigma, \lambda) < \E[Y^2]$. Following~\eqref{eqn:W-into-rank-one-perturbations}, 
we get 
\begin{equation}
\label{eqn:inequality-for-W}
	  \mathrm{D} J(\Sigma,\lambda)[W]
  = \sum_i \,\lambda_i \mathrm{D} J(\Sigma,\lambda)[v_i v_i^\top] \ge 0.
\end{equation}
If in addition $W \neq 0$ and $J(\Sigma, \lambda) < \E[Y^2]$, then at least one eigenvalue satisfies $\lambda_i > 0$, 
and for such an index $i$ we have $\mathrm{D} J(\Sigma,\lambda)[v_i v_i^\top] > 0$. 
Thus the inequality~\eqref{eqn:inequality-for-W} is strict in this case.  
This completes the proof of Theorem~\ref{theorem:sharpness-general}.

\subsection{Proof of Theorem~\ref{theorem:sharpness}}
Fix $\lambda > 0$ and a minimizer $\Sigma_*$ of $J(\cdot,\lambda)$.
By Theorem~\ref{theorem:sharpness-general}, we have 
\begin{equation*}
	\Range(\Sigma_*) \subset S. 
\end{equation*}
By Assumption~\ref{assumption:non-degeneracy}, $\E[Y|X] \neq 0$ with positive probability. 
Lemma~\ref{lemma:upper-bound-U} then implies $\mathcal{J}(I,\lambda) < \E[Y^2]$ for the identity matrix $I$.  Thereby, 
$J(I,\lambda)=\mathcal{J}(I,\lambda) < \E[Y^2]$, and by optimality, 
\begin{equation*}
	J(\Sigma_*, \lambda) < \E[Y^2].
\end{equation*} 

By Theorem~\ref{theorem:sharpness-general}, we deduce that for every $W \in T_{\Sym_+^d}(\Sigma_*)$ with $\Range(W) \subset S^\perp$ and $W \neq 0$, 
\begin{equation}
	\mathrm{D} J(\Sigma_*, \lambda)[W] > 0.
\end{equation}
Since $W \mapsto \mathrm{D}J(\Sigma_*,\lambda)[W]$ is continuous and homogeneous in $W$, this implies the existence of $\rho > 0$ such that 
$\mathrm{D} J(\Sigma_*, \lambda)[W] \ge \rho \norm{W}$ for all such $W$. Theorem~\ref{theorem:sharpness} follows.

\section{Uniform Approximation of Empirical Objective} 
\label{sec:uniform-convergence}
\indent\indent
We establish uniform convergence of $J_n$ to $J$, together with their 
gradients, over the entire domain $\Sym_+^d$, under the assumption that $X$ has continuous distribution.
This strengthens Lemma~\ref{lemma:law-of-large-number-constrained-set}, which yields uniform 
convergence only on bounded parameter sets without assuming continuity.

\vspace{.5em} 
%

\begin{theorem}
\label{theorem:uniform-convergence}
Let Assumption~\ref{assumption:X-continuous} hold. Let $\lambda > 0$. Almost surely as $n \to \infty$: 
\vspace{.5em} 

\begin{itemize}
\item[$\mathsf{(a)}$]
\begin{equation*}
	\sup_{\Sigma \in \Sym_+^d} 
	\bigl|J_n(\Sigma, \lambda) - J(\Sigma, \lambda)\bigr|
	\;\longrightarrow\; 0.
\end{equation*}

\item[$\mathsf{(b)}$]
\begin{equation*}
	\sup_{\Sigma \in \Sym_+^d} 
	\bigl\|
	\nabla J_n(\Sigma, \lambda)
	-
	\nabla J(\Sigma, \lambda)
	\bigr\|
	\;\longrightarrow\; 0.
\end{equation*}
\end{itemize}
\end{theorem}

\subsection{Proof of Uniform Convergence of the Objective}
\indent\indent
In this section, we prove part $\mathsf{(a)}$ of Theorem~\ref{theorem:uniform-convergence}.
By Lemma~\ref{lemma:law-of-large-number-constrained-set}, we already know that for every finite constant $M < \infty$, almost surely as $n \to \infty$
\begin{equation}
	\sup_{\Sigma:\,\|\Sigma\| \le M} 
	\bigl|J(\Sigma, \lambda) - J_n(\Sigma, \lambda)\bigr|
	\;\longrightarrow\; 0.
\end{equation}
Thus the only remaining difficulty lies in ruling out pathological behavior of $J_n$ and $J$
over $\Sigma$ with large norms, which will be handled using the assumption that $X$ has a 
continuous distribution.
Define
\begin{equation}
\label{eqn:U-U-n-definition}
	W(M) :=  \sup_{\Sigma: \norm{\Sigma} \ge M}
	\E\!\left[\psi^2\!\left(\|X - X'\|_{\Sigma}^2\right)\right], \qquad 
	W_n(M) :=  \sup_{\Sigma: \norm{\Sigma} \ge M}
	\E_n\!\left[\psi^2\!\left(\|X - X'\|_{\Sigma}^2\right)\right].
\end{equation}
To control the gap between $J$ and $J_n$ uniformly over large $\Sigma$, we use
\vspace{.5em}

\begin{lemma}
\label{lemma:deterministic-bounds}
With probability one, for all $M < \infty$,
\begin{equation*}
	\sup_{\Sigma: \norm{\Sigma} \ge M} \bigl|J(\Sigma, \lambda) - J_n(\Sigma, \lambda)\bigr| \le 
		\bigl| \E_n[Y^2] - \E[Y^2] \bigr|
	+ 
	\lambda^{-1} \E[Y^2] \cdot \sqrt{W(M)}
	+ 
	\lambda^{-1} \E_n[Y^2]  \cdot \sqrt{W_n(M)}.
\end{equation*}
\end{lemma} 
\begin{proof} 
By Lemma~\ref{lemma:bounds-on-Sigma}, we have upper and lower bounds on $J$: 
\begin{equation*}
	\E[Y^2] \ge J(\Sigma, \lambda) 
	\ge 
	\E[Y^2] 
	- 
	\lambda^{-1} 
	 \E[Y^2]  \sqrt{\E[\psi^2\!\left(\|X - X'\|_{\Sigma}^2\right)]},
	\qquad 
	\forall \Sigma \in \Sym_+^d.
\end{equation*}
Applying the same result to the empirical measure yields the corresponding bounds for $J_n$:
\begin{equation*}
	\E_n[Y^2]  \ge J_n(\Sigma, \lambda) 
	\ge 
	\E_n[Y^2] 
	- 
	\lambda^{-1} 
	 \E_n[Y^2]  \sqrt{\E_n[\psi^2\!\left(\|X - X'\|_{\Sigma}^2\right)]},
	\qquad 
	\forall \Sigma \in \Sym_+^d.
\end{equation*}
Combining these bounds, we get, 
\begin{equation*}
	\bigl|J(\Sigma, \lambda) - J_n(\Sigma, \lambda)\bigr|
	\le 
	\bigl| \E_n[Y^2] - \E[Y^2] \bigr|
	+ 
	\lambda^{-1} \E[Y^2] 
	\sqrt{\E\!\left[\psi^2\!\left(\|X - X'\|_{\Sigma}^2\right)\right]}
	+ 
	\lambda^{-1} \E_n[Y^2] 
	\sqrt{\E_n\!\left[\psi^2\!\left(\|X - X'\|_{\Sigma}^2\right)\right]},~~\forall \Sigma.
\end{equation*}
Taking the supremum over matrices $\Sigma$ with $\|\Sigma\| \ge M$ yields the result. 
\end{proof} 

Here comes a crucial observation. Since $\psi$ is monotone decreasing (as $\psi(z) = \int_0^\infty e^{-tz} \mu(dt)$), 
\begin{equation}
\label{eqn:identities-for-U-U_n}
\begin{split} 
	W(M) = \sup_{\Sigma: \norm{\Sigma} \ge M} \E\!\left[\psi^2\!\left(\|X - X'\|_{\Sigma}^2\right)\right]
		&= \sup_{\Sigma: \norm{\Sigma} = M} \E\!\left[\psi^2\!\left(\|X - X'\|_{\Sigma}^2\right)\right] \\
	W_n(M) = \sup_{\Sigma: \norm{\Sigma} \ge M} \E_n\!\left[\psi^2\!\left(\|X - X'\|_{\Sigma}^2\right)\right]
		&= \sup_{\Sigma: \norm{\Sigma} = M} \E_n\!\left[\psi^2\!\left(\|X - X'\|_{\Sigma}^2\right)\right]
\end{split} .
\end{equation}
Thus the supremum is attained over the \emph{compact} set $\{\Sigma : \|\Sigma\| = M\}$.
This reduction allows us to establish the almost sure convergence $W_n(M) \to W(M)$ as $n \to \infty$.

\vspace{.5em} 
\begin{lemma}
\label{lemma:convergence-U-n-U}
For every $M < \infty$, we have almost surely,
\begin{equation*}
	\lim_{n \to \infty} W_n(M) = W(M).
\end{equation*}
\end{lemma} 
\begin{proof}
By the triangle inequality applied to~\eqref{eqn:identities-for-U-U_n}, it suffices to show that
\begin{equation}
\label{eqn:uniform-convergence-psi-2}
\lim_{n \to \infty} \sup_{\Sigma: \norm{\Sigma} = M} |\E_n[\psi^2\!\left(\|X - X'\|_{\Sigma}^2\right)] - \E[\psi^2\!\left(\|X - X'\|_{\Sigma}^2\right)]| = 0. 
\end{equation}
Let $\eps> 0$ be a small number. The set $\{\Sigma: \norm{\Sigma} = M\}$ is compact, so there exists a finite $\eps$-net 
$\Sigma_1, \Sigma_2, \ldots, \Sigma_N$ such that any $\Sigma$ in this set is $\eps$-close to some $\Sigma_i$. By the law of 
large number applied to $\psi^2\!\left(\|X - X'\|_{\Sigma_i}^2\right)$, we have almost surely, 
\begin{equation}
\label{eqn:pointwise-limit}
	\lim_{n \to \infty} |\E_n[\psi^2\!\left(\|X - X'\|_{\Sigma_i}^2\right)] - \E[\psi^2\!\left(\|X - X'\|_{\Sigma_i}^2\right)]| = 0,~~~\forall i = 1, 2, \ldots.
\end{equation}
Note that $z \mapsto \psi^2(z)$ is Lipschitz, so there is a finite (non-random) $L < \infty$ such that,
\begin{equation*}
	|\psi^2\!\left(\|v\|_{\Sigma}^2\right)  - \psi^2\!\left(\|v\|_{\Sigma'}^2\right) | \le L \norm{\Sigma- \Sigma'} \norm{v}^2~~~\forall v, \Sigma, \Sigma'.
\end{equation*}
This implies, almost surely, over all choices of $\Sigma$, $\Sigma_i$, 
\begin{equation}
\label{eqn:equicontinuity}
\begin{split} 
	|\E_n[\psi^2\!\left(\|X-X'\|_{\Sigma}^2\right)] - \E_n[\psi^2\!\left(\|X-X'\|_{\Sigma_i}^2\right)]| &\le L \norm{\Sigma- \Sigma_i} \E_n[\norm{X-X'}^2] \\
	|\E_n[\psi^2\!\left(\|X-X'\|_{\Sigma}^2\right)] - \E_n[\psi^2\!\left(\|X-X'\|_{\Sigma_i}^2\right)]| &\le L \norm{\Sigma- \Sigma_i} \E_n[\norm{X-X'}^2].
\end{split} 
\end{equation} 
Since $\E_n[\norm{X-X'}^2] \to \E[\norm{X-X'}^2] < \infty$ as $n \to \infty$, and $\eps > 0$ is arbitrary, with~\eqref{eqn:pointwise-limit} 
and~\eqref{eqn:equicontinuity} we obtain
the uniform law of large number~\eqref{eqn:uniform-convergence-psi-2}.  This completes the proof of the lemma. 
\end{proof}

\vspace{.5em} 

Finally, using the assumption that $X$ is continuous, we get

\vspace{.5em} 
\begin{lemma}
\label{lemma:limit-U}
\begin{equation*}
	\lim_{M \to \infty} W(M) = 0.
\end{equation*}  
\end{lemma} 
\begin{proof} 
Since $X$ has a continuous distribution, $\P(X \neq X') = 1$. Because $\psi(r) \to 0$ as $r \to \infty$, we have
$\psi(\|X - X'\|_\Sigma^2) \to 0$ almost surely as $\Sigma \to\infty$. 
Moreover, since $\psi$ is uniformly bounded, Lebesgue’s dominated convergence theorem yields
\begin{equation*}
	\lim_{\Sigma \to \infty} \E\!\left[\psi^2\!\left(\|X - X'\|_{\Sigma}^2\right) \right] = 0.
\end{equation*}  
This implies that $\lim_{M \to \infty} W(M) = 0$. 
\end{proof}

We are ready to finish the proof of Theorem~\ref{theorem:uniform-convergence}$\mathsf{(a)}$. 
We proceed via a standard $\varepsilon$--$\delta$ argument.
Fix an arbitrary $\varepsilon > 0$. By Lemma~\ref{lemma:limit-U}, there exists a finite $M_\varepsilon < \infty$ such that
\[
	W(M_\varepsilon) \le \varepsilon^2.
\]
By Lemma~\ref{lemma:convergence-U-n-U}, almost surely 
\[
	\lim_{n \to \infty} W_n(M_\varepsilon) =  W(M_\varepsilon) \le \varepsilon^2.
\]
By Lemma~\ref{lemma:deterministic-bounds}, together with the law of large number
$\E_n[Y^2] \to \E[Y^2]$, it follows that almost surely
\[
	\limsup_{n\to \infty} \sup_{\Sigma: \norm{\Sigma} \ge M_\eps} \bigl|J(\Sigma, \lambda) - J_n(\Sigma, \lambda)\bigr| 
		\le 2 \lambda^{-1} \E[Y^2] \varepsilon.
\]
On the other hand, by Lemma~\ref{lemma:law-of-large-number-constrained-set}, for this fixed $M_\varepsilon$ we have almost surely
\begin{equation*}
	\lim_{n\to \infty} \sup_{\Sigma:\,\|\Sigma\| \le M_\eps} 
	\bigl|J(\Sigma, \lambda) - J_n(\Sigma, \lambda)\bigr|
	= 0.
\end{equation*}
Combining the two bounds yields almost surely
\begin{equation*}
	\limsup_{n\to \infty} \sup_{\Sigma} \bigl|J(\Sigma, \lambda) - J_n(\Sigma, \lambda)\bigr| \le 2 \lambda^{-1} \E[Y^2] \varepsilon.
\end{equation*}
Take a sequence $\varepsilon_k \downarrow 0$. Since the bound holds for all $k$
with probability one, letting $k \to \infty$ completes the proof.

\subsection{Proof of Uniform Convergence of the Gradient}
\indent\indent
In this section, we prove part $\mathsf{(b)}$ of Theorem~\ref{theorem:uniform-convergence}.
By Lemma~\ref{lemma:law-of-large-number-constrained-set}, we already know that for every finite constant $M< \infty$, almost surely as $n \to \infty$, 
\begin{equation}
	\sup_{\Sigma:\,\|\Sigma\| \le M} 
	\norm{|\grad J(\Sigma, \lambda) - \grad J_n(\Sigma, \lambda)}
	\;\longrightarrow\; 0.
\end{equation}
Thus the only remaining difficulty lies in controlling the behavior of $\grad J_n$ and $\grad J$
over $\Sigma$ with large norms, which will be handled using the assumption that $X$ has a 
continuous distribution.

Define 
\begin{equation}	
\begin{split} 
	V(M) &:= \sup_{\Sigma: \norm{\Sigma} \ge M} \E[ |\psi^\prime|^2(\norm{X-X'}_\Sigma^2) \cdot \norm{(X-X')}^4]  \\
	V_n(M) &:= \sup_{\Sigma: \norm{\Sigma} \ge M} \E[ |\psi^\prime|^2(\norm{X-X'}_\Sigma^2) \cdot \norm{(X-X')}^4]
\end{split} .
\end{equation}

\vspace{.5em}
\begin{lemma}
\label{lemma:deterministic-bounds-V}
With probability one, for all $M < \infty$, 
\begin{equation*}
	\sup_{\Sigma: \norm{\Sigma} \ge M} \norm{\grad J_n(\Sigma, \lambda) - \grad J(\Sigma, \lambda)} \le 
		\frac{1}{\lambda} \E[Y^2] \cdot \sqrt{V(M)} + \frac{1}{\lambda} \E_n[Y^2] \cdot \sqrt{V_n(M)}.
\end{equation*}
\end{lemma}


\begin{proof} 
By definition, 
\begin{equation*}
	\grad J(\Sigma, \lambda) = - \frac{1}{\lambda} \E\left[R_{\Sigma, \lambda}(X, Y) R_{\Sigma, \lambda}(X', Y') \psi^\prime(\norm{X-X'}_\Sigma^2) (X-X') (X-X')^\top\right]. 
\end{equation*} 
Taking the norm and applying the triangle inequality gives
\begin{equation}
\label{eqn:operator-norm-of-grad-J}
	\norm{\grad J(\Sigma, \lambda)} \le \frac{1}{\lambda} \E\left[ |R_{\Sigma, \lambda}(X, Y) R_{\Sigma, \lambda}(X', Y')| \cdot | \psi^\prime(\norm{X-X'}_\Sigma^2)| \cdot \norm{(X-X')}^2\right]. 
\end{equation}
Since $(X',Y')$ is an independent copy of $(X,Y)$, we have 
\begin{equation*}
	\E[|R_{\Sigma, \lambda}(X, Y) R_{\Sigma, \lambda}(X', Y')|^2] = \E[R^2_{\Sigma, \lambda}(X, Y)] \cdot \E[R^2_{\Sigma, \lambda}(X', Y')]. 
		=  \E[R^2_{\Sigma, \lambda}(X, Y)]^2.
\end{equation*}
By Lemma~\ref{lemma:bounds-on-Sigma}, $\E[R_{\Sigma, \lambda}^2(X, Y)] \le J(\Sigma, \lambda) \le \E[Y^2]$ for all $\Sigma$. Thereby, we have 
\begin{equation*}
	\E[|R_{\Sigma, \lambda}(X, Y) R_{\Sigma, \lambda}(X', Y')|^2] \le (\E[Y^2])^2,~~~\forall \Sigma \in {\Sym_+^d}.
\end{equation*}
Applying the Cauchy–Schwarz inequality to the right-hand side of~\eqref{eqn:operator-norm-of-grad-J}, we obtain
\begin{equation*}
	\norm{\grad J(\Sigma, \lambda)} \le \frac{1}{\lambda} \E[Y^2] \cdot \sqrt{\E[ |\psi^\prime|^2(\norm{X-X'}_\Sigma^2) \cdot \norm{(X-X')}^4]}.
\end{equation*}
Taking the supremum over $\Sigma$ satisfying $\norm{\Sigma} \ge M$,  we conclude that for every finite $M < \infty$,
\begin{equation*}
	\sup_{\Sigma: \norm{\Sigma} \ge M} \norm{\grad J(\Sigma, \lambda)} \le
			 \frac{1}{\lambda} \E[Y^2] \cdot \sqrt{V(M)}.
\end{equation*}
By the same argument, replacing $(J, \E, V)$ with $(J_n, \E_n, V_n)$ yields the empirical counterpart
\begin{equation*}
	\sup_{\Sigma: \norm{\Sigma} \ge M} \norm{\grad J_n(\Sigma, \lambda)} \le
			 \frac{1}{\lambda} \E_n[Y^2] \cdot \sqrt{V_n(M)}.
\end{equation*}
The result follows by applying the triangle inequality. 
\end{proof} 

Here comes a crucial observation. Since $|\psi^\prime|$ is monotone decreasing (as $|\psi^\prime|(z) = \int_0^\infty t e^{-tz} \mu(dt)$), 
\begin{equation}
\label{eqn:identities-for-V-V_n}
\begin{split} 
	V(M) &:= \sup_{\Sigma: \norm{\Sigma} = M} \E[ |\psi^\prime|^2(\norm{X-X'}_\Sigma^2) \cdot \norm{(X-X')}^4]  \\
	V_n(M) &:= \sup_{\Sigma: \norm{\Sigma} = M} \E_n[ |\psi^\prime|^2(\norm{X-X'}_\Sigma^2) \cdot \norm{(X-X')}^4]
\end{split} .
\end{equation}
In other words, the supremum is attained over the \emph{compact} set $\{\Sigma : \|\Sigma\| = M\}$.
This reduction allows us to establish the almost sure convergence $V_n(M) \to V(M)$ as $n \to \infty$.

\vspace{.5em} 
\begin{lemma}
\label{lemma:convergence-V-n-V}
For every $M < \infty$, we have almost surely,
\begin{equation*}
	\lim_{n \to \infty} V_n(M) = V(M).
\end{equation*}
\end{lemma} 

\begin{proof} 
The proof follows the same $\varepsilon$-net and uniform law of large numbers argument as in 
Lemma~\ref{lemma:convergence-U-n-U}. 
The only change is to note that $z \mapsto |\psi^\prime|^2(z)$ is Lipschitz on $z\ge 0$, and that the bounds now involve $\norm{X - X'}^6$, 
which remains integrable under our moment assumptions.
\end{proof} 

Finally, using the assumption that $X$ is continuous, we get

\vspace{.5em} 
\begin{lemma}
\label{lemma:limit-V}
\begin{equation*}
	\lim_{M \to \infty} V(M) = 0.
\end{equation*}  
\end{lemma} 
\begin{proof} 
The proof follows the same argument as in Lemma~\ref{lemma:limit-U}.
Since $X$ is continuous, $\P(X \neq X') = 1$. Because $\psi^\prime(r) \to 0$ as $r \to \infty$, we have
$|\psi^\prime|^2(\norm{X-X'}_\Sigma^2) \cdot \norm{(X-X')}^4 \to 0$ almost surely as $\Sigma \to\infty$. 
Since $\psi^\prime$ is uniformly bounded, and $X$ has finite fourth moments, by dominated convergence theorem,
\begin{equation*}
	\lim_{\Sigma \to \infty} \E\!\left[|\psi^\prime|^2(\norm{X-X'}_\Sigma^2) \cdot \norm{(X-X')}^4 \right] = 0.
\end{equation*}  
This implies $\lim_{M \to \infty} V(M) = 0$ as desired. 
\end{proof}

We are ready to finish the proof of Theorem~\ref{theorem:uniform-convergence}$\mathsf{(b)}$ using a standard 
$\varepsilon$-$\delta$ argument.
Fix an arbitrary $\varepsilon > 0$. By Lemma~\ref{lemma:limit-V}, there exists a finite $M_\varepsilon < \infty$ such that
\[
	V(M_\varepsilon) \le \varepsilon^2.
\]
By Lemma~\ref{lemma:convergence-V-n-V}, almost surely 
\[
	\lim_{n \to \infty} V_n(M_\varepsilon) =  V(M_\varepsilon) \le \varepsilon^2.
\]
By Lemma~\ref{lemma:deterministic-bounds-V}, together with the law of large number
$\E_n[Y^2] \to \E[Y^2]$, it follows that almost surely
\[
	\limsup_{n\to \infty} \sup_{\Sigma: \norm{\Sigma} \ge M_\eps} \norm{ \grad J(\Sigma, \lambda) - \grad J_n(\Sigma, \lambda)} 
		\le 2 \lambda^{-1} \E[Y^2] \varepsilon.
\]
On the other hand, by Lemma~\ref{lemma:law-of-large-number-constrained-set}, for this fixed $M_\varepsilon$ we have almost surely
\begin{equation*}
	\lim_{n\to \infty} \sup_{\Sigma:\,\|\Sigma\| \le M_\eps} 
	\norm{ \grad J(\Sigma, \lambda) - \grad J_n(\Sigma, \lambda)} 
	= 0.
\end{equation*}
Combining the two bounds yields almost surely
\begin{equation*}
	\limsup_{n\to \infty} \sup_{\Sigma} \norm{ \grad J(\Sigma, \lambda) - \grad J_n(\Sigma, \lambda)}  \le 2 \lambda^{-1} \E[Y^2] \varepsilon.
\end{equation*}
Take a sequence $\varepsilon_k \downarrow 0$. Since the bound holds for all $\eps_k$, letting $k \to \infty$ completes the proof.

\section{Empirical Minimizers}
\label{sec:empirical-minimizers} 

\subsection{Existence}
Theorem~\ref{theorem:existence-of-global-minimizers-empirical} shows that when $X$ is continuous and $\E[Y|X] \neq 0$, the empirical objective $J_n$ attains a minimizer with high probability. In particular, it ensures that the ``escape to infinity" of minimizers is unlikely to happen, providing a finite-sample counterpart to the population result in Theorem~\ref{theorem:existence-of-global-minimizers}.

\vspace{.5em} 
\begin{theorem}
\label{theorem:existence-of-global-minimizers-empirical}
Let Assumption~\ref{assumption:X-continuous} and~\ref{assumption:non-degeneracy} hold. Then, for every $\lambda > 0$,
\begin{equation*}
	\lim_{n \to \infty} \P \left(J_n(\cdot, \lambda)~\text{admits a minimizer}~\Sigma_{n, *} \in \Sym_+^d\right) = 1.
\end{equation*}
\end{theorem}

\begin{proof} 
Since $X$ is continuous, Theorem~\ref{theorem:existence-of-global-minimizers} implies that $\lim_{\Sigma \to \infty} J(\Sigma, \lambda) = \E[Y^2]$.  
Since $\E[Y|X] \neq 0$ with positive probability, Lemma~\ref{lemma:upper-bound-U} implies 
$J(I, \lambda) = \mathcal{J}(I,\lambda) < \E[Y^2]$ for the identity matrix $I$. Consequently, there must exist a finite 
$M < \infty$ such that 
\begin{equation*}
	\inf_{\Sigma: \norm{\Sigma} \ge M} J(\Sigma, \lambda) > J(I, \lambda).
\end{equation*}
By Theorem~\ref{theorem:uniform-convergence}, almost surely as $n \to \infty$
\begin{equation*}
	\sup_{\Sigma \in \Sym_+^d} |J_n(\Sigma, \lambda) - J(\Sigma, \lambda)| \longrightarrow 0.
\end{equation*} 
This implies 
\begin{equation*}
	\lim_{n \to \infty} \P\left(\inf_{\Sigma: \norm{\Sigma} \ge M} J_n(\Sigma, \lambda) > J_n(I, \lambda)\right) = 1.
\end{equation*} 
On this event $\inf_{\Sigma: \norm{\Sigma} \ge M} J_n(\Sigma, \lambda) > J_n(I, \lambda)$, any minimizer of $J_n$ 
over the compact set $\{\Sigma: \norm{\Sigma} \le M\}$ is necessarily a global minimizer over the entire domain $\Sym_+^d$. 
Since $J_n$ is continuous, this guarantees that such a minimizer exists, completing the proof. 
\end{proof}


\vspace{.5em}

\subsection{Dimension Reduction} 

\indent\indent
This section presents the main theoretical findings of this paper: that the \emph{finite-sample} 
objective $J_n$ admits \emph{exact} low-rank global minimizers, achieving dimension reduction 
despite the absence of any explicit rank or nuclear-norm regularization on $\Sigma$. 

In particular, Theorem~\ref{thm:finite-sample-recovery} shows that for small $\lambda \le \lambda_0$, (i) every global minimizer 
$\Sigma_{n,*}$ of $J_n(\cdot, \lambda)$ exactly recover the correct subspace dimension $\dim(S)$ 
with high probability, and (ii) the corresponding low-rank subspaces spanned by $\Sigma_{n,*}$ 
converge to the true subspace $S$ as $n \to \infty$. Since $S$ captures the essential degree of freedom
of $X$ in predicting $Y$, in this precise sense we say that dimension reduction occurs at the finite sample levels.

\vspace{.5em} 
\begin{theorem}
\label{thm:finite-sample-recovery}
Let Assumptions~\ref{assumption:X-continuous}--~\ref{assumption:symmetry-of-RKHS} 
hold.  Then \vspace{.5em} 
\begin{itemize}
\item[$\mathsf{(a)}$] 
For every $\lambda > 0$, 
\begin{equation*}
	\lim_{n \to \infty} \P( \rank(\Sigma_{n, *}) \le \dim(S),~~~\forall \Sigma_{n, *}~\text{minimizing}~J_n(\cdot, \lambda)) = 1.
\end{equation*}
\item[$\mathsf{(b)}$] 
For a given $0 < \lambda \le \lambda_0$, where $\lambda_0$ is specified in Theorem~\ref{thm:population-recovery}$\mathsf{(b)}$, 
\begin{equation*}
	\lim_{n \to \infty} \P( \rank(\Sigma_{n, *}) = \dim(S),~~~\forall \Sigma_{n, *}~\text{minimizing}~J_n(\cdot, \lambda)) = 1.
\end{equation*} 
\item[$\mathsf{(c)}$] 
For a given $0 < \lambda \le \lambda_0$, where $\lambda_0$ is specified in Theorem~\ref{thm:population-recovery}$\mathsf{(b)}$, as $n \to \infty$,
\begin{equation*}
	\sup\left\{\norm{\Pi_{\Range(\Sigma_{n, *})}- \Pi_S} \mid \Sigma_{n, *}~\text{minimizing}~J_n(\cdot, \lambda) \right\} = o_P(1).
\end{equation*}  
\end{itemize}
\end{theorem} 
The $U$-formulation, Theorem~\ref{theorem:finite-sample-result-intro}, follows from Theorem~\ref{thm:finite-sample-recovery} via the reparameterization $\Sigma = U^\top U$.
In below, we prove Theorem~\ref{thm:finite-sample-recovery}. Our proof shows that the
sharpness property (Theorem~\ref{theorem:sharpness}) and uniform convergence (Theorem~\ref{theorem:uniform-convergence}) 
are the primary mechanisms that drive this result. 


%


\subsection{Proof of Theorem~\ref{thm:finite-sample-recovery}}
\subsubsection{Preparatory Lemmas}

\indent\indent
We start with a first-order stationarity condition for the empirical minimizer. 

\vspace{.5em} 
\begin{lemma} 
\label{lemma:empirical-minimizer-stationarity}
Let $\Sigma_{n,*}$ be a minimizer of $J_n(\cdot,\lambda)$. Then 
\begin{equation*}
	\mathrm{D}J(\Sigma_{n, *}, \lambda) [W] \ge 0,~~~\forall W \in \T_{\Sym_+^d} (\Sigma_{n, *}).
\end{equation*}
In particular, we have $\mathrm{D}J(\Sigma_{n, *}, \lambda) [vv^\top] = 0$ for all $v \in \Range(\Sigma_{n, *})$.
\end{lemma} 
\begin{proof}
Let $W \in \T_{\Sym_+^d} (\Sigma_{n, *})$. By definition, there exists $\eps>0$ such that $\Sigma_{n,*} + tW \in \Sym_+^d$
for all $t \in [0,\eps]$. By the optimality of $\Sigma_{n,*}$, it follows that $J_n(\Sigma_{n,*},\lambda) \le J_n(\Sigma_{n,*} + tW,\lambda)$. Hence, 
\begin{equation*}
	\mathrm{D} J_n(\Sigma_{n, *}, \lambda)[W] = 
		\lim_{t \downarrow 0^+} \frac{1}{t} \left(J_n(\Sigma_{n, *} + tW, \lambda) - J_n(\Sigma_{n, *}, \lambda)\right) \ge 0.
\end{equation*}
Let $v \in \Range(\Sigma_{n, *})$. Then both $vv^\top$ and $-vv^\top$ 
belong to $T_{\Sym_+^d} (\Sigma_{n, *})$. The above inequality implies  
\begin{equation*}
	\mathrm{D} J_n(\Sigma_{n, *}, \lambda)[vv^\top] \ge 0, \qquad \mathrm{D} J_n(\Sigma_{n, *}, \lambda)[-vv^\top] \ge 0.
\end{equation*}
Since $\mathrm{D} J_n(\Sigma_{n, *}, \lambda)[\cdot]$ is linear, it follows that $\mathrm{D} J_n(\Sigma_{n, *}, \lambda)[vv^\top] = 0$ 
for all $v \in \Range(\Sigma_{n, *})$.
\end{proof} 

The next result shows that empirical minimizers remain bounded and converge to population minimizers. The 
argument relies primarily on uniform convergence of the objective in Theorem~\ref{theorem:uniform-convergence}.

\vspace{.5em} 
\begin{lemma}
\label{lemma:minimizers-bounded-convergence}
Let $\{\Sigma_{n,*}\}_{n\ge 1}$ be any sequence such that each $\Sigma_{n,*}$ is a minimizer of $J_n(\cdot,\lambda)$. 
Then the sequence $\{\Sigma_{n,*}\}_{n\ge 1}$ is bounded. 
Moreover, every cluster point $\Sigma$ of $\Sigma_{n, *}$ is a minimizer of $J(\cdot, \lambda)$. 
\end{lemma} 

\vspace{.5em} 
\begin{proof} 
By Theorem~\ref{theorem:uniform-convergence}, we have the uniform convergence
\begin{equation*}
	\sup_{\Sigma \in \Sym_+^d} \left|J_n(\Sigma, \lambda) - J(\Sigma, \lambda)\right| \to 0.
\end{equation*} 
Since $\Sigma_{n,*}$ minimizes $J_{n}(\cdot,\lambda)$, using uniform convergence yields
\begin{equation}
\label{eqn:minimum-convergence}
	 \lim_{n \to \infty} J_{n}(\Sigma_{n,*}, \lambda) = 
	 	\lim_{n \to \infty} \min_{\Sigma \in \Sym_+^d} J_n(\Sigma, \lambda) =  \min_{\Sigma \in \Sym_+^d} J(\Sigma, \lambda).
\end{equation}

Suppose for the sake of contradiction that there exists a subsequence $\{n_k\}$ such that 
$\|\Sigma_{n_k,*}\| \to \infty$. By Theorem~\ref{theorem:existence-of-global-minimizers}, 
$\lim_{\Sigma \to \infty} J(\Sigma, \lambda) = \E[Y^2]$, so uniform convergence would imply 
\begin{equation*}
	\min_{\Sigma \in \Sym_+^d} J(\Sigma, \lambda)= \lim_{k \to \infty} J_{n_k}(\Sigma_{n_k,*}, \lambda) = 
		\lim_{k \to \infty} J(\Sigma_{n_k,*}, \lambda) = \E[Y^2].
\end{equation*}
However, because $\E[Y|X] \neq 0$, Lemma~\ref{lemma:lower-bound-U} implies $\min_{\Sigma \in \Sym_+^d} J(\Sigma, \lambda) \le J(I, \lambda) = \mathcal{J}(I,\lambda) < \E[Y^2]$, a contradiction! 
Therefore, $\{\Sigma_{n,*}\}_{n\ge 1}$ must be bounded.

Finally, let $\Sigma$ be any cluster point of $\{\Sigma_{n,*}\}_{n\ge 1}$, and let $\Sigma_{n_k,*} \to \Sigma$. 
Uniform convergence gives
\begin{equation*}
	J(\Sigma,\lambda) = \lim_{k \to \infty} J(\Sigma_{n_k, *}, \lambda) 
	= \lim_{k\to\infty} J_{n_k}(\Sigma_{n_k,*},\lambda)
	=  \min_{\Sigma \in \Sym_+^d} J(\Sigma, \lambda)
\end{equation*}
where the last equality follows from~\eqref{eqn:minimum-convergence}. 
This shows that $\Sigma$ minimizes $J(\cdot,\lambda)$.
\end{proof} 

Finally, we mention a standard perturbation result from matrix analysis, which ensures stability of subspaces under rank-preserving perturbations.

\vspace{.5em} 
\begin{lemma}
\label{lemma:standard-matrix-perturbation}
Let $\Sigma_n$ denote a sequence such that $\Sigma_n \to \Sigma$ as $n \to \infty$. 
Suppose $\rank(\Sigma_n) = \rank(\Sigma)$ for all $n$.  
Then, as $n \to \infty$, 
\begin{equation*}
	\matrixnorm{\proj_{\Range(\Sigma_n)} - \proj_{\Range(\Sigma)}} \to 0.
\end{equation*}
\end{lemma}

\vspace{.5em} 
\begin{proof}
This follows from Davis-Kahan Theorem~\cite[Section V]{StewartSu90}. 
\end{proof}

\vspace{.5em} 

\subsubsection{Main Argument}
\indent\indent
We start by proving part $\mathsf{(a)}$. We will show that almost surely, for any sequence 
$\{\Sigma_{n,*}\}_{n\ge 1}$ in which each $\Sigma_{n,*}$ is a minimizer of $J_n(\cdot,\lambda)$, we have
\begin{equation}
\label{eqn:rank-limit}
	\limsup_{n \to \infty} \rank(\Sigma_{n,*}) \le \dim(S).
\end{equation}
This would imply with probability one, there exists a finite $N< \infty$ such that $\rank(\Sigma_{n, *}) \le \dim(S)$ for all $n \ge N$, 
and for every minimizer $\Sigma_{n,*}$ of $J_n(\cdot,\lambda)$. Consequently, 
\begin{equation*}
\lim_{n \to \infty} 
\P\!\left(
\rank(\Sigma_{n, *}) \le \dim(S)
\ \text{for all minimizers } \Sigma_* \text{ of } J_n(\cdot,\lambda)
\right) = 1,
\end{equation*} 
which establishes the desired claim.

We now prove the rank limit in~\eqref{eqn:rank-limit}. It suffices to show that for any arbitrary subsequence $\{n_k\}$, 
and any arbitrary choice of $\Sigma_{n_k,*}$ such that each $\Sigma_{n_k,*}$ is a minimizer of $J_{n_k}(\cdot,\lambda)$, we have
\begin{equation}
	\limsup_{k \to \infty} \rank(\Sigma_{n_k,*}) \le \dim(S).
\end{equation}
Suppose for the sake of contradiction that this fails. Then there will be a subsequence $\{n_{k_j}\}$ where
$\rank(\Sigma_{n_{k_j},*}) > \dim(S)$. By Lemma~\ref{lemma:minimizers-bounded-convergence}, this 
subsequence $\{\Sigma_{n_{k_j},*}\}$ is bounded and therefore admits a cluster point $\Sigma_*$, which 
by the same lemma, must be a minimizer of $J(\cdot,\lambda)$. 

We now invoke the sharpness property. By Theorem~\ref{theorem:sharpness}, 
there exists $\rho > 0$ such that
\begin{equation*}
	\mathrm{D} J(\Sigma_*, \lambda)[W] \ge \rho \norm{W}, \qquad \forall W \in T_{\Sym_+^d}(\Sigma_*)~\text{with}~ \Range(W) \subset S^\perp.
\end{equation*}
Moreover, for any $v \in S^\perp$, the matrix $vv^\top$ belongs to
$T_{\Sym_+^d}(\Sigma_*)$ with $\Range(vv^\top) \subset S^\perp$. Hence,
\begin{equation*}
	\mathrm{D} J(\Sigma_*, \lambda)[vv^\top] \ge \rho \norm{v}^2,~~~\forall v\in S^\perp.
\end{equation*} 
Since $\Sigma_{n_{k_j}, *} \to \Sigma_*$, and $\mathrm{D} J(\cdot, \lambda)$ is continuous, we get for all large enough index $j$, 
\begin{equation}
\label{eqn:sharpness-empirical-one}
	\mathrm{D} J(\Sigma_{n_{k_j}, *}, \lambda)[vv^\top] \ge \frac{\rho}{2} \norm{v}^2,~~~\forall v\in S^\perp.
\end{equation}
On the other hand, by Lemma~\ref{lemma:empirical-minimizer-stationarity}, we have 
\begin{equation}
\label{eqn:sharpness-empirical-two}
	\mathrm{D} J(\Sigma_{n_{k_j}, *}, \lambda)[vv^\top] = 0,~~~\forall v\in \Range(\Sigma_{n_{k_j}, *}).
\end{equation} 
Since $\rank(\Sigma_{n_{k_j},*}) > \dim(S)$, a dimension-count argument implies that
\[
\Range(\Sigma_{n_{k_j},*}) \cap S^\perp \neq \{0\}.
\]
Thus there exists a nonzero vector $v \in \Range(\Sigma_{n_{k_j},*}) \cap S^\perp$. For such a $v$, the last two displays
\eqref{eqn:sharpness-empirical-one} and \eqref{eqn:sharpness-empirical-two} contradict each other. This contradiction completes the proof.

We next prove part $\mathsf{(b)}$. By Theorem~\ref{thm:population-recovery}, 
there exists $\lambda_0 < \infty$ such that for every $0 < \lambda \le \lambda_0$, 
\begin{equation}
\label{eqn:range-Sigma_*-S}
	\Range(\Sigma_*) = S~~~\forall \Sigma_*~\text{minimizing}~J(\cdot, \lambda).
\end{equation}
We prove that for the same $\lambda_0$, almost surely, for any sequence 
$\{\Sigma_{n,*}\}_{n\ge 1}$ in which each $\Sigma_{n,*}$ is a minimizer of $J_n(\cdot,\lambda)$ where $0 < \lambda \le \lambda_0$, we have
\begin{equation}
\label{eqn:rank-inf-limit}
	\liminf_{n \to \infty} \rank(\Sigma_{n,*}) \ge \dim(S).
\end{equation}
To see this~\eqref{eqn:rank-inf-limit}, recall that any cluster point $\Sigma_*$ of $\Sigma_{n,*}$ must 
be a minimizer of $J(\cdot,\lambda)$ by Lemma~\ref{lemma:minimizers-bounded-convergence}. 
By~\eqref{eqn:range-Sigma_*-S}, this minimizer satisfies $\rank(\Sigma_*) = \dim(S)$. Given that the rank 
function $\Sigma \mapsto \rank(\Sigma)$ is lower-semicontinuous, the limit~\eqref{eqn:rank-inf-limit} follows. 

Combining~\eqref{eqn:rank-limit} and~\eqref{eqn:rank-inf-limit}, these two results imply with probability one, 
\begin{equation*}
	\lim_{n \to \infty} \rank(\Sigma_{n,*}) = \dim(S).
\end{equation*}
Therefore, there exists a finite $N< \infty$ such that $\rank(\Sigma_{n, *}) = \dim(S)$ for all $n \ge N$, 
and for every minimizer $\Sigma_{n,*}$ of $J_n(\cdot,\lambda)$. Consequently, 
\begin{equation*}
\lim_{n \to \infty} 
\P\!\left(
\rank(\Sigma_{n, *}) = \dim(S)
\ \text{for all minimizers } \Sigma_* \text{ of } J_n(\cdot,\lambda)
\right) = 1,
\end{equation*} 
which establishes the desired conclusion.

Finally, we prove part $\mathsf{(c)}$. Every cluster point $\Sigma$ of the sequence ${\Sigma_{n,*}}$ satisfies
$\Range(\Sigma) = S$. Moreover, by the argument in part $\mathsf{(b)}$, we have already shown that 
$\rank(\Sigma_{n, *}) = \dim(S)$ for all sufficiently large $n$, almost surely. Therefore, the convergence of the 
associated range subspaces
\begin{equation*}
	\sup\left\{\norm{\Pi_{\Range(\Sigma_{n, *})}- \Pi_S} \mid \Sigma_{n, *}~\text{minimizing}~J_n(\cdot, \lambda) \right\} = o_P(1)
\end{equation*}  
follows directly from Lemma~\ref{lemma:standard-matrix-perturbation}. This completes the proof.

\section{Invariance and Relaxation of Assumptions} 
\label{sec:invariance-and-relaxation}
\indent\indent 
At this point, the main results, namely
Theorems~\ref{theorem:population-result-intro}
and~\ref{theorem:finite-sample-result-intro} have been established.
We now highlight an invariance property of the model that plays an
important conceptual role: it allows the assumptions imposed in these theorems 
to be relaxed without affecting the conclusions.

\subsection{Invariance under Linear Reparameterization} 
\indent\indent 
Let $(X,Y)\sim \P$, and define $(X',Y')\sim\P'$ by
\begin{equation*}
Y' = Y,
\qquad
X' = AX,
\end{equation*}
where $A \in \mathbb{R}^{d \times d}$ is invertible. Let $S \subset \mathbb{R}^d$
denote the central mean subspace of $Y$ given $X$ under $\P$, and let $S'$
denote the corresponding one under $\P'$. Consider the variational problems $\text{(P)}$ and $\text{(P')}$
\begin{equation*}
\begin{split}
\text{(P)}\quad  &\min_{U} \min_{f} \; \E\!\left[(Y - f(UX))^2\right] + \lambda \|f\|_H^2, \\
\text{(P')}\quad  &\min_{U} \min_{f} \; \E'\!\left[(Y' - f(UX'))^2\right] + \lambda \|f\|_H^2,
\end{split}
\end{equation*}
where $\E[\cdot]$ and $\E'[\cdot]$ denote expectation under $\P$ and $\P'$,
respectively. The following lemma summarizes the invariance properties of both the
predictive subspace and the associated variational problems.

\vspace{.5em} 

\begin{lemma}
\label{lem:affine-invariance}
With notation as above, the central mean subspaces satisfy
\begin{equation*}
S' = A^{-\top} S,~~\text{where}~~A^{-\top} S := \{ A^{-\top} v : v \in S \} \subset \mathbb{R}^d.
\end{equation*}
Moreover, $(U_*, f_*)$ is a minimizer of \textnormal{(P)}, 
if and only if $(U_* A^{-1}, f_*)$ is a minimizer of \textnormal{(P')}.
\end{lemma}

\subsection{Generalization via Linear Reparameterization} 
\indent\indent
We now introduce a relaxed version of Assumption~\ref{assumption:independence},
allowing an invertible linear reparameterization of the covariates. 
Assumption~\ref{assumption:independence} corresponds to the special case $A = I$.

\renewcommand{\theassumption}{1.3'}

\vspace{.5em} 
\begin{assumption}
\label{assumption:independence-new}
There exists a invertible linear map $A$ such that, for
$X' := AX$, the random vectors $\Pi_S X'$ and $\Pi_{S^\perp} X'$ are
probabilistically independent. Moreover, $\Pi_{S^\perp} X'$ is
non-degenerate, in the sense that $\Cov(\Pi_{S^\perp} X')$ is not the zero matrix.
\end{assumption} 

\vspace{.5em}
This assumption permits, for example,
$X \sim \normal(\E[X], \Cov(X))$ with an arbitrary non-degenerate covariance
matrix.

\vspace{.5em}

\begin{theorem}
\label{theorem:dependence}
Let $(X, Y) \sim \P$ satisfy Assumptions~\ref{assumption:X-continuous},
\ref{assumption:non-degeneracy},~\ref{assumption:independence-new}. Suppose 
further that Assumption
\ref{assumption:symmetry-of-RKHS} holds. Then the conclusions of
Theorems~\ref{theorem:population-result-intro} and~\ref{theorem:finite-sample-result-intro}  hold.
\end{theorem}

\begin{proof} 
By Assumption~\ref{assumption:independence-new}, there exists an invertible
linear map $A$ such that, for $X' := AX$ and $Y' := Y$, the random vectors
$\Pi_S X'$ and $\Pi_{S^\perp} X'$ are probabilistically independent, and
$\Pi_{S^\perp} X'$ is non-degenerate. Let $\P'$ denote the law of $(X',Y')$.
By Lemma~\ref{lem:affine-invariance}, the corresponding central mean subspace
under $\P'$ is $S' = A^{-\top} S$, where $A^{-\top} S := \{ A^{-\top} v : v \in S \} \subset \mathbb{R}^d$.

Since $A$ is invertible, $(X',Y') \sim \P'$ satisfies
Assumptions~\ref{assumption:X-continuous} and~\ref{assumption:non-degeneracy}.
Moreover, by construction,
Assumption~\ref{assumption:independence} holds for $(X',Y')$ with respect to
$S'$. Therefore, the conclusions of
Theorems~\ref{theorem:population-result-intro}
and~\ref{theorem:finite-sample-result-intro} apply to the distribution $\P'$. In particular, 
if $(U'_*, f'_*)$ denotes a minimizer of the empirical or population objective under $\P'$, 
then $U'_*$ recovers $S'$ in the sense stated in those theorems.
\end{proof} 

\vspace{.5em} 

\begin{remark}
\emph{
Conceptually, the linear map $A$ should be viewed as a change of
coordinates. The invariance properties above reflect the intrinsic nature of the variational
objective and the associated predictive subspaces. 
In this paper, invariance is
verified through explicit reparameterization. A fully differential–geometric 
formulation makes this invariance transparent; see~\cite{LiRu25}.
}
\end{remark} 

\vspace{-.25em} 

\newcommand{\subjectto}{{\rm subject~to}}
\section{Numerical Experiments} 
\label{sec:numerical-experiments}
Code is available at
\url{
https://github.com/tinachentc/kernel-learning-in-ridge-regression
}.

\subsection{Synthetic Datasets} 
\indent\indent
This section documents additional simulation evidence, showing that compositional kernel 
model favors low-rank solution $U_{n, *}$ in finite samples, despite the absence of explicit 
low-rank penalties. 


\vspace{-.25em} 
\subsubsection{Setting and Methodology} 
\label{sec:setting-and-methodology}
Our experiments are of the following form. 

\begin{enumerate}
\item We generate $n$  i.i.d. samples $\{(X_i, Y_i)\}_{i=1}^n \sim \P$ according to 
	\begin{equation*}
		Y = F(X) + \epsilon
	\end{equation*}
	where the noise term $\epsilon \sim \normal(0, \sigma^2)$ is independent from $X \sim \P_X$. 
\item We solve the finite-sample compositional kernel learning problem in its
$\Sigma$-formulation:
	\begin{equation}
	\label{eqn:general-kernel-learning-objective}
		\minimize_{\Sigma}~~~J_n(\Sigma, \lambda)
			~~\subjectto~~\Sigma \succeq 0.
	\end{equation} 
	The expression for $J_n$ follows directly from the formula derived in Section~\ref{sec:illustration-under-finite-atomic-measures}: 
\begin{equation*}	
	J_n(\Sigma, \lambda) =
	\lambda\,
	\y^\top
	\bigl(\K_\Sigma + n\lambda I_n\bigr)^{-1}
	\y,
\end{equation*}
where $\K_\Sigma  \in \R^{n\times n}$ has entries $(\K_\Sigma)_{ij} = \psi\!\bigl(\norm{X_i-X_j}_\Sigma^2\bigr)$,
and $\y = (Y_1,\ldots,Y_n)^\top \in \R^n$.	
\end{enumerate}

\vspace{.5em} 
The main purpose of the experiments is to vary the signal form $F$ and $\P_X$ 
to investigate when $\rank(\Sigma_{n, *}) \le \dim(S)$ and 
$\rank(\Sigma_{n, *}) = \dim(S)$ hold with 
high probability.

\vspace{.3cm}
\noindent\noindent
\textbf{Function $F$.}~~
The function $F$ takes one of the following forms: (a) $F(x) = x_1 + x_2 + x_3$ (b) $F(x) = x_1 x_2$, 
(c) $F(x) = 0.1(x_1 + x_2 + x_3)^3 + \tanh(x_1 + x_3 + x_5)$, (d) $F(x)= 2(x_1 + x_2) + (x_2 + x_3)^2 + (x_4 - 0.5)^3$, 
(e) $F(x) = 0$.

\vspace{.3cm}
\noindent\noindent
\textbf{RKHS $H$.}~~
We take $H$ to be the Gaussian RKHS with kernel $(x, x') \mapsto \exp(-\norm{x-x'}^2)$. 

\vspace{.3cm}
\noindent\noindent
\textbf{Algorithm}.~~
We minimize $J_n$ using projected gradient descent with projection onto the semidefinite cone $\Sym_+^d$, employing an Armijo line search for stepsize selection. The algorithm is initialized at $(1/d) I$ and terminated when the ratio between the difference of consecutive iterates, measured by the Frobenius norm, and the stepsize is below a threshold $\Delta > 0$.

\vspace{.3cm}
\noindent\noindent
\textbf{Size $n, d, \sigma, \Delta$.}~~ For all our simulations, $n = 300$, $d = 50$, $\sigma = 0.1$,
$\Delta = 10^{-3}$. 

\subsubsection{Results} 
We first list our claims and then provide the source of evidence that supports our claims.

\vspace{.5em}
\begin{enumerate}
\item[(1)] The phenomenon $\rank(\Sigma_{n, *}) \le \dim(S)$ occurs with high probability 
under Assumptions~\ref{assumption:X-continuous}--~\ref{assumption:symmetry-of-RKHS}.
This is consistent with the guarantee in Theorem~\ref{theorem:finite-sample-result-intro}.
\item[(2)] The phenomenon $\rank(\Sigma_{n, *}) \le \dim(S)$ can occur in the solution with high probability even if 
$X$ is discrete. This shows that Assumption~\ref{assumption:X-continuous} is not necessary for the phenomenon to hold. 
\item[(3)] The phenomenon disappears when $\E[Y|X] \equiv 0$. Hence, Assumption~\ref{assumption:non-degeneracy} is necessary. 
\item[(4)] The phenomenon $\rank(\Sigma_{n,*}) \le \dim(S)$ can occur in the solution with high probability even if 
the covariates of $X$ are dependent. As shown in Theorem~\ref{theorem:dependence}, Assumption~\ref{assumption:independence} is not required for this phenomenon to hold. Our simulations illustrate this behavior in dependent settings.
\end{enumerate}
\vspace{.5em}

Below we make each claim concrete with evidence. 

\vspace{.5em}
\noindent\noindent
\textbf{Claim (1).} We generate covariates according to $X \sim \normal(0, I)$ For all functions $F$ in Section~\ref{sec:setting-and-methodology}(a)--(d), we observe that $\rank(\Sigma_{n,*}) \le \dim(S)$ with high probability; see Figure~\ref{plot:claim1} in Appendix~\ref{sec:additional-plots}. 

\vspace{.3cm}
\noindent\noindent
\textbf{Claim (2).} We consider a discrete setting where each coordinate of $X$ is independently distributed as $\mathrm{Bernoulli}(0.5)$. 
For all functions $F$ in Section~\ref{sec:setting-and-methodology}(a)--(d), we observe that $\rank(\Sigma_{n,*}) \le \dim(S)$ with high probability; see Figure~\ref{plot:claim3} in Appendix~\ref{sec:additional-plots}.

\vspace{.3cm}
\noindent\noindent
\textbf{Claim (3).} 
When $F \equiv 0$ and $X \sim \normal(0,I)$, we observe that $\Sigma_{n,*}$ is full rank across a wide range of $\lambda$ values in our experiments; see Figure~\ref{plot:claim4} in Appendix~\ref{sec:additional-plots}.

\vspace{.3cm}
\noindent\noindent
\textbf{Claim (4).} We generate covariates according to $X \sim \normal(0,C)$ with $C_{ij}=0.5^{|i-j|}$, thus inducing dependence among coordinates.
For all choices of $F$ in Section~\ref{sec:setting-and-methodology}(a)--(d), the inequality $\rank(\Sigma_{n, *}) \le \dim(S)$ occurs with high probability; see Figure~\ref{plot:claim2} in Appendix~\ref{sec:additional-plots}.

\subsection{Real Data Experiments}
\indent\indent
We evaluate several models on the SVHN (Street View House Numbers) dataset~\cite{NetzerWaCoBiWuNg11}, a real-world image dataset for digit classification consisting of approximately 100,000 labeled images of house numbers captured in natural scenes. The task is to identify the digit at the center of each image, which typically appears amid neighboring digits and complex backgrounds.

\begin{figure}[!htb]
\centering
\begin{subfigure}[b]{.9\columnwidth}
\centering
\includegraphics[width=\textwidth]{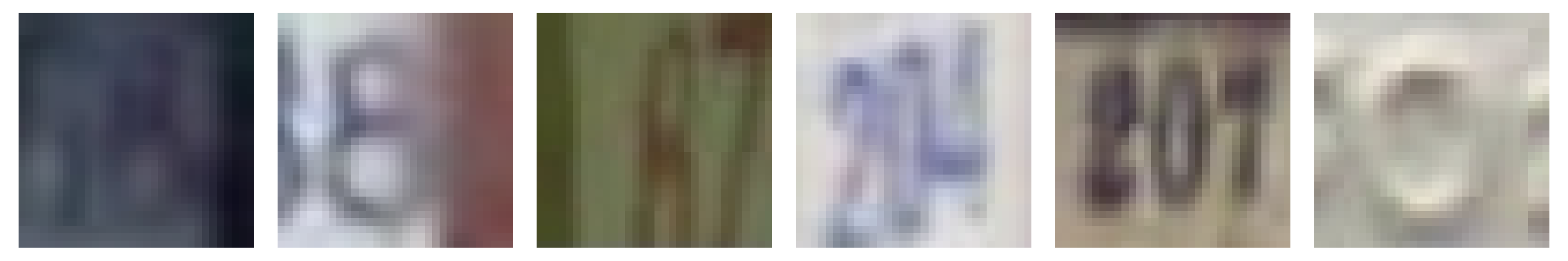}
\end{subfigure}
\caption{Sample images of 8 and 0 from SVHN data.}
\label{fig:sample}
\end{figure}

We focus on a binary classification task distinguishing digits $0$ and $8$. Example images are shown in Figure~\ref{fig:sample}. The training set contains $n = 9993$ samples with input dimension $d = 3 \times 32^2$, corresponding to RGB images of size $32 \times 32$. The test set contains $n_{\mathrm{test}} = 3404$ samples.


\vspace{-1.5em}

This task naturally favors compositional architectures that separate representation learning from prediction. In such models, the representation map acts as a feature extractor, projecting images onto a lower-dimensional subspace that emphasizes the central digit while suppressing peripheral content. The predictor is then applied to this reduced representation to perform classification.

\textbf{Goal}. The goal of these experiments is to compare how \emph{classical kernel ridge regression},
which uses a fixed representation, and \emph{compositional kernel models}, which learn a
representation $U$, differ in their learned structure and generalization behavior on a
real-world vision task. This comparison is designed to \emph{isolate the effect of learning
the representation} $U$ within an otherwise classical kernel framework. Two-layer ReLU neural
networks are included mainly as a reference point for intuition.


\subsubsection{Feature Learning via Compositions}
\indent\indent
To isolate the effect of feature learning through $U$, we compare the compositional kernel model with standard kernel ridge regression.
Recall the compositional kernel model
\begin{equation*}
\begin{split}
      ~~~\text{Compositional kernel model}:~~~~~~~ \min_U \min_f \E[(Y- f(U X))^2] + \lambda \norm{f}_{H}^2.
\end{split}
\end{equation*}
We compare this formulation with standard kernel ridge regression with bandwidth selection,
\begin{equation*}
    \text{KRR with bandwidth selection}:~~~~\min_{h \ge 0} \min_f  \E[(Y- f(X))^2] + \lambda \norm{f}_{H_h}^2.
  \end{equation*}
In our experiments, $H$ is the Gaussian RKHS induced by the kernel $k(x, x') = \exp\!\left(-\|x-x'\|^2\right)$, 
while $H_h$ is the the Gaussian RKHS induced by $k_h(x,x') = \exp\!\left(-h \|x-x'\|^2\right)$ where $h$ is a scalar bandwidth parameter.
The key distinction between these two approaches is that the compositional kernel model learns an additional data representation through the matrix $U$, whereas the other relies on an isotropic kernel and adapts only the scalar bandwidth parameter $h$.
  
\begin{figure}[!htb]
\centering
\begin{subfigure}[b]{\columnwidth}
\centering
\includegraphics[width=.8\textwidth]{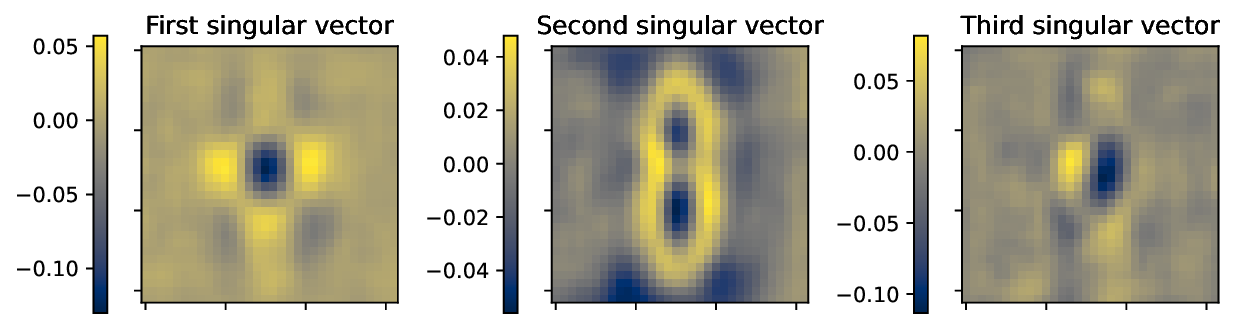}
\caption{compositional kernel model}
\label{fig:svhn_kernel}
\end{subfigure}
\vspace{1em}
\begin{subfigure}[b]{\columnwidth}
\centering
\includegraphics[width=.8\textwidth]{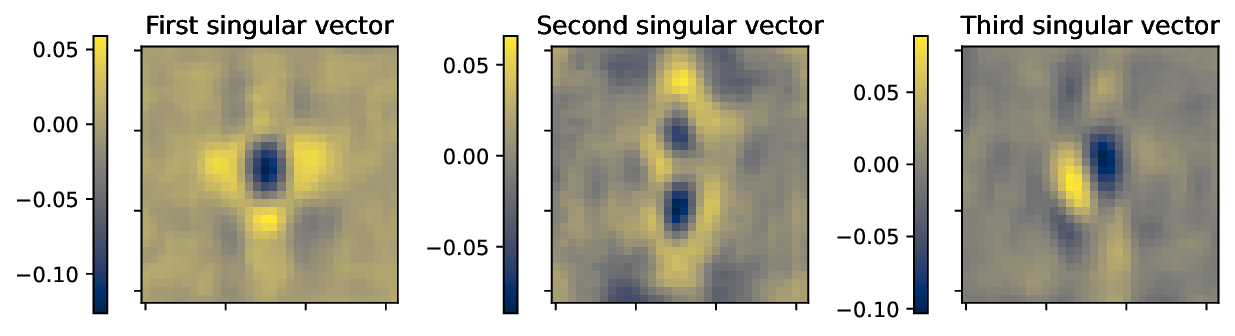}
\caption{two-layer ReLU neural network}
\label{fig:svhn_nn}
\end{subfigure}

\caption{The first row displays the top three left singular vectors of $U$ of a top-performing compositional kernel model, which 
achieves a prediction accuracy of $92.83\%$. The second row shows 
the top three singular vectors of the bottom layer weight matrix from a top-performing 
two-layer ReLU neural network on SVHN data,  with a prediction accuracy of $90.80\%$. 
Each vector technically has a dimension of $3 \times 1024$ (where $3$ stands for the red, green, and blue 
channels, and $1024$ stands for the size of the image); for illustration, we only present the subvector 
of dimension $1024$ that corresponds to the red channel. Both models seem to learn visually interpretable features. 
}
\label{fig:singularvectors}
\end{figure}
  

Empirically, compositional kernel model substantially outperforms standard kernel ridge regression. On the full dataset, test error drops from $12.0\%$ to $7.2\%$, a relative reduction of approximately $40\%$. This improvement is stable under subsampling: even with only $20\%$ of the training data, the error decreases from $20\%$ to $11\%$, reflecting a similar relative 45\% reduction. These improvements cannot be attributed to bandwidth tuning alone and reflect the benefit of learning representations $U$.

To build further intuitions, we visualize the top three singular vectors of the \emph{learned} matrix $U$.
This visualization reveals that the learned top singular vectors capture patterns that highlight structural differences between digits $0$ and $8$.

%

\begin{figure}[!htb]
    \centering
    \begin{subfigure}[b]{0.4\textwidth}
        \includegraphics[width=\textwidth]{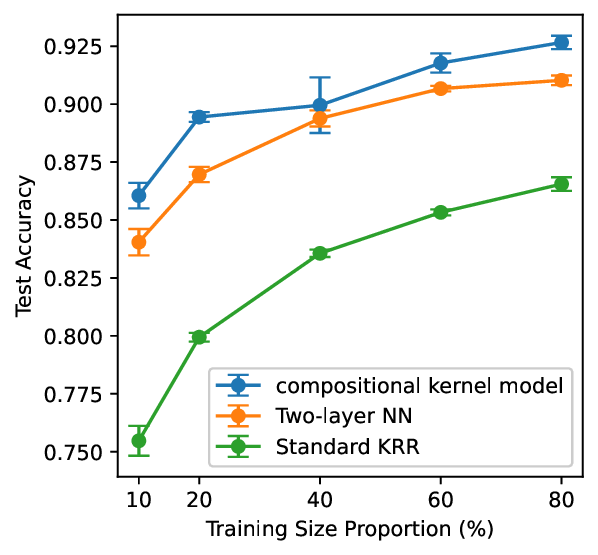}
    \end{subfigure}
    \caption{
Prediction accuracy for classifying digits $0$ and $8$. For each method—compositional kernel models, two-layer neural networks, and standard kernel ridge regression—we subsample the SVHN training set to proportions of $10\%$, $20\%$, $40\%$, $60\%$, and $80\%$, train the corresponding model on each subsample, and evaluate performance on a fixed test set. Each subsampling proportion is repeated 5 times, and the average accuracy is reported with standard error bars.}
    \label{fig:svhn_ssub}
\end{figure}

\subsubsection{The prevalence of low-rank structure in learned feature matrices}
\indent\indent
We observe a pronounced low-rank bias in the feature matrices $U$ learned by compositional kernel model, providing empirical evidence that the model favors parsimonious representations that is consistent with Theorem~\ref{theorem:finite-sample-result-intro}. In our experiments, we solve the minimization problem using the $\Sigma$-formulation via projected gradient descent, as described in Section~\ref{sec:setting-and-methodology}; this choice is justified by the equivalence between the $U$- and $\Sigma$-formulations established earlier.


\begin{figure}[!htb]
    \centering
    \begin{subfigure}[b]{0.4\textwidth}
        \includegraphics[width=\textwidth]{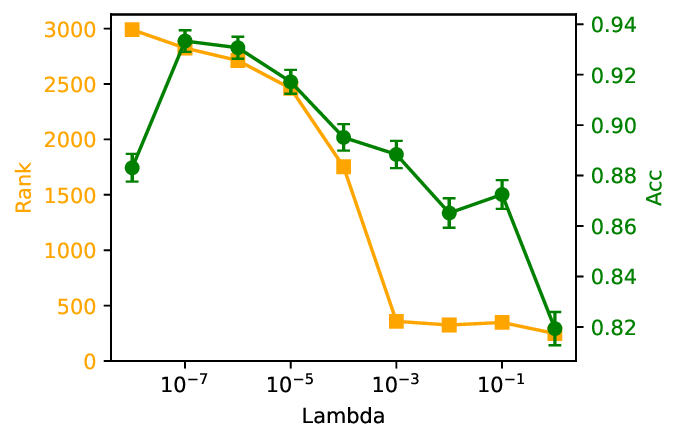}
    \end{subfigure}
    \caption{Prediction accuracy (green curve) and rank of bottom-layer weight matrix $U$ (orange curve) from the two-layer 
    kernel ridge regression model as functions of 
    $\lambda$, with $\lambda$ values taken from the set $\{10^{-8},10^{-7},10^{-6},10^{-5},10^{-4},10^{-3},10^{-2},10^{-1},1$\}.}
    \label{fig:svhn_select}
\end{figure}

Figure~\ref{fig:svhn_select} shows prediction accuracy and the rank of the learned matrix $U$ as functions of $\lambda$.
Across a wide range of regularization levels, the learned representation remains low-rank while maintaining strong predictive performance.
In particular, as $\lambda$ increases, the rank of $U$ decreases sharply---from roughly $3000$ to below $500$—while accuracy remains around $89\%$, exceeding that of standard kernel ridge regression.

Examining the singular spectrum of a top-performing solution with $\lambda = 10^{-7}$ reveals rapid decay in singular values, indicating that most predictive information is concentrated in a small number of singular vector directions. Indeed, retaining only the top five singular vectors of $U$ and retraining $f$ preserves over $90\%$ of the original accuracy, demonstrating that performance is driven by a \emph{low-dimensional feature subspace}.

For comparison, the bottom-layer weights of a two-layer ReLU network exhibit slower singular value decay and a higher effective rank. While both models appear to rely on low-dimensional subspaces for predictive performance, compositional kernel model achieves a slightly more compact representation, as reflected by its sharper spectral decay and lower effective rank.

\begin{figure}[!htb]
    \centering
    \begin{minipage}[t]{0.4\textwidth}
        \centering
        \includegraphics[width=\textwidth]{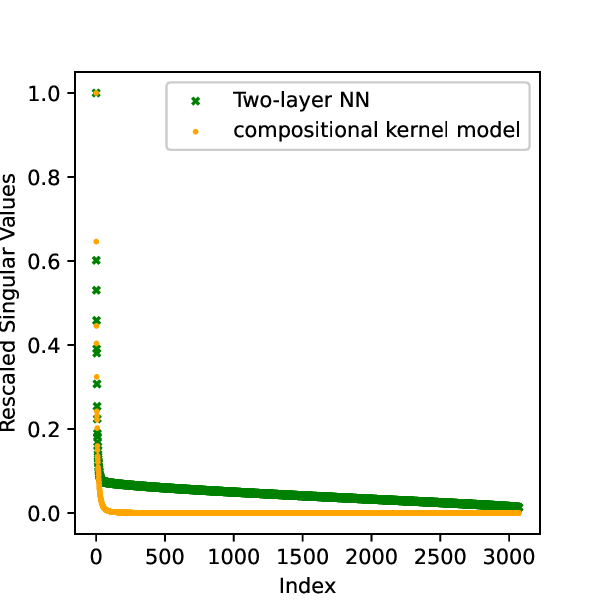}

    \end{minipage}%
    \begin{minipage}[t]{0.6\textwidth}
        \centering
        \vspace{-6cm}  
        \begin{tabular}{lc}
            \hline
            Method & Accuracy (\%) \\
            \hline
            Compositional KRR & 92.83\\
            Standard KRR &  87.96 \\
            Two-layer NN & 90.80\\
            Compositional KRR (Top 1 Singular Vector) &  61.13\\
            Two-layer NN (Top 1 Singular Vector) & 55.49\\            
            Compositional KRR (Top 5 Singular Vectors) &  90.63\\
            Two-layer NN (Top 5 Singular Vectors) & 89.31\\  
            Compositional KRR (Top 50 Singular Vectors) &  92.83\\
            Two-layer NN (Top 50 Singular Vectors) & 90.69\\  
            \hline
        \end{tabular}

    \begin{tabular}{lc}
    \hline
        Method & Effective Rank \\
    \hline
       	Compositional KRR & 7.87\\
	Standard KRR & 3072\\
        Two-layer NN & 132.01\\
    \hline
    \end{tabular}

    \end{minipage}
    \vspace{.1em} 
    \caption{The compositional kernel model displays similar prediction performance to two-layer neural networks with more compressed feature representations: (a) singular values of the bottom layer weight matrices $U$, where orange dots correspond to our model and green crosses correspond to neural networks (b) prediction accuracy for different methods (c) effective rank (sum of singular values divided by the maximum singular value) for different methods on SVHN dataset.}
    \label{fig:svhn_results}
\end{figure}


\section{Further Discussions}
\label{sec:further-discussions}
\indent\indent
This section provides additional perspective and is not required for the
understanding of the main results.

\subsection{A Variational Perspective of Regularization} 
\label{sec:variational-perspective-of-regularization}
\indent\indent
The results of this paper establish why low-dimensional predictive minimizers emerge in \emph{finite-sample} compositional kernel objectives through a variational calculus perspective. In this section, we place these findings in a broader context and discuss how they relate to classical regularization-based viewpoints in statistics and machine learning.

In statistics and machine learning, regularization is commonly interpreted as techniques to enforce structural simplicity in solutions, typically by adding penalty terms to the loss function \cite{HastieTiFr09}. Classical examples include $\ell_1$ penalties to promote sparsity and nuclear norm penalties to induce low-rank structure. The standard intuition is that such structures are inherently unstable under perturbations and must therefore be imposed explicitly~\cite{HastieTiFr09}.

This paper shows, however, that assessing sparsity or low-rank structure solely through the presence or absence of explicit penalty terms does not capture the full picture. Our perspective is grounded in sensitivity analysis and set identifiability techniques from the optimization literature~\cite{Wright93,Lewis02,HareLe07,Drusvyatskiys13,DrusvyatskiyLe14}, which address the following question: if  minimizers of $J$ lie in a structured set $\mathsf{S}$, under what conditions will minimizers of smoothly perturbed objectives $\tilde{J}$ also lie in $\mathsf{S}$?

A recurring answer is that such stability is governed by a \emph{sharpness} property of the objective $J$ around $\mathsf{S}$: roughly speaking, the objective must grow at least linearly when moving away from $\mathsf{S}$. In classical nonlinear programming, when $\mathsf{S}$ is an identifiable surface of a convex constraint set, sharpness appears as \emph{strict complementarity slackness condition} \cite{Wright93}. In smooth settings where $\mathsf{S}$ is a submanifold of a Euclidean space, sharpness is formalized through the notion of \emph{partial smoothness} \cite{Lewis02}. For more general sets, the appropriate notion is that of an \emph{identifiable set} \cite{Drusvyatskiys13,DrusvyatskiyLe14}.

Viewed through this lens, the finite-sample behavior of compositional kernel objectives can be interpreted as a concrete instance of the broader variational principle described above. Specifically, empirical minimizers are expected to exhibit low-dimensional structure if two conditions are met: (i) \emph{sharpness} of the population objective relative to an underlying structured set $\mathsf{S}$, and (ii) \emph{uniform convergence} of the empirical objective to its population counterpart in value and gradient. In our setting, where $\mathsf{S}$ denotes the set of low-rank matrices, the established sharpness properties for compositional kernel models offer a concrete explanation for the observed stability of low-rank structure under \emph{finite sampling}---a form of perturbation---despite the absence of explicit penalty.


\subsection{Parameterization and Algorithmic Identification}
\indent\indent
Careful readers may notice that, in the numerical experiments, we implement the algorithm using the equivalent $\Sigma$-parameterization rather than the $U$-parameterization. This choice is fundamental; it is dictated by algorithmic considerations, rather than by any difference in the low-rank structure of the underlying minimizers. We explain the rationale for this choice below.

A useful concept for interpreting our numerical experiments is that of \emph{active set identification} in optimization~\cite{Bertsekas03, FacchineiFiKa98, OberlinWr06}. Let $J(\theta)$ be an objective function and let $\mathsf{S}$ denote a structured subset of the parameter space (e.g., the set of low-rank matrices). A numerical algorithm $\{\theta^{(t)}\}_{t\in \N}$ is said to identify the set $\mathsf{S}$ if there exists a finite iteration index $T$ such that all subsequent iterates satisfy
\begin{equation*}
	\theta^{(t)} \in \mathsf{S},~~~~\forall t\ge T. 
\end{equation*}
In other words, the iterates enter $\mathsf{S}$ in finite iterations and remain in $\mathsf{S}$ thereafter. Identification is a stronger notion than algorithm asymptotic convergence to $\mathsf{S}$ (i.e., $\mathrm{dist}(\theta^{(t)},\mathsf{S}) \to 0$ as $t \to \infty$), which does not imply that the iterates ever enter $\mathsf{S}$ in finite iteration $t$. 

This notion directly underlies our implementation choices. In the experiments, our goal is to observe \emph{finite-time} identification of $\mathsf{S}$, since algorithms can only be run for finitely many iterations; in our setting, $\mathsf{S}$ corresponds to the set of low-rank matrices.

It is well known in variational analysis and optimization that active set identification typically requires more than the mere existence of structured minimizers~\cite{Bertsekas03, BurkeMo88, Burke90, AlKy91}. In particular, finite-time identification is closely tied to a sharpness property of the objective relative to $\mathsf{S}$ at the minimizer, together with an numerical algorithm whose updates are compatible with this geometry. 

For this reason, we optimize the objective using the $\Sigma$-parameterization, in which the objective $J(\cdot, \lambda)$ exhibits sharpness relative to the set of low-rank matrices, enabling projected gradient methods to identify low-rank structure in finitely many iterations~\cite{Wright93}. By contrast, this sharpness is not reflected in the $U$-parameterization, and standard gradient-based methods do not identify low-rank solutions in finite time, even when such solutions exist. 

Retrospectively, the mathematical equivalence between the $U$- and $\Sigma$-parameterizations is essential: without this analysis, the low-rank recovery phenomenon would be \emph{experimentally difficult to detect/verify} using standard numerical methods (e.g., gradient descent) applied to $U$-formulation.

\subsection{Analogy to Two-layer Neural Networks} 
\label{sec:analogy-to-two-layer-NN}
\indent\indent
The proposed compositional kernel formulation bears a structural analogy to two-layer neural networks. Recall 
the compositional kernel formulation
\begin{equation*}
	\min_U \min_f \E[(Y - f(UX))^2] + \lambda \norm{f}_H^2
\end{equation*}
where $U \in \R^{d \times d}$ is a linear transformation. 

By comparison, a (population) two-layer neural network can be written in the form
\begin{equation*}
	\min_a \min_W \E[(Y - \langle a, \sigma(WX) \rangle)^2] + \lambda \norm{a}^2
\end{equation*} 
Here, $W \in \mathbb{R}^{N \times d}$ denotes the weight matrix of the first (linear) layer, $\sigma(\cdot)$ is a fixed activation function applied elementwise, and $a \in \mathbb{R}^N$ contains the second-layer weights. 

In both formulations, the first layer is a linear map ($x \mapsto Ux$ or $x \mapsto Wx$) that defines a learnable representation of the input, while the second layer performs prediction on the transformed features. Conditional on the first-layer parameters $U$ or $W$, both objectives are quadratic in the second-layer variables $f$ or $a$. In this analogy, the rank of $W$ plays a role analogous to the rank of $U$.


This naturally raises the question of how the results presented here may relate to two-layer neural networks. A more detailed discussion of this connection is provided in the postscript of~our companion work~\cite{LiRub25}.  At present, most  results established in this paper—particularly those concerning population-level identifiability, sharpness, and finite-sample structure recovery—have no known analogues proven for two-layer neural networks at comparable levels of generality. 

One may therefore adopt either of the following viewpoints: \vspace{.5em}
\begin{itemize}
\item 
\emph{Toy model perspective}. The compositional kernel problem can be viewed as a simplified, analytically tractable surrogate for two-layer neural networks. From this perspective, the present results provide insight into the feature-learning mechanisms—such as sharpness and identifiability—that may also be present in neural network models, and motivate future theoretical and numerical investigations in that setting. 
\vspace{.5em} 
\item 
\emph{Model-specific structure perspective}. Alternatively, the compositional kernel problem may possess mathematical structure that is absent in neural networks, while neural networks benefit from highly developed numerical optimization schemes. From this viewpoint, neural network methods may offer practical algorithmic tools that could potentially be adapted to solve the compositional kernel problems efficiently.
\end{itemize}  \vspace{.5em}
We leave the assessment of these viewpoints to the reader.

\section{Open Questions} 
\label{sec:open-questions}
\indent\indent
This paper sits within the broader emerging program on understanding \emph{feature learning} in
the kernel ridge regression framework~\cite{FukumizuBaJo04, FukumizuBaJo09, Adit24A, FollainBa24, RadhakrishnanBeDr25, HuangLaDePoRo25, ZhuDaDrFa25, LiRu25, LiRub25, RuanLiJo25}.  This line of work treats kernel ridge regression as a \emph{canonical} baseline, and studies how the introduction of additional structure gives rise to feature learning---not because it captures all aspects of modern learning systems, but because its variational structure is sufficiently clean to allow underlying mechanisms to be \emph{isolated} and analyzed precisely. With this perspective in mind, important questions remain:

\vspace{.5em} 
\begin{itemize}
\item (Variational Geometry) 
Theorem~\ref{theorem:sharpness-general} establishes a linear growth (sharpness property) of the population objective $J$ at its global minimizer $\Sigma_*$ in directions normal to the low-rank matrix manifold $\mathcal M$, viewed as a submanifold of $\Sym^d$.  A  complementary question concerns the behavior of $J$ along directions tangent to $\mathcal M$, along which we conjecture a subquadratic growth.

\vspace{.5em}
\item (Stationary Points) 
This paper focuses on the structure of global minimizers. Under a Gaussianity assumption on $\Pi_{S^\perp} X$, it is shown in~\cite{LiRub25} that any nontrivial stationary point $\Sigma$---that is, satisfying $J(\Sigma,\lambda) < \mathbb{E}[Y^2]$---obeys $\Range(\Sigma) \subset S$ at the population level. Whether stationary points recover the signal subspace $S$ (namely, $\Range(\Sigma) = S$), and how their structure (e.g., low-rankness) behaves in finite samples, remain open questions.


 \vspace{.5em} 
\item (Nonsmooth RKHS)
This paper focuses on settings in which the RKHS $H$ is smooth. A natural question is whether allowing $H$ to be nonsmooth enhances signal recovery properties at stationary points. In the special case where $U$ is restricted to be diagonal, it is shown that the use of a \emph{nonsmooth} Laplace kernel leads to improved signal recovery at stationary points; see~\cite{RuanLiJo25}. \vspace{.5em}
\item (Other RKHS Structures)
An additional open question concerns the role of RKHS induced by inner-product kernels, such as $K(x,x')=\exp(x^\top x')$. Understanding whether such kernels admit analogous geometric properties remains an interesting direction for future work. \vspace{.5em}
\item (Statistical Rates) 
The sharpness property stabilizes the low-rank structure of \emph{finite-sample} minimizers, in a manner analogous to the effect of an explicit nuclear-norm penalty. A natural question is whether this form of structural stability can be leveraged to obtain quantitative generalization error bounds that depend on the intrinsic, rather than ambient, dimension.

\vspace{.5em} 
\item (Statistical Mechanics) 
This paper focuses on the fixed-dimensional regime in which $d$ is fixed and $n \to \infty$; understanding high-dimensional asymptotic behavior where $d/n \to \gamma \in (0, \infty)$ remains an open direction. In the special case where $U$ is restricted to be diagonal,  algorithmic approaches based on iteratively reweighted schemes have been proposed, together with a precise characterizations of the resulting generalization error~\cite{ZhuDaDrFa25}. Extending such approaches to the full matrix setting $U$ or other algorithms could be an interesting direction for future work.

\vspace{.5em}
\item (Low-rank Formulations) 
This paper studies the unconstrained setting $U \in \mathbb{R}^{d \times d}$ and shows that low-rank solution structure emerges in finite samples without explicit rank constraints. A natural question is whether explicitly assuming $U \in \mathbb{R}^{k \times d}$ with $k \le d$ yields statistical or computational benefits. Recent works derived promising statistical rates with reduced dependence on the ambient dimension~\cite{FollainBa24,HuangLaDePoRo25}. It would be interesting to investigate the statistical optimality of these approaches and understand further how to balance statistical adaptivity against computational cost required to achieve it.


\vspace{.5em} 
\item (Descent Algorithms)
An important open question concerns the design of efficient descent algorithms that reliably locate statistically meaningful solutions. 
In the optimization literature, sharpness properties are often leveraged to obtain favorable convergence behavior, e.g.,~\cite{BurkeDe02, BurkeDe05}. A Riemannian optimization algorithm has recently been proposed in~\cite{LiRub25}; however, how such methods (or others) can be connected to, or benefit from, the sharpness properties established in this work remains an open question.

\vspace{.5em} 
\item (Other Algorithms) 
Recent work on the Average Gradient Outer Product (AGOP) method~\cite{Adit24A} demonstrates strong empirical performance across a range of applications, see, e.g.,~\cite{Beagleholeetal25}. However, since AGOP is not derived from an explicit variational principle, its relationship to the underlying variational landscape and the structure of minimizers remains unclear~\cite{RadhakrishnanBeDr25}. Clarifying this connection is a natural and promising direction for future work. Conversely, it may also be of interest to investigate whether insights from the variational formulation in this paper can be leveraged to further improve the numerical performance of AGOP in practice.
\end{itemize}

\vspace{.5em} 
Ultimately, we wish to understand feature learning as a statistical and a variational phenomenon, by treating the compositional kernel models as variational objects governed by the interplay between
the data distribution $\mathbb{P}$, the geometry of the RKHS $H$, sampling effects $n$,
and the matrix structures of $U$ (diagonal, low-rank, or full-rank).
This interplay gives rise to a wide range of mathematically
nontrivial questions, with potential implications for algorithmic design.

\section{Acknowledgements} 
\indent\indent
The authors would like to thank X.Y. Han for enriching discussions on the experiments. 
Feng Ruan thanks Lijun Ding for insightful comments on an early draft of the manuscript; Ying Cui, Jiajin Li, and Yiping Lu for discussions that prompted clarification of future directions and motivated the inclusion of Section~\ref{sec:open-questions}; Dmitriy Drusvyatskiy, Shuo Huang, and Adit Radhakrishnan for clarifying the related work~\cite{RadhakrishnanBeDr25, ZhuDaDrFa25, HuangLaDePoRo25}; and Basil Saeed and Lanqiu Yao for helpful discussions on presentation.

\bibliographystyle{amsalpha} 
\bibliography{sn-bibliography}

\newpage
\appendix 

\section{Deferred Proofs} 
\subsection{Proof of Lemma~\ref{lemma:continuous-embedding}}
\label{sec:proof-of-lemma-continuous-embedding}
By Cauchy-Schwartz, we have 
\begin{equation*}
	\norm{f}_{L_\infty}^2 \le \|\hat{f}\|_{L_1}^2 \le 
		 \int_{\R^d} \frac{|\hat{f}|^2(\omega)}{k(\omega)} d\omega \int_{\R^d} k(\omega) d\omega = \norm{f}_H^2.
\end{equation*} 
Since $\hat{f} \in L_1$, the continuity of $f$ follows from Lebesgue's dominated convergence theorem. 
By the Riemann-Lebesgue lemma, we get $\lim_{x \to \infty} f(x) = 0$.

\subsection{Proof of Lemma~\ref{lemma:denseness-of-H-in-L-2}}
\label{sec:proof-of-lemma-denseness}
Suppose this is not true. There is nonzero $f \in L_2(\nu)$ such that $\int f h d\nu = 0$ 
for all $h \in H$. In particular, taking $h(\cdot) = \kernel(\cdot - z) \in H$, we obtain
$\int f(x) \kernel(x-z)\nu(dx) = 0$ for all $z \in \R^d$.  Using the Fourier inversion formula, 
$\kernel(x-z) =  \int k(\omega) e^{2\pi i \langle x-z, \omega\rangle} d\omega$, we get
\begin{equation*}
\begin{split} 
	\iint f(x)  k(\omega) e^{2\pi i \langle x-z, \omega\rangle} d\omega \nu(dx) = 0~~~\forall z \in \R^d.
\end{split} 
\end{equation*} 
Define the Fourier transform of $f$ under $\nu$ as $F(\omega) := \int f(x) e^{-2\pi i \langle x, \omega \rangle} \nu(dx)$. Then the above equation becomes
\begin{equation*}
	\int F(-\omega) k(\omega) e^{-2\pi i\langle z, \omega\rangle}  d\omega = 0~~~\forall z \in \R^d.
\end{equation*} 
That is, the Fourier transform of $F(-\omega) k (\omega)$ is identically zero. 
Since $k(\omega) > 0$, this implies $F(\omega) = 0$ almost everywhere $\omega$, so $f = 0$ in $L_2(\nu)$—a contradiction.

\subsection{Proof of Lemma~\ref{lem:uniform-gap}}
\label{section:proof-uniform-gap}
For a subspace $T \subset \R^d$, define $\eps_0(T) :=  \E\bigl[(Y - \E[Y | \Pi_T X])^2\bigr] -  \E\bigl[(Y - \E[Y | X])^2\bigr]$. By definition, 
\begin{equation*}
	\eps_0(T) = \E[(\E[Y|\proj_T X]- \E[Y|X])^2] \ge 0
\end{equation*}
To prove Lemma~\ref{lem:uniform-gap}, we need to show: 
\begin{equation}
\label{eqn:general-nonsense}
	\inf \{\eps_0(T) | \dim(T) < \dim(S)\} > 0.
\end{equation}

By the defining minimality property of the central mean subspace $S$, \[
\E[Y | \Pi_T X] \neq \E[Y | X]
\quad \text{with positive probability}
\]
for every subspace $T$ satisfying $\dim(T) < \dim(S)$. Consequently,
$\varepsilon_0(T) > 0$ for all such subspaces.

Assume for contradiction that the infimum in~\eqref{eqn:general-nonsense} equals zero.
Then there exists a sequence of subspaces $\{S_n\}$ with $\dim(S_n) < \dim(S)$
such that $\varepsilon_0(S_n) \to 0$. 
Note the set of orthogonal projection matrices in $\mathbb{R}^{d\times d}$ with rank at most
$k$, i.e., 
$
	\left\{\Pi \in \R^{d\times d}: \Pi = \Pi^\top, \Pi^2 = \Pi, {\rm trace}(\Pi) \le k \right\}
$
is compact for each $k \le d$. Thus, by passing to a subsequence if necessary, we may assume
$\proj_{S_n} \to \proj_{S_\infty}$ holds for some 
subspace $S_\infty$ where $\dim(S_\infty) < \dim(S)$. 

For each $n$, define
\begin{equation*}
	\phi_n(X) :=  \E[Y|\proj_{S_n} X]- \E[Y|X]. 
\end{equation*}
The sequence $\{\phi_n\}$ is uniformly bounded in $\mathcal{L}_2(\mathbb{P})$. By the
Banach--Alaoglu theorem, after possibly passing to a further subsequence,
$\phi_n \rightharpoonup \phi_\infty$ weakly in $\mathcal{L}_2(\mathbb{P})$, meaning that
\[
\lim_{n\to\infty} \E[\phi_n(X) h(X)] = \E[\phi_\infty(X) h(X)]
\quad \text{for all } h \in \mathcal{L}_2(\mathbb{P}).
\]

We claim that \emph{(a)} $\phi_\infty(X) = \E[Y| \Pi_{S_\infty} X] - \E[Y | X]$, and  
\emph{(b)} $\E[\phi_\infty(X)^2] = 0$.

\emph{Proof of (a)}. 
Define $\gamma_n(X) := \E[Y | \Pi_{S_n} X]$. For any continuous and uniformly
bounded function $h : \mathbb{R}^d \to \mathbb{R}$, the defining property of
conditional expectation gives
\[
\E\!\left[(Y - \gamma_n(X))\, h(\Pi_{S_n} X)\right] = 0 .
\]
Let $\gamma_\infty$ be a weak-$*$ accumulation point of $\{\gamma_n\}$ in
$\mathcal{L}_2(\mathbb{P})$. Since $\Pi_{S_n} \to \Pi_{S_\infty}$, we have
$h(\Pi_{S_n} X) \to h(\Pi_{S_\infty} X)$ strongly in $\mathcal{L}_2(\mathbb{P})$
for every such $h$. Because $\{\gamma_n\}$ is uniformly bounded in
$\mathcal{L}_2(\mathbb{P})$ and $\gamma_n \rightharpoonup \gamma_\infty$ weakly,
we may pass to the limit to obtain
\begin{equation}
\label{eqn:orthogonality-condition}
\E\!\left[(Y - \gamma_\infty(X))\, h(\Pi_{S_\infty} X)\right] = 0 .
\end{equation}
Since~\eqref{eqn:orthogonality-condition} holds for all continuous and uniformly
bounded function $h$, it characterizes conditional expectation, implying
$\gamma_\infty(X) = \E[Y | \Pi_{S_\infty} X]$ in $\mathcal{L}_2(\mathbb{P})$.
Therefore,
$\phi_\infty(X) = \E[Y| \Pi_{S_\infty} X] - \E[Y | X]$.

\emph{Proof of (b).}
By assumption,
\[
\E[\phi_n(X)^2] = \varepsilon_0(S_n) \to 0 .
\]
By Cauchy--Schwarz,
$\E[\phi_n(X)\phi_\infty(X)] \to 0$. Since $\phi_n \rightharpoonup \phi_\infty$
weakly in $\mathcal{L}_2(\mathbb{P})$, we get
$
\E[\phi_\infty(X)^2] = 0$. 

%

\emph{Conclusion.}
We have shown that
\[
\varepsilon_0(S_\infty) = \E[\phi_\infty(X)^2] = 0
\quad \text{with} \quad \dim(S_\infty) < \dim(S),
\]
contradicting the defining minimality of the central mean subspace $S$. Hence,
the infimum in~\eqref{eqn:general-nonsense} must be strictly positive, completing
the proof of Lemma~\ref{lem:uniform-gap}.


\subsection{Proof of Lemma~\ref{lemma:basic-convex-geometry}}
\label{sec:proof-of-lemma-basic-convex-geometry}
\indent\indent
Let $\Sigma$ satisfy $\Range(\Sigma) \subset S$.  We wish to establish two sets are equal $T_1 = T_2$ where
\begin{equation*}
\begin{split} 
	T_1&:= \left\{W: W \in T_{\Sym_+^d}(\Sigma) \, , \, \Range(W)\subset S^\perp\right\} \\
	T_2&:= \left\{W: W \in \Sym_+^d \, , \, \Range(W)\subset S^\perp\right\}.
\end{split} 
\end{equation*}

Let $W \in T_1$. By definition, there exists $s > 0$ such that
$\Sigma + s W \in \Sym_+^d$. Hence,
$v^\top (\Sigma + s W) v \ge 0$ for all $v \in \mathbb{R}^d$. 
In particular, if $v \in S^\perp$, then since
$\Range(\Sigma) \subset S$, we have $v^\top \Sigma v = 0$, and therefore
\[
v^\top W v \ge 0, \qquad \forall v \in S^\perp.
\] 
Since $W$ is symmetric and obeys $\Range(W) \subset S^\perp$, it follows that
$Wv = 0$ for all $v \in S$, and therefore
\[
v^\top W v = (P_{S^\perp} v)^\top W (P_{S^\perp} v) \ge 0, \qquad \forall  v \in \mathbb{R}^d.
\]
Thus, $W$ is positive semidefinite, and hence $W \in T_2$.

Conversely, let $W \in T_2$.  Since both $\Sigma$ and $W$ belong to $\Sym_+^d$,
we have $\Sigma + sW \in \Sym_+^d$ for all $s > 0$, which implies $W \in T_{\Sym_+^d}(\Sigma)$ and hence, $W \in T_1$.

\section{Failure of minimizer existence for a discrete distribution}
\label{sec:discrete-counterexample}
\indent\indent
In this section, we illustrate that the continuity assumption on $X$ in
Theorem~\ref{theorem:existence-of-global-minimizers} is necessary by exhibiting
a discrete distribution $\P$ for which the objective $J(\cdot,\lambda)$ does not attain
its infimum.

Let $d=1$ and consider the discrete distribution
\[
\P\big((X,Y)=(+1,+2)\big)
=
\P\big((X,Y)=(-1,-2)\big)
=
1/2.
\]
We consider the Gaussian kernel $K(x,x')=e^{-(x-x')^2}$.
Following the formula~\eqref{eqn:J-formula-finite-atoms}, 
for $\Sigma\ge 0$ we first construct the associated kernel matrix and response vector by
\[
\mathbf{K}_\Sigma
=
\begin{pmatrix}
1 & e^{-4\Sigma} \\
e^{-4\Sigma} & 1
\end{pmatrix},
\qquad
\mathbf{y}=\begin{pmatrix} 2 \\ - 2 \end{pmatrix}.
\]
Substituting these expressions into~\eqref{eqn:J-formula-finite-atoms} yields the explicit formula
\[
J(\Sigma,\lambda)
=
\lambda\,\mathbf{y}^\top(\mathbf{K}_\Sigma+2\lambda I)^{-1}\mathbf{y}
=
\frac{8\lambda}{(1+2\lambda)-e^{-4\Sigma}}.
\]

Since $e^{-4\Sigma}$ is strictly decreasing in $\Sigma$, the function
$\Sigma\mapsto J(\Sigma,\lambda)$ is strictly decreasing and satisfies
\[
\lim_{\Sigma\to\infty} J(\Sigma,\lambda)
=
\frac{8\lambda}{1+2\lambda} = \inf_{\Sigma \ge 0} J(\Sigma,\lambda).
\]
In particular, $J(\cdot,\lambda)$ does not attain its infimum over $\Sym_+^1 = \{\Sigma \in \R: \Sigma \ge 0\}$,
and any minimizing sequence necessarily satisfies $\Sigma\to\infty$.
This example demonstrates that when $\P$ is discrete, minimizers of $J(\cdot,\lambda)$
may escape to infinity, and hence the continuity assumption on $X$ in
Theorem~\ref{theorem:existence-of-global-minimizers} is essential.

\section{Linear RKHS $H$}
\label{sec:linear-vs-nonlinear} 
The conclusions of this paper fail when $H$ is a linear RKHS. In this case,
\begin{equation*}
	H = \bigl\{ x \mapsto \langle w, x \rangle : w \in \R^d \bigr\}.
	\qquad 
	\|f\|_H = \|w\| ,
\end{equation*}
Recall the variational problem
\begin{equation*}
	\inf_U \min_{f \in H} \; \E[(Y - f(UX))^2] + \lambda \|f\|_H^2.
\end{equation*}
In the linear RKHS case, the objective reduces to
\begin{equation*}
	\inf_U \min_{w} \; \E[(Y - \langle w, UX \rangle)^2] + \lambda \|w\|^2 .
\end{equation*}
This formulation is fundamentally ill-posed as a joint optimization problem in $(w,U)$. 
Indeed, for any $c \neq 0$, the reparametrization $(w,U) \mapsto (c w, U/c)$ leaves the prediction
$\langle w, UX \rangle$ unchanged, while scaling the regularization term by $c^2$.
Hence, the objective value can be made arbitrarily close to its infimum,
\begin{equation*}
	\min_{v \in \R^d} \E[(Y - \langle v, X \rangle)^2],
\end{equation*}
without being attained. As a result, the linear RKHS setting 
provides no meaningful notion of minimizers, let alone structural recovery.

\newpage

\section{Additional Results on Synthetic Experiments} 
\label{sec:additional-plots}
In this section, we document additional numerical evidence from the synthetic experiments 
that support the claims (i)--(iv) we made in the main text.

\begin{figure}[!htb]
\centering
  \begin{subfigure}{.8\linewidth}
    \centering
    \includegraphics[width=0.38\linewidth]{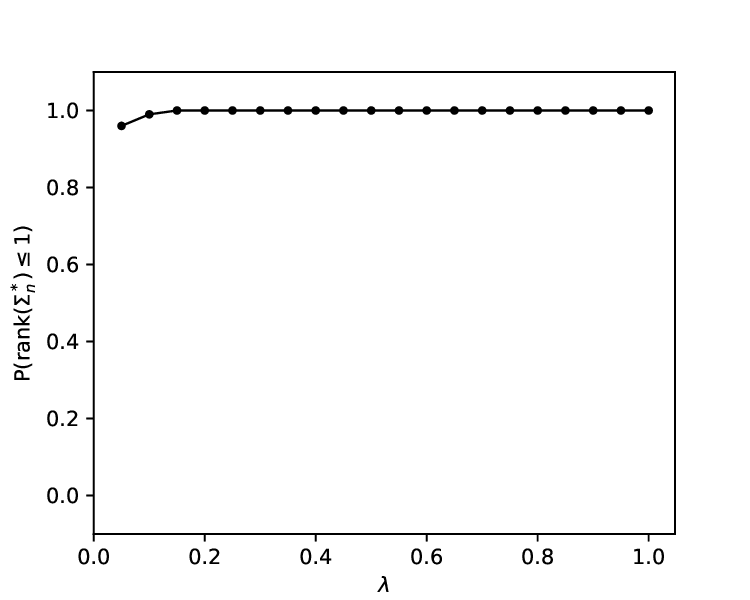}
    \includegraphics[width=0.38\linewidth]{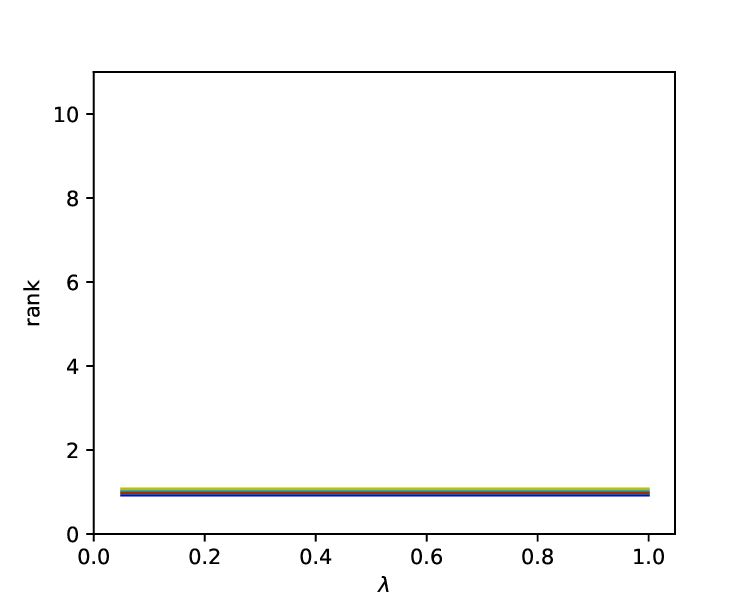}
    \caption{$F(x) = x_1 + x_2 + x_3$.}
  \end{subfigure}
  
  \begin{subfigure}{.8\linewidth}
    \centering
    \includegraphics[width=0.38\linewidth]{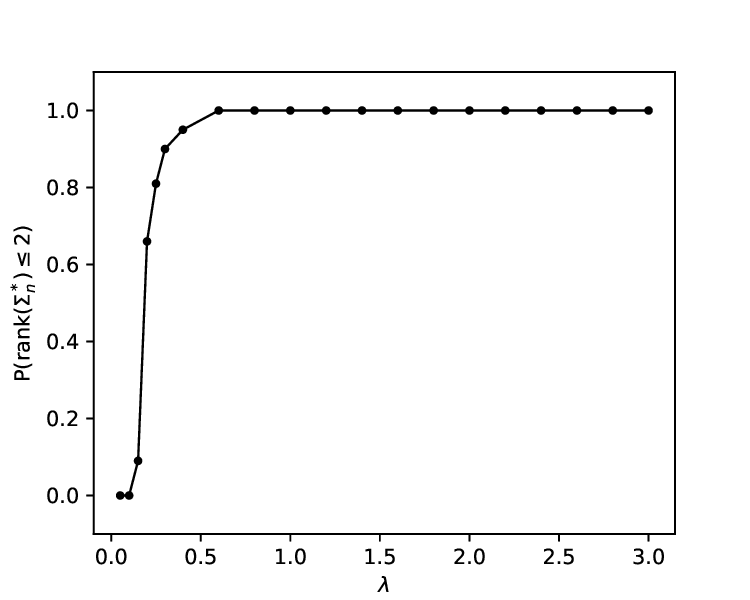}
    \includegraphics[width=0.38\linewidth]{fig/sp_simulation61_Gau_seed100_n300_d50_e0.1_r0.0_c0.0_iter2000_alpha0.001_beta0.5_lr0.1_l0.05_3.0sp5.eps}
    \caption{$F(x) = 0.1(x_1 + x_2 + x_3)^3 + \tanh(x_1 + x_3 + x_5)$.}
  \end{subfigure}
  
  \begin{subfigure}{.8\linewidth}
    \centering
    \includegraphics[width=0.38\linewidth]{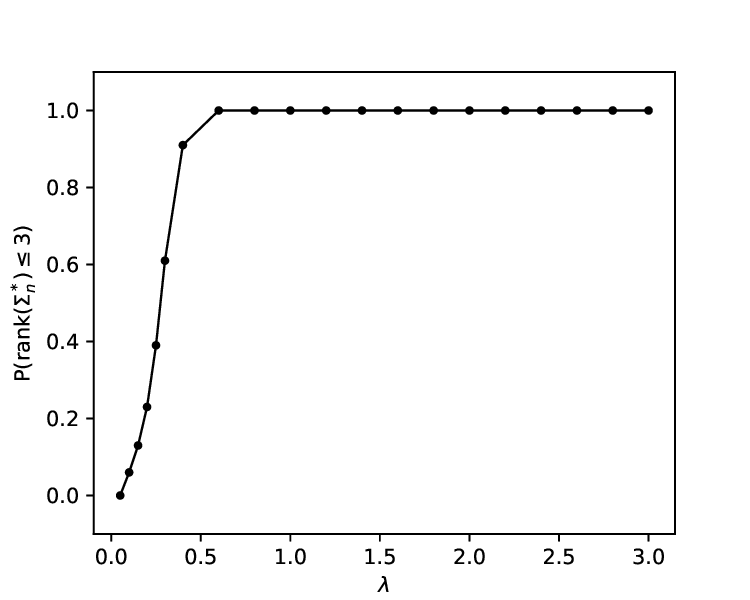}
    \includegraphics[width=0.38\linewidth]{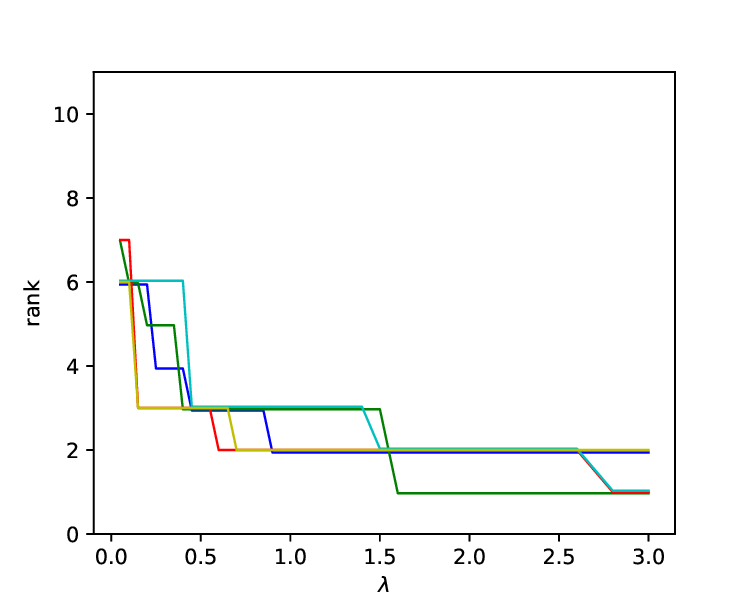}
    \caption{$F(x)= 2(x_1 + x_2) + (x_2 + x_3)^2 + (x_4 - 0.5)^3$.}
  \end{subfigure}

\caption{Plots for Claim (i). 
For each row, the left panel shows the empirical probability of $\rank(\Sigma_{n, *}) \le \dim(S)$ over $100$ repeated
experiments for different $\lambda$ values, while the right panel illustrates the rank of 
$\Sigma_{n, *}$ as $\lambda$ varies, using $5$ example pairs of $(X, y)$.
}
\label{plot:claim1}
\end{figure}

\begin{figure}[!htb]
\centering
  \begin{subfigure}{\linewidth}
    \centering
    \includegraphics[width=0.45\linewidth]{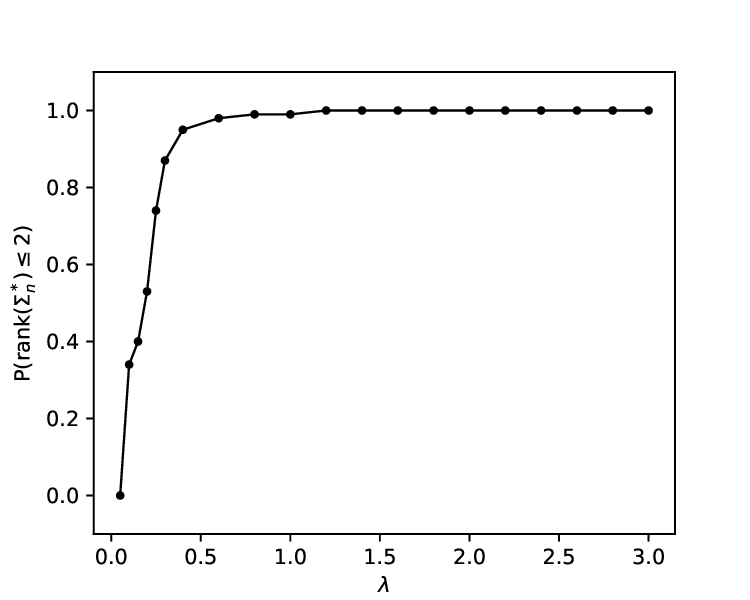}
    \includegraphics[width=0.45\linewidth]{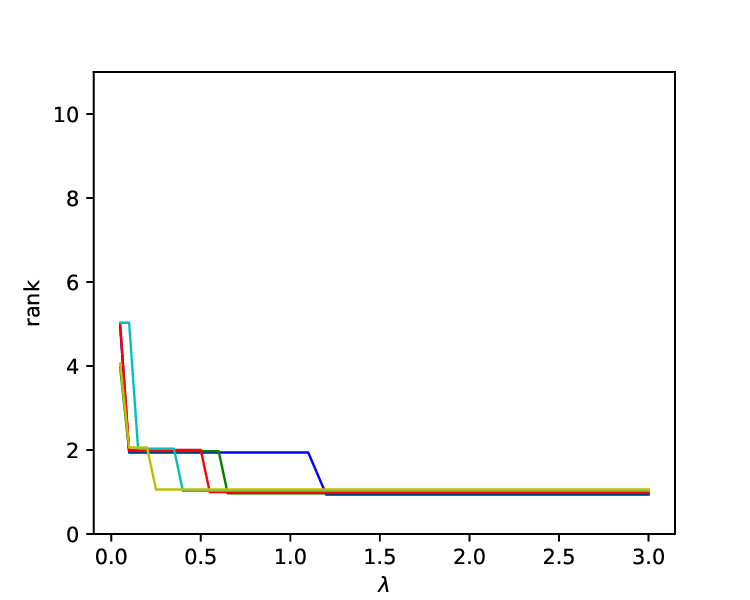}
    \caption{Discrete $X$, $F(x) = 0.1(x_1 + x_2 + x_3)^3 + \tanh(x_1 + x_3 + x_5)$.}
  \end{subfigure}
  

\caption{
Plots for Claim (ii). 
Experimental setting: the covariate $X$ has independent coordinates $X_i$, with each $X_i$ 
obeying $\P(X_i = 1) = \P(X_i = 0) = 0.5$, 
and the response follows $Y = F(X) + \normal(0, \sigma^2)$ where $\sigma = 0.1$. 
Here, we set $n = 300$ and $d = 50$.
The left panel shows the empirical probability of $\rank(\Sigma_{n, *}) \le \dim(S)$ over $100$ repeated
experiments for different $\lambda$ values. The right panel displays how the rank of the solution 
$\Sigma_{n, *}$ changes with different $\lambda$ values, using $5$ example pairs of $(X, y)$.
}
\label{plot:claim3}
\end{figure}

\begin{figure}[!htbp]
\centering  
\includegraphics[width=0.5\linewidth]{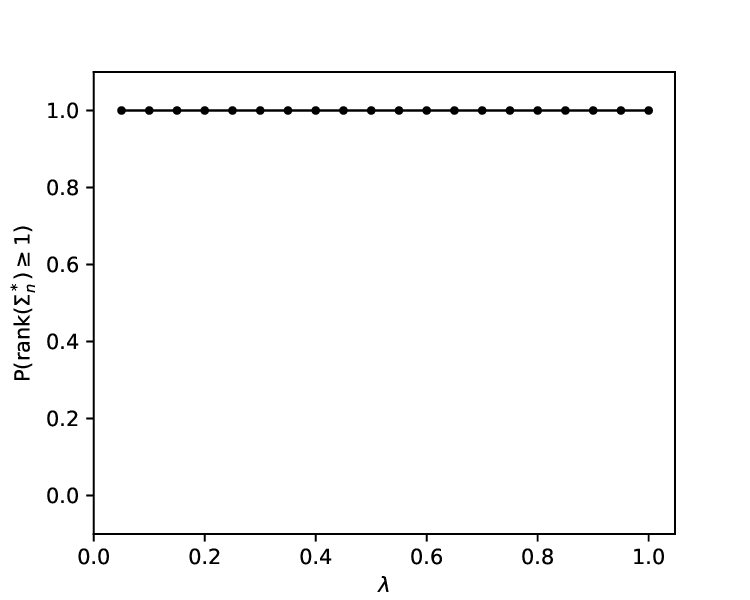}
\caption{Plot for Claim (iii). 
In this experiment, $\E[Y|X] \equiv 0$ where $X\sim \normal(0,I)$, and $Y=\epsilon$, where $\epsilon\sim \normal(0,\sigma^2)$ 
is independent of $X$ with $\sigma=0.1$. Here, we choose $n = 300$ and $d = 50$.  
The solution is always full rank: $\rank(\Sigma_{n, *}) = 50$ for every $\lambda > 0$ in the $100$ repeated experiments. 
}
\label{plot:claim4}
\end{figure}

%

\begin{figure}[!htb]
\centering
  \begin{subfigure}{\linewidth}
    \centering
    \includegraphics[width=0.33\linewidth]{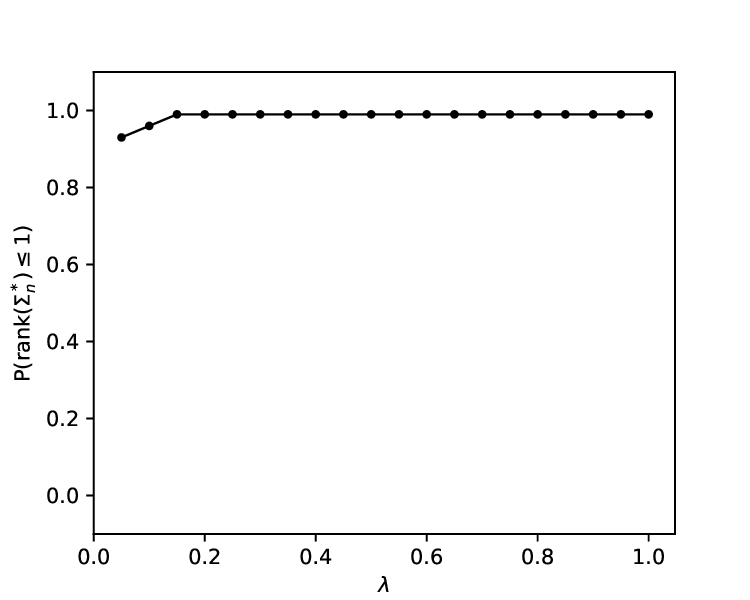}
    \includegraphics[width=0.33\linewidth]{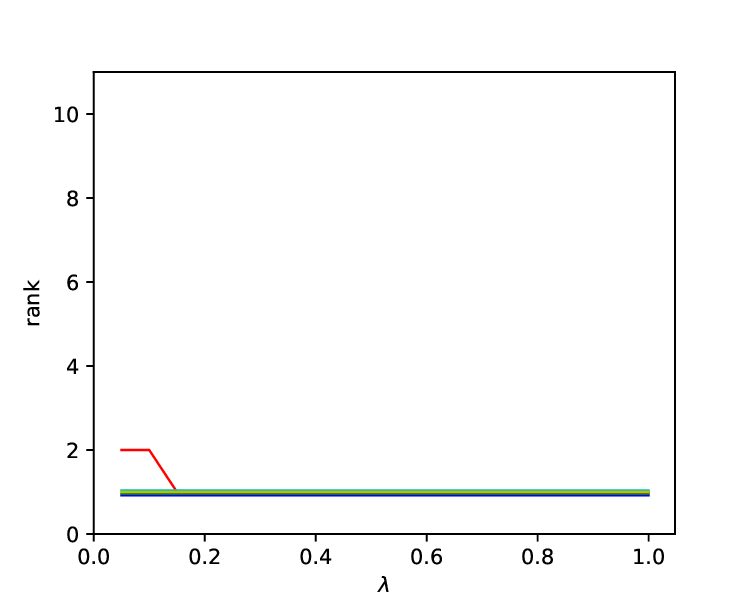}
    \caption{$F(x) = x_1 + x_2 + x_3$.}
  \end{subfigure}

    \begin{subfigure}{\linewidth}
    \centering
    \includegraphics[width=0.33\linewidth]{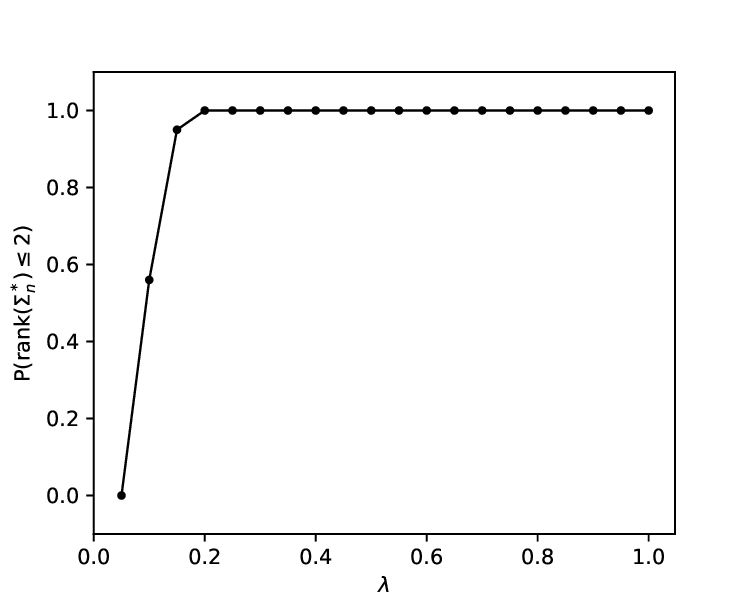}
    \includegraphics[width=0.33\linewidth]{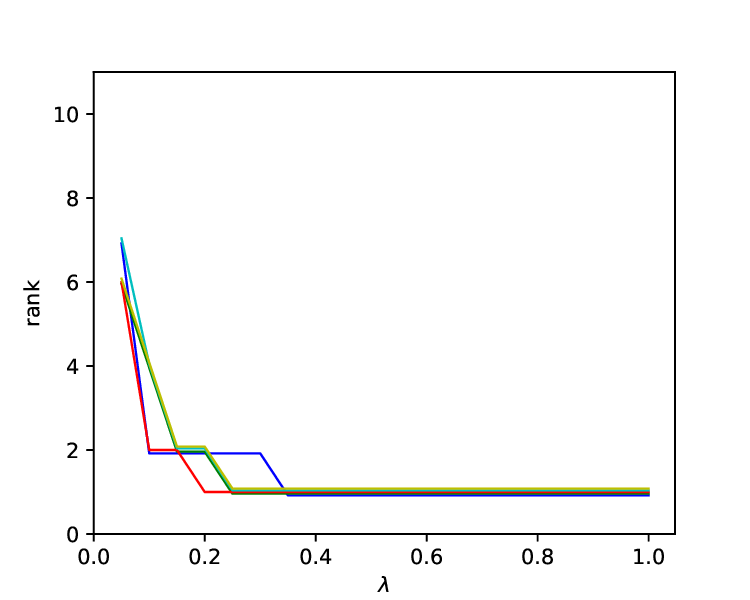}
    \caption{$F(x) = x_1 x_2$.}
  \end{subfigure}

    \begin{subfigure}{\linewidth}
    \centering
    \includegraphics[width=0.33\linewidth]{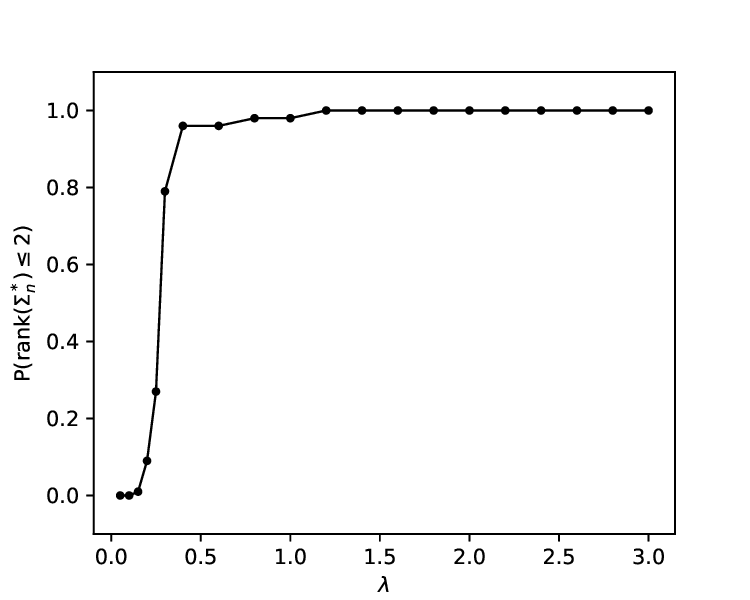}
    \includegraphics[width=0.33\linewidth]{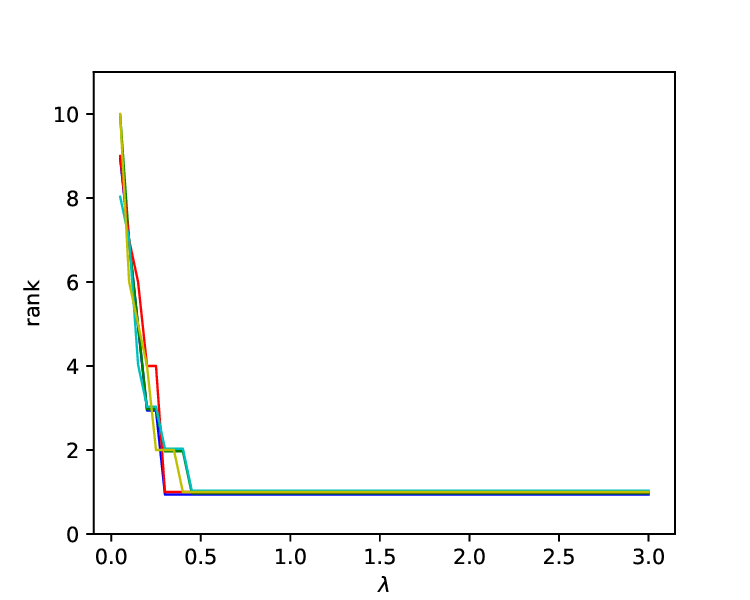}
    \caption{$F(x) = 0.1(x_1 + x_2 + x_3)^3 + \tanh(x_1 + x_3 + x_5)$.}
  \end{subfigure}

    \begin{subfigure}{\linewidth}
    \centering
    \includegraphics[width=0.33\linewidth]{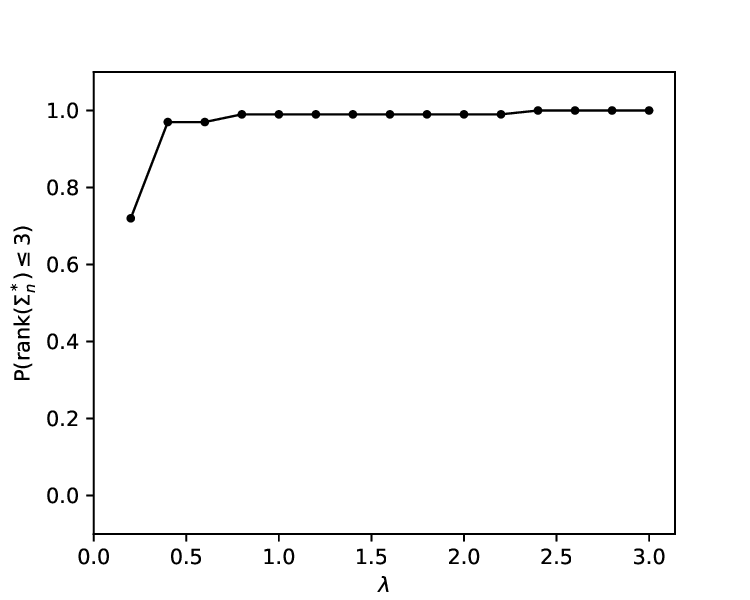}
    \includegraphics[width=0.33\linewidth]{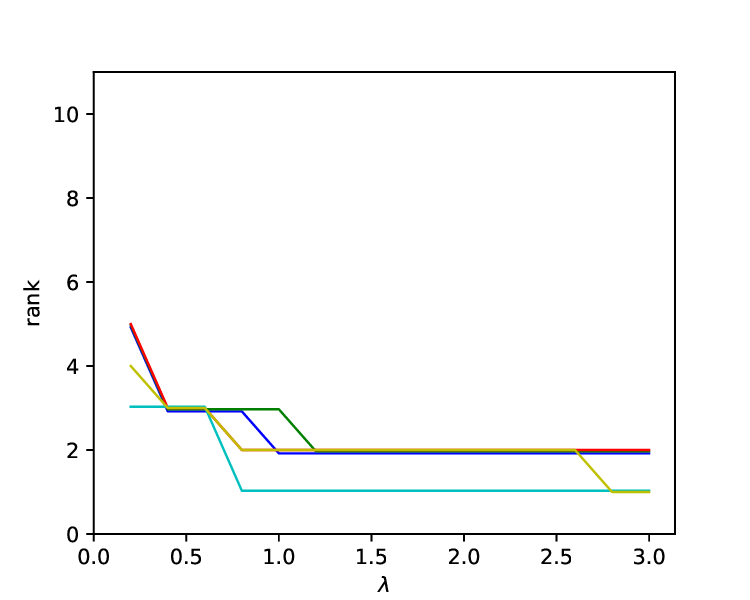}
    \caption{$F(x)= 2(x_1 + x_2) + (x_2 + x_3)^2 + (x_4 - 0.5)^3$.}
  \end{subfigure}
\caption{Plots for Claim (iv). 
Experimental setting: the covariate $X \sim \normal(0, C)$, with $C_{i,j} =  0.5^{|i-j|}$ and
the response follows $Y = F(X) + \normal(0, \sigma^2)$ where $\sigma = 0.1$. Here, we set $n = 300$ and $d = 50$.\\
For each row, the left panel shows the empirical probability of $\rank(\Sigma_{n, *}) \le \dim(S)$ over $100$ repeated
experiments for different $\lambda$ values. The right panel displays how the rank of the solution 
$\Sigma_{n, *}$ changes with different $\lambda$ values, using $5$ example pairs of $(X, y)$.
}
\label{plot:claim2}
\end{figure}

\end{document}